\documentclass{siamart190516}

 \usepackage{amsmath, amsfonts, amsopn}
 \usepackage{amssymb}
 \usepackage{graphicx}
 \usepackage{epstopdf}
 \usepackage{color}
 \usepackage{tcolorbox}
 \usepackage{algorithm}
 \usepackage{algpseudocode}
 \usepackage{url}
 \usepackage[perpage]{footmisc}

\newcommand{\R}{\mathbb{R}}
\newcommand{\C}{\mathbb{C}}
\newcommand{\K}{\mathcal{K}}
\DeclareMathOperator*{\argmin}{arg\,min}

\title{Towards understanding CG and GMRES through examples}

\author{Erin Carson\footnotemark[1], J\"{o}rg Liesen\footnotemark[2] \and
Zden{\v e}k Strako{\v s}\footnotemark[3]}

\begin{document}

\maketitle

\footnotetext[1]{Faculty of Mathematics and Physics, Charles University, Prague,
Czech Republic. Supported by the Exascale Computing Project (17-SC-20-SC), a collaborative effort of the U.S. Department of Energy Office of Science and the National Nuclear Security Administration, and by the European Union (ERC, inEXASCALE, 101075632). Views and opinions expressed are those of the
authors only and do not necessarily reflect those of the European Union or the European Research Council. Neither the European Union nor the granting authority can be held responsible for them.\\
e-mail: {\tt carson@karlin.mff.cuni.cz}.}
\footnotetext[2]{Institute of Mathematics, TU Berlin, Stra{\ss}e des 17. Juni 136, 10623 Berlin, Germany.\\ 
e-mail: {\tt liesen@math.tu-berlin.de}.}
\footnotetext[3]{Faculty of Mathematics and Physics, Charles University, Prague,
Czech Republic. \\
e-mail: {\tt strakos@karlin.mff.cuni.cz}.}

\begin{abstract}
When the conjugate gradient (CG) method for solving linear algebraic systems was formulated about 70 years ago by Lanczos, Hestenes, and Stiefel, it was considered an iterative process possessing a mathematical finite termination property. With the deep insight of the original authors, CG was placed into a very rich mathematical context, including links with Gauss quadrature and continued fractions. The optimality property of CG was described via a normalized weighted polynomial least squares approximation to zero. This highly nonlinear problem explains the adaptation of CG iterates to the given data. Karush and Hayes immediately considered CG in infinite dimensional Hilbert spaces  and investigated its superlinear convergence. Since then, the view of CG, as well as other Krylov subspace methods developed in the meantime, has changed. Today these methods are considered primarily as computational tools, and their behavior is typically characterized using linear upper bounds, or heuristics based on clustering of eigenvalues. Such simplifications limit the mathematical understanding of Krylov subspace methods, and also negatively affect their practical application.

This paper offers a different perspective. Focusing on CG and the generalized minimal residual (GMRES) method, it presents mathematically important as well as practically relevant phenomena that uncover their behavior \emph{through a discussion of computed examples}. These examples provide an easily accessible approach that enables understanding of the methods, while pointers to more detailed analyses in the literature are given. This approach allows readers to choose the level of depth and thoroughness appropriate for their intentions. Some of the points made in this paper illustrate well known facts. Others challenge mainstream views and explain existing misunderstandings. Several points refer to recent results leading to open problems. We consider CG and GMRES crucially important for the mathematical understanding, further development, and practical applications also of other Krylov subspace methods. The paper additionally addresses the motivation of preconditioning.

\end{abstract}

\begin{keywords}
Krylov subspace methods, method of moments, CG method, GMRES method, iterative methods, convergence analysis, rounding error analysis, polynomial approximation problems
\end{keywords}

\begin{AMS}
15A60, 65F10, 65F35
\end{AMS}

\section{Introduction}\label{sec:intro}

Taking the 1952 landmark paper of Hestenes and Stiefel~\cite{HesSti52} on the conjugate gradient method (CG) as their historical starting point, Krylov subspace methods for solving linear algebraic systems $Ax=b$ have been around for more than 70 years. The work of Hestenes and Stiefel, together with the papers from the early 1950s by Lanczos~\cite{Lan50,Lan52,Lan53}, which approach the same topic from a slightly different perspective, crowned the effort of many researchers aiming to find new ways of solving systems of linear algebraic equations using contemporary computers. These works were immediately followed by the description and investigation of the CG method in the infinite dimensional Hilbert space setting, for example by Karush~\cite{Kar52} and Hayes~\cite{Hay54}.\footnote{Krylov subspaces are named after A.~N.~Krylov, who in 1931 published a paper~\cite{Kry31} that described a method for solving a secular equation determining the frequency of small oscillations of mechanical systems. Algebraic formulations of his method were subsequently published by Luzin in 1931~\cite{Luz31} and Gantmacher in 1934~\cite{Gan34}. A detailed account of the life and work of A.~N.~Krylov, and of the early development of Krylov subspace methods can be found, e.g., in~\cite[Section~2.5.7]{LieStr13}.} 

Over the last 70 years, tens of thousands of research articles on the derivation, analysis, or applications of Krylov subspace methods have been published by authors coming from the most diverse scientific backgrounds. Numerous different Krylov subspace methods have been developed. Their names typically give very brief descriptions of their algorithmic properties, and the methods then are referred to by acronyms that abbreviate these names. Examples include, besides CG, the biconjugate gradient (BiCG) method~\cite{Fle76} (with the idea of biconjugation going back to Lanczos~\cite{Lan50,Lan52}), the stabilized biconjugate gradient (BiCG-Stab) method~\cite{Van92}, the full orthogonalization method (FOM)~\cite{Saa81}, the minimal residual (MINRES) method~\cite{PaiSau75}, the generalized minimal residual (GMRES) method~\cite{SaaSch86}, the quasi-minimal residual (QMR) method~\cite{FreNac91}, and related methods like the induced dimension reduction (IDR) method~\cite{WesSon80}. The majority of publications on Krylov subspace methods focus on algorithmic techniques, i.e., on construction of methods or preconditioners. A minority focus on mathematical principles and analysis of the methods in a wider context. Comprehensive treatments of the area can be found in many monographs and survey articles that are devoted entirely, or at least to a large extent, to Krylov subspace methods; see,  e.g.,~\cite{GreBook97,LieStr13,MeuBook06,MeuDui20,Saa03,VanBook03} and~\cite{AshManSay90,FreGolNac92,GolVdV97,LieTic04,MeuStr06}, respectively. 

Instead of adding yet another technical and algorithmic overview or taxonomy of Krylov subspace methods, this paper approaches the methods through investigation of important and practically relevant phenomena that uncover their mathematical foundations. This leads to understanding their behavior, and it also allows the clarification of persisting misunderstandings and issues that still remain open.  We use \emph{computed examples} for this purpose. Many of these examples can be found scattered throughout the existing literature. But their presentation and organization in this paper represents an entirely different approach than in the monographs and survey papers mentioned above. We have organized everything around what we consider \emph{the main points for the understanding of Krylov subspace methods}. 

We believe that such an approach and the perspective offered by this paper can help students, researchers, and practitioners to gain additional insights into Krylov subspace methods, and that in this way insight can be gained more easily than through extensive technical expositions. The text gives many references to the literature containing more detailed descriptions and analysis. We have selected works that underline the presented points and provide technical and historical background for the given arguments. Readers can thus get involved in the subject in the way and to the depth they find appropriate.

We focus on the CG method and the GMRES method, which have evolved as the standard iterative methods for solving linear algebraic systems with symmetric positive definite and general (nonsymmetric) matrices, respectively. Their understanding is a prerequisite for understanding other Krylov subspace methods, as well as for recent and possible future developments that involve them.  In each of the computed examples, we first describe the setup as transparently as possible. We then describe the phenomena observable in the computed figures, followed by an explanation that usually contains pointers to further research literature. Throughout the paper we purposely use simple data (matrices and right-hand sides), so that the examples can be easily reproduced. Clearly, difficulties that can be observed on small model problems usually do not disappear in problems coming from real-world applications. On the other hand, numerical efficiency demonstrated on small model problems may not materialize in large-scale computations. We have this in mind and aspire not to introduce distorted arguments.

Before starting with CG and GMRES, we would like to present a point that is often neglected and that is most important for the understanding of Krylov subspace methods. We consider a linear algebraic system $Ax=b$ and an initial vector $x_0$. Krylov subspace methods are based on projecting onto (some variant of) the Krylov subspace $\K_k(A,r_0)={\rm span}\{r_0,Ar_0,\dots,A^{k-1}r_0\}$, $k=1,2,\dots$, where $r_0=b-Ax_0$; see, e.g., \cite{Saa03}. Obviously,  they are nonlinear in $A$ because Krylov subspaces are formed using repeated multiplications with the matrix $A$. But this form of nonlinearity can be attributed to any iteration process repeatedly applying an iteration matrix. In particular, it can obviously be attributed to the Chebyshev method; see, e.g., \cite{GolVar61a,GolVar61b}. A principally different and
more substantial nonlinearity with respect to $A$ as well as to $r_0$ comes from using {\em projections onto Krylov subspaces}, which is equivalent to enforcing some form of {\em optimality} of the approximate solution {\em that requires  adaptation to the data} $A$ and $r_0$ at each iteration step; see, e.g.,~\cite{LieStr13} for a comprehensive treatment. This essential point is frequently overlooked with far reaching consequences. In some cases the nonlinear behavior of a Krylov subspace method is simply, and mathematically incorrectly, identified with widely known linear convergence bounds.{\footnote{For symmetric positive definite matrices this means that CG is considered equal in performance to the Chebyshev method with a priori knowledge of the edges of the spectrum. Sometimes the distinction between CG and the Chebyshev method is indeed reduced only to the fact that CG does not need any such a priori information. But the distinction consists of their principally different optimality properties, which is not linked in the Chebyshev method to an adaptation to the data; see  Section~\ref{sec:CG} below. This is clear from the early papers by Hesteness, Stiefel, and Lanczos (see, in particular, \cite{Lan53,Lan52}) mentioned above, although the description of the Chebyshev method came later.}}
In other cases, the nonlinear behavior is accepted but viewed primarily as an obstacle for the analysis. However, we claim the following:

\medskip
\begin{tcolorbox}
\emph{Main point: The nontrivial nonlinearity is the main mathematical asset as well as the beauty of Krylov subspace methods, since it requires the methods to adapt to the hidden inner structure of the problem to be solved. This can lead to a significant speedup of the convergence in comparison with (linear) iterative methods that do not adapt to the problem.}
\end{tcolorbox}

\medskip
Because of their nonlinearity, which is rooted in their optimality requiring adaptation to the data, Krylov subspace methods can show their advantages particularly in solving hard practical problems, where their behavior can not be understood using linear upper bounds. Neglecting the nonlinearity hampers further investigation of intricate and open problems, which however is critically needed for advancing the theory as well as the successful practical application of the methods. The main point stated above will therefore reappear in many of the computed examples below.

The paper is organized as follows. In Section~\ref{sec:CG} we consider the CG method, and in Section~\ref{sec:GMRES} the GMRES method. Both sections start with a brief description of the methods (mathematical properties and standard implementations), followed by the computed examples. Section~\ref{sec:Conclusions} contains concluding remarks. The Appendix contains comments on the modern relevance of early works on Krylov subspace methods, and on the topic of preconditioning.

\smallskip
\paragraph{Notation and conventions} Throughout the paper we consider real linear algebraic systems for simplicity of notation. Most results can be easily extended to the complex case. We use $N$ for the matrix size, i.e., $A\in\R^{N\times N}$, and~$k$ denotes the iteration number. Usually the right-hand side~$b$ is a normalized vector of ones, and the initial approximate solution is $x_0=0$. We use the term ``mathematical'' to refer to cases where computations are performed exactly, i.e., in infinite precision, and the term ``computational'' to refer to finite precision computations. In some experiments we compare the infinite and finite precision behavior of algorithms. Finite precision computations are performed using the standard double precision arithmetic in MATLAB. Unless otherwise specified, the mathematical (infinite precision) behavior is then simulated using the Advanpix Multiprecision Computing Toolbox for MATLAB\footnote{\url{https://www.advanpix.com}}.

\section{The CG Method}\label{sec:CG}

The CG method is well defined for any linear algebraic system $Ax=b$ with a symmetric positive definite matrix $A\in\R^{N\times N}$ and right-hand side $b\in\R^N$. If $x_0\in\R^N$ is an initial approximation, and $d=d(A,r_0)$ is the grade\footnote{The grade of $r_0$ with respect to $A$ is defined as the maximal dimension of the Krylov space generated by $A$ and $r_0$.} of the initial residual $r_0=b-Ax_0$ with respect to $A$, then at every step $k=1,2,\dots,d$ the CG method constructs a uniquely determined approximation
\begin{equation}\label{eqn:CGorth}
x_k\in x_0+\K_k(A,r_0)\quad\mbox{such that}\quad r_k\perp \K_k(A,r_0),
\end{equation}
where $\K_k(A,r_0):={\rm span}\{r_0,Ar_0,\dots,A^{k-1}r_0\}$ is the $k$th Krylov subspace generated by $A$ and $r_0$. Mathematically the method terminates with $x_d=x$.

There are many equivalent formulations of the task that is solved by the CG method. For example, at step $k$ the CG method determines the solution of the simplified Stieltjes moment problem (see~\cite{PozStr19}) or, equivalently, it determines the $k$-point Gauss quadrature of the Riemann-Stieltjes integral defined by $A$ and $r_0$; see, e.g.,~\cite[Section~3.5]{LieStr13} and~\cite[Section~5.2]{MSB15} for overviews.
These connections were pointed out by Hestenes and Stiefel in~\cite{HesSti52} and they are important for understanding the mathematical as well as the computational behavior of CG.

There are also many mathematically equivalent algorithms that realize the projection process \eqref{eqn:CGorth}. The most popular variant is the original formulation of Hestenes and Stiefel~\cite{HesSti52}, shown in Algorithm \ref{alg:hscg}. This algorithm recursively updates coupled 2-term recurrences for the approximate solution $x_{k+1}$ and residual $r_{k+1}$, as well as the auxiliary ``search direction'' vector $p_{k+1}$. As it turns out, this variant is also preferable computationally; see, e.g., \cite{Rei71, GutStr00}.

\begin{algorithm}[htbp!]
	\caption{Conjugate Gradient (2-term recurrence variant) \label{alg:hscg}}
	\begin{algorithmic}[1]
		\Require{Symmetric positive definite matrix $A\in\R^{N\times N}$; right-hand side $b$; initial approximation $x_0$; convergence tolerance $\tau$; maximum number of iterations $n_{\max}$.}
		\State{$r_0 = b-Ax_0$}
		\State{$p_0 = r_0$}
		\For{$k = 0, 1, 2,\ldots,n_{\max} $}
			\State{$\alpha_k = (r_k^T r_k)/(p_k^T Ap_k)$}
			\State{$x_{k+1} = x_k + \alpha_k p_k$}
			\State{$r_{k+1} = r_k - \alpha_k Ap_k$}
			\State{Test for convergence using tolerance $\tau$. If satisfied, then return $x_{k+1}$ and stop.}
			\State{$\beta_{k+1} = (r_{k+1}^T r_{k+1})/(r_k^T r_k)$}
			\State{$p_{k+1} = r_{k+1} + \beta_{k+1}p_k$}
		\EndFor
	\end{algorithmic}
\end{algorithm}

Let $A=Q\Lambda Q^T$, with $\Lambda={\rm diag}(\lambda_1,\dots,\lambda_N)$ and $0<\lambda_1\leq\cdots\leq \lambda_N$, be an orthogonal diagonalization of $A$.\footnote{CG is mathematically invariant under orthogonal transformations of the basis in~$\R^N$. In particular, one can study its mathematical behavior using the basis formed by the orthonormalized eigenvectors of $A$, i.e., using the diagonal matrix $\Lambda$ instead of $A$. The results of Greenbaum~\cite{Gre89} allow to view finite precision CG computations (apart from a small inaccuracy) as exact CG for a particular larger matrix having all its eigenvalues close to the eigenvalues of $A$. Combining this with experimental demonstrations in~\cite{GreStr92}, and with further theoretical arguments in~\cite[Section~5.9.1]{LieStr13} and~\cite{MeuStr06}, one can study also the computational behavior of CG using the diagonal matrix $\Lambda$.} We can represent the initial residual $r_0$ by its components in the individual eigenvectors of $A$, stored in the columns of $Q$, as $r_0=Q[\eta_1,\dots,\eta_N]^T$. The approximation $x_k\in x_0+\K_k(A,r_0)$ that is uniquely determined by the orthogonality condition in \eqref{eqn:CGorth} satisfies the (equivalent) optimality property
\begin{align}
\|x-x_k\|_A &= \min_{p\in P_k(0)}\|p(A)(x-x_0)\|_A=\min_{p\in P_k(0)}\left(\sum_{i=1}^N \eta_i^2
\frac{p(\lambda_i)^2}{\lambda_i}\right)^{1/2}, \label{eqn:CGmin}
\end{align}
where $P_k(0)$ denotes the set of polynomials of degree at most $k$ with value~$1$ at the origin; see, e.g.,~\cite[Section~5.6]{LieStr13}. Thus, in every step the CG method solves a certain weighted polynomial approximation problem on the discrete set $\{\lambda_1,\dots,\lambda_N\}$. Moreover, if $\theta_1^{(k)}, \dots , \theta_k^{(k)}$ are the $k$ roots of the polynomial providing the minimum in~(\ref{eqn:CGmin}), then we can easily get
\begin{equation}
\|x - x_k\|_A^2 =  \sum_{i=1}^N\prod_{\ell=1}^k \left(1-\frac{\lambda_i}{\theta_\ell^{(k)}}\right)^2 \frac{\eta_i^2}{\lambda_i},
\label{eqn:CG_Ritz_values}
\end{equation}
which establishes the relationship of the roots of the minimizing polynomial in~(\ref{eqn:CGmin}), called also the Ritz values, with the eigenvalues of the matrix $A$.

Note that
$$\sum_{i=1}^N \frac{\eta_i^2}{\lambda_i} = r_0^T A^{-1}r_0 =(x-x_0)^T A (x-x_0)=\|x-x_0\|_A^2.$$
Therefore, maximizing over the values $p(\lambda_i)$ in the minimization problem on the right-hand side of \eqref{eqn:CGmin} and dividing by $\|x-x_0\|_A$ gives the upper bound
\begin{equation}
\frac{\|x-x_k\|_A}{\|x-x_0\|_A}\leq
\min_{p\in P_k(0)}\max_{1\leq i\leq N} |p(\lambda_i)|
\label{eqn:polynomial_Worst_case_bound}
\end{equation}
It is important to note that the polynomial min-max approximation problem on the right-hand side of \eqref{eqn:polynomial_Worst_case_bound} only depends on $A$, but not on $r_0$.

If $d(A)$ denotes the degree of the minimal polynomial of $A$, then $d(A)\geq d$. It was shown by Greenbaum~\cite{Gre79}, that for any given symmetric positive definite matrix~$A\in\R^{N\times N}$ the bound \eqref{eqn:polynomial_Worst_case_bound} is sharp in the sense that for every step $k\leq d(A)$ there exists an initial residual $r_0$ so that equality holds. Thus, for the given matrix $A$ the value of the polynomial min-max approximation problem is an attainable \emph{worst-case bound} on the relative $A$-norm of the error in the CG method at every step $k \leq d(A)$. (The step $k=d(A)$ is trivial.)
Moreover, for every $k=1,\dots,d(A)-1$, there exist $k+1$ distinct eigenvalues
$\widehat{\lambda}_1,\dots,\widehat{\lambda}_{k+1}$ of $A$, such that
\begin{equation}\label{eqn:CG_eigbd2}
\min_{p\in P_k(0)}\max_{1\leq j\leq N}\,|p(\lambda_j)|=
\left(\sum_{i=1}^{k+1}\prod_{\substack{j=1\\ j\neq i}}^{k+1}
\frac{\widehat{\lambda}_j}{|\widehat{\lambda}_j-\widehat{\lambda}_i|}\right)^{-1}.
\end{equation}
The value of the worst-case bound for CG in step $k$ is thus expressed in terms of a subset of $k+1$ particular eigenvalues of~$A$. This subset is determined by the min-max polynomial approximation problem on the set $\{{\lambda}_1,\dots,{\lambda}_{N}\}$ for the polynomials of degree $k$ that are normalized at the origin.\footnote{The whole paper~\cite{Gre79} is worth reading. It illustrates the peculiarities of the polynomial approximation problem on a discrete set of points, and it derives the weights for which the normalized $k$th degree min-max polynomial on a subset $\{\widehat{\lambda}_1,\dots,\widehat{\lambda}_{k+1}\}$ of the set $\{\lambda_1,\dots,\lambda_N\}$ is equal to the normalized weighted least squares polynomial approximation to zero on $\{\widehat{\lambda}_1,\dots,\widehat{\lambda}_{k+1}\}$.}

Replacing the discrete set $\{\lambda_1,\dots,\lambda_N\}$ by the continuous interval $[\lambda_1,\lambda_N]$ and using Chebyshev polynomials on this interval yields (with a small additional simplification)
\begin{align}\label{eqn:rhoAbound}
\min_{p\in P_k(0)}\max_{1\leq i\leq N} |p(\lambda_i)|\leq
2\left(\frac{\sqrt{\kappa(A)}-1}{\sqrt{\kappa(A)}+1}\right)^k,\quad
\kappa(A)=\frac{\lambda_N}{\lambda_1}.
\end{align}
This bound represents a substantial simplification. While (\ref{eqn:CG_eigbd2}) gives a sharp upper bound on the relative CG error $A$-norms for the given (fixed) spectrum and any initial residual, the right hand side in (\ref{eqn:rhoAbound}) gives a nearly sharp upper bound for {\em arbitrary} eigenvalues $\{\lambda_2,\dots,\lambda_{N-1}\}$ in the interval $[\lambda_1, \lambda_N]$ and any initial residual.

Combining (\ref{eqn:rhoAbound}) with \eqref{eqn:CGmin} results in the frequently stated convergence bound
\begin{align}\label{eqn:rhoAbound1}
\frac{\|x-x_k\|_A}{\|x-x_0\|_A}\leq 2\left(\frac{\sqrt{\kappa(A)}-1}{\sqrt{\kappa(A)}+1}\right)^k.
\end{align}
We will sometimes refer to \eqref{eqn:rhoAbound1} as the \emph{$\kappa(A)$-bound}. This bound implies that if the condition number $\kappa(A)$ is  small, then a fast reduction of the $A$-norm of the error in the CG method can be expected. This bound does \emph{not} imply, however, that a large condition number results in slow convergence of CG. In particular, preconditioners that provide smaller condition numbers than others do not necessarily lead to faster convergence. A convincing example in~\cite{GerMarNieStr19} uses standard PDE test problem and contributes towards opening a new line of research combining the PDE operator context with the algebraic matrix context arising from discretization; see further comments on this point below.

Also note that the $\kappa(A)$-bound for CG is a \emph{linear bound for a nonlinear process}. A comparison with the value of the polynomial min-max approximation problem in~(\ref{eqn:CG_eigbd2}), which gives the worst-case CG value in step $k$ for the given spectrum of the matrix~$A$, shows that neglecting the eigenvalue distribution of $A$ in the interval $[\lambda_1, \lambda_N]$ can mean a substantial loss of information. Similarly, a comparison with the actual minimization problem~(\ref{eqn:CGmin}), which is solved by CG applied to the linear system $Ax=b$ with the initial approximation $x_0$, shows that the size of the components $\eta_j$ of $r_0$ in the invariant subspaces of $A$ can be important; see also~(\ref{eqn:CG_Ritz_values}).

In the examples that follow, we will frequently make use of a certain class of diagonal matrices which is often used in the literature to illuminate the behavior of CG; see, e.g., \cite{GreStr92}. For given $N\geq 3$, $0<\lambda_1<\lambda_N$, and $\rho>0$ we define
\begin{equation}
A={\rm diag}(\lambda_1,\lambda_2,\dots,\lambda_{N-1},\lambda_N)\;\;\mbox{with}\;\;\lambda_i = \lambda_1 + \left( \frac{i-1}{N-1}\right) (\lambda_N -\lambda_1)\rho^{N-i},
\label{eq:diagmatrix}
\end{equation}
for $i=2,\dots,N-1$.
The parameter $\rho$ determines the eigenvalue distribution of $A$. When $\rho=1$, the eigenvalues are equally spaced between $\lambda_1$ and $\lambda_N$. As $\rho$ becomes smaller, the eigenvalues accumulate towards $\lambda_1$. As mentioned above, we can use a diagonal matrix mathematically and computationally without any loss of generality.

\subsection{Mathematical behavior of CG for different eigenvalue distributions}\label{sec:cgeigdist}
\phantom{}
\smallskip
\begin{tcolorbox}
\emph{Main point: The CG optimality property (see \eqref{eqn:CGmin}) depends on the positions of the individual eigenvalues. Therefore CG adapts without any a priori information not only to the spectral interval, but in a significant (and nonlinear) way also to the distribution of the inner eigenvalues. Acceleration of CG convergence is more pronounced for matrices with outlying eigenvalues, and is different when the outliers are small or large.} 
\end{tcolorbox}

\emph{Setup:} We consider the behavior of CG in exact arithmetic for matrices having three different eigenvalue distributions. All matrices are diagonal with $N=30$, $\lambda_1=0.1$, and $\lambda_N=10^3$. The first matrix is a slight modification of \eqref{eq:diagmatrix},  with $\lambda_i = \lambda_N - \frac{(i-1)}{(N-1)}(\lambda_N - \lambda_1)\rho^{N-i}$ for $i=2,\ldots,N-1$ and $\rho=0.6$, so that the eigenvalues accumulate on the right side of the spectrum. The second matrix is \eqref{eq:diagmatrix} with $\rho=0.6$, so that its eigenvalues accumulate to the left side of the spectrum, and the third matrix is \eqref{eq:diagmatrix} with $\rho=1$, so that its eigenvalues are equally spaced. In all cases we use $b=[1,\dots,1]^T/\sqrt{N}$, and $x_0=0$.

\begin{figure}[t!]
\includegraphics[width=0.43\textwidth]{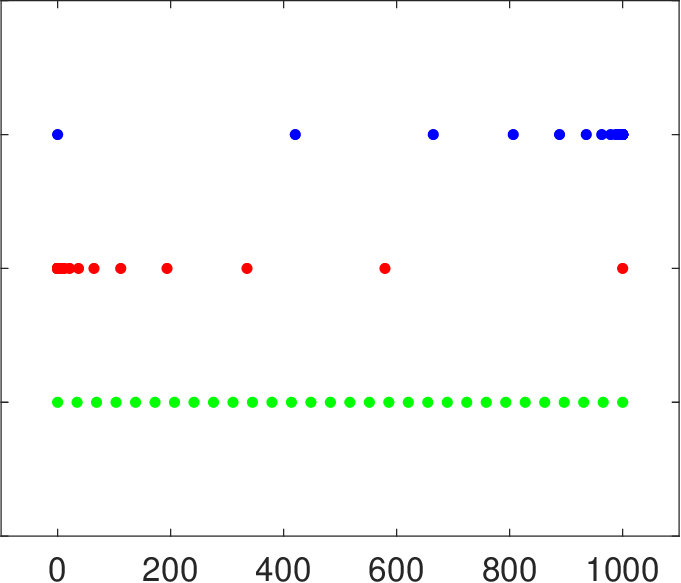}\hfill
\includegraphics[width=0.475\textwidth]{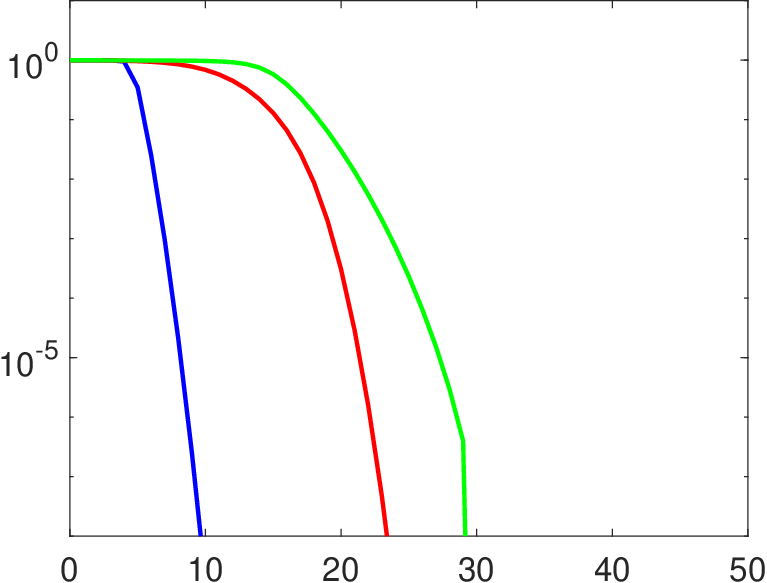}
\caption{Left: Three distributions of 30 eigenvalues in $[0.1,\,10^3]$.
Right: The relative error in the $A$-norm for exact CG applied to the corresponding linear algebraic systems.}\label{fig:cgeigdist}
\end{figure}

\begin{figure}[t!]
\centering
\includegraphics[trim={4cm 3cm 4cm 3cm},clip,width=0.575\textwidth]{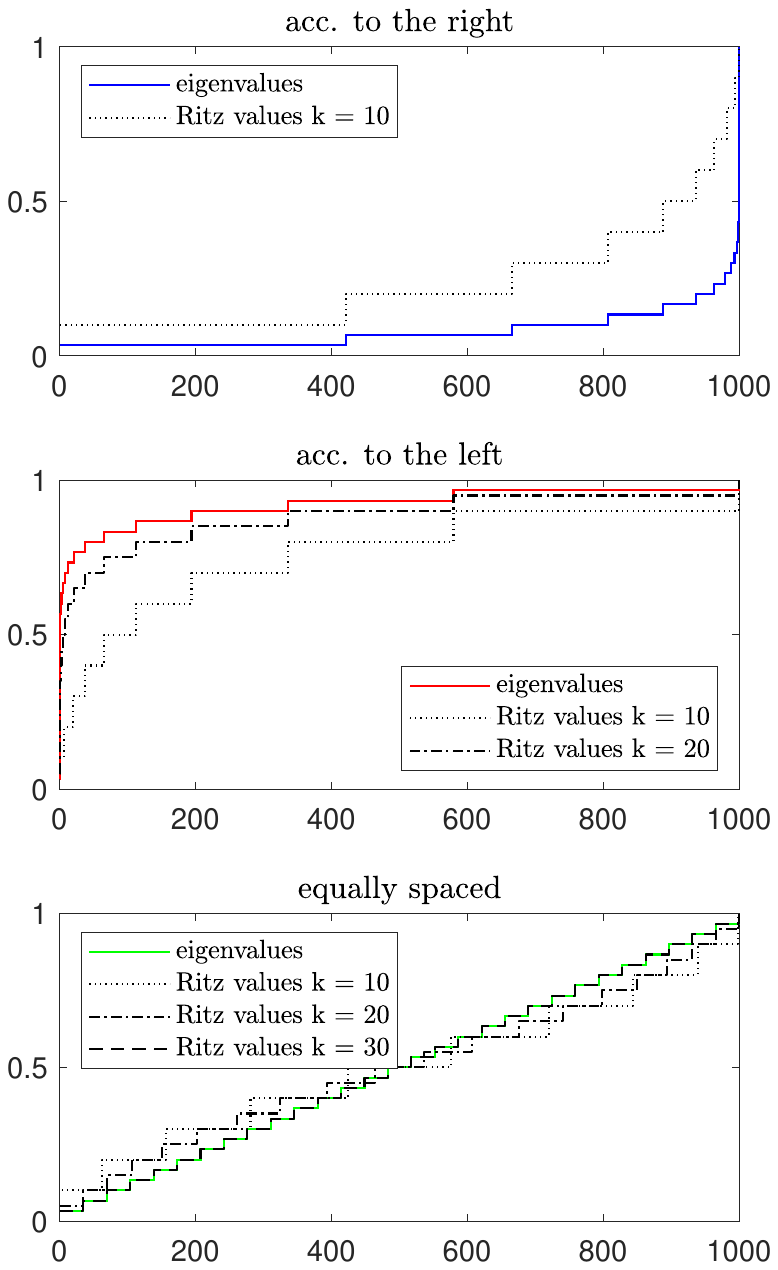}
\caption{Cumulative spectral density plots for Figure \ref{fig:cgeigdist}.}\label{fig:specdensexact}
\end{figure}

\smallskip
\emph{Observations:} The eigenvalue distributions are shown in the left part of Figure~\ref{fig:cgeigdist}. The right part of Figure~\ref{fig:cgeigdist} shows the mathematical behavior of CG. For the matrix with eigenvalues accumulated to the right (blue), CG  converges fastest. For the matrix with eigenvalues accumulated to the left (red), CG converges significantly slower. For the matrix with equally spaced eigenvalues (green), CG converges the slowest.  In Figure \ref{fig:specdensexact} we show cumulative spectral density (CSD) plots using the stepwise functions with points of increase at Ritz values and the size of the vertical steps equal for each Ritz value (see \cite[Appendix C]{LinSaaYan16})\footnote{In order to illustrate the position of Ritz values (that can form tight clusters), we approximate the CSD defined by the spectrum of the matrix $A$ using the CSD associated with the Ritz values. We purposefully do not use an approximation of the CSD defined by the spectrum of $A$ via the Riemann-Stieltjes distribution functions associated with the Gauss quadrature (see, e.g., \cite[Section 3.5]{LieStr13} and \cite[Section 3.2]{MSB15}), since this would not allow to observe forming clusters of Ritz values. This will be important in the next two subsections.}.

\smallskip
\emph{Explanation:}
As can be seen from \eqref{eqn:CGmin}, for equal components of the initial residual in the invariant subspaces, the $A$-norm of the error minimizes the sum of the squared values of the CG polynomial divided by the associated eigenvalues. For an accumulation of the eigenvalues to the right, the CG polynomial approximates by its roots (Ritz values) the small outlying eigenvalues within a few iterations, which takes care for the associated part of the sum. For the rest of the eigenvalues that are large and close to each other, the values of the polynomial do not need to be be so small, because their squares are divided by the large eigenvalues. Therefore fast convergence (as if the small outlying eigenvalues were nonexistent) will occur within a few iterations; see \cite[Theorem 5.6.9]{LieStr13} and the enlightening paper by van der Sluis and van der Vorst~\cite{SluVor86}.\footnote{Although the paper~\cite{SluVor86} assumes exact arithmetic and focuses on spectra with small outlying eigenvalues, it also comments on the case of large outlying eigenvalues, and on the difference between the mathematical and computational CG behavior in such cases; see \cite[p.~559]{SluVor86} and the detailed discussion with further references in \cite[pp.~279-280]{LieStr13}.} Section 5.6.4 of \cite{LieStr13}, called ``Outlying Eigenvalues and Superlinear Convergence'', recalls further closely related results by Lanczos, Rutishauser, Jennings, and others. The arguments above also explain why in this case the convergence rate becomes fast even when the CSD is not yet closely approximated.

For the eigenvalues accumulated to the left, the large outliers are also well approximated by the Ritz values within a few iterations. However, since for the bulk of the small eigenvalues the CG polynomial must place many roots close to the left end of the spectrum in order to make up for the division of its squared values by the small eigenvalues, the acceleration of convergence appears much later. Therefore also the CSD must be closely approximated in order to significantly decrease the CG error. For the equally spaced eigenvalues the CSD seems visually well approximated. But a closer look reveals rather slow convergence of the Ritz values to the individual eigenvalues, which proceeds from both edges of the spectrum. For more on the convergence of Ritz values in this case see \cite{DruKni95,GreStr92}, \cite[Lecture 36]{TreBau97}, and, in the asymptotic case, \cite{Kui00}. Further interesting points will occur when the same experiments will be performed in finite precision arithmetic; see Section \ref{sec:fpcg} below.

\subsection{Mathematical behavior of CG for matrices with clustered eigenvalues}\label{sec:cgclust}
\phantom{}
\smallskip
\begin{tcolorbox}
\emph{Main point: A spectrum localized in $\ell$ tight clusters does not mean reaching mathematically a good CG approximation to the solution in $\ell$ steps. The position of clusters is essential. }
\end{tcolorbox}

\smallskip
\emph{Setup:} We generate three auxiliary diagonal matrices via \eqref{eq:diagmatrix} with different  parameters $\rho$ to control the eigenvalue distributions. All matrices have $N=10$, $\lambda_1=0.1$, and $\lambda_N=10^3$. The first matrix uses $\rho=0.6$ and a slight modification of \eqref{eq:diagmatrix} so that eigenvalues accumulate to the right. The second matrix uses $\rho=0.6$ with eigenvalues accumulated to the left. The third matrix uses $\rho=1.0$, which gives equally spaced eigenvalues. For each auxiliary matrix, we then construct a new matrix of size $N=100$, which is used in the experiment, by replacing each of the eigenvalues (diagonal entries) by a tight cluster of $10$ eigenvalues with spacing $10^{-12}$. Thus our matrices have 10 clusters, each with 10 eigenvalues, with cluster diameter $O(10^{-11})$. In each case we use $b = [1,\dots,1]^T/\sqrt{N}$ and $x_0=0$.

\begin{figure}
\centering
\includegraphics[width=0.55\textwidth]{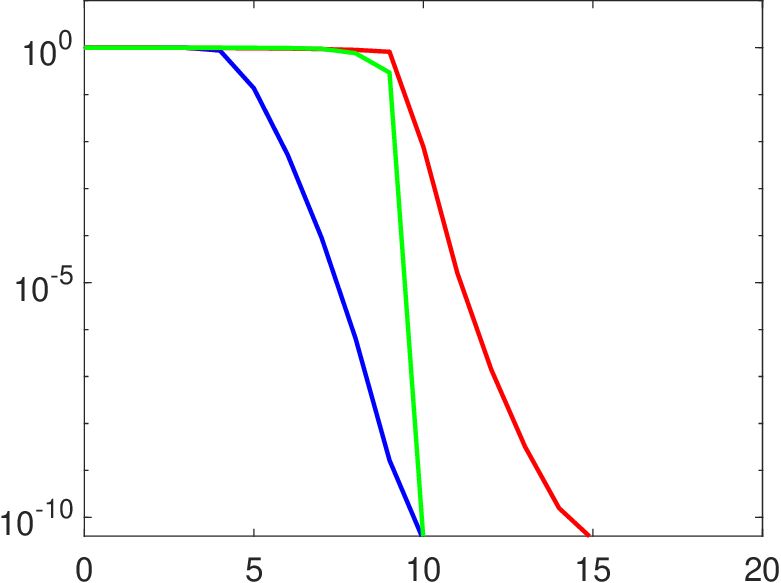}
\caption{The relative error in the $A$-norm for exact CG run on three problems with matrices having different distributions of 10 eigenvalue clusters, where each cluster contains 10 eigenvalues. }\label{fig:cgclust}
\end{figure}

\begin{figure}
\centering
\includegraphics[trim={4cm 3cm 4cm 3cm},clip,width=0.575\textwidth]{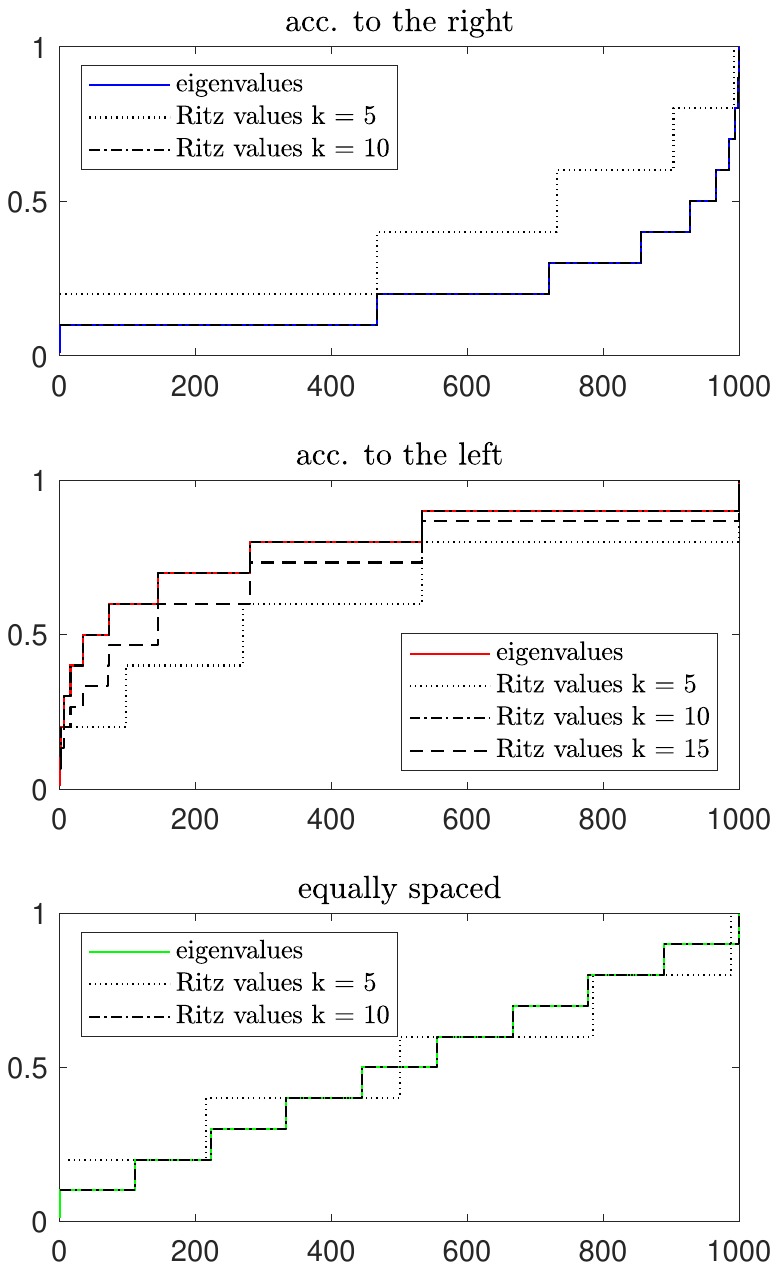}
\caption{Cumulative spectral density plots for Figure \ref{fig:cgclust}.}\label{fig:specdensclust}
\end{figure}

\smallskip
\emph{Observations:}
In Figure \ref{fig:cgclust} we plot the convergence of exact CG for the three problems. Accompanying CSD plots are given in Figure \ref{fig:specdensclust}. The matrix has in each case 10 tight clusters of eigenvalues. When the clusters of eigenvalues are accumulated to the right (blue) and when they are the equally spaced (green), the relative error in the $A$-norm reaches in 10 iterations the level below $10^{-10}$. When the clusters are accumulated to the left (red), the relative error in the $A$-norm makes no progress in 10 iterations. Note that this behavior contradicts the widespread general claims about clustering of eigenvalues and CG convergence, which ignore the positions of clusters.\footnote{We use on purpose a very small and simple example. It is easy to find examples for which the same observations are substantially more pronounced.}

\smallskip
\emph{Explanation:}
In the first 10 iterations the CG polynomials for all cases place a single Ritz value in each cluster. For the clusters accumulated to the right as well as equally spaced this is sufficient for approximating the minimal polynomial of the matrix (which is of degree 100) {\em in the sense of \eqref{eqn:CG_Ritz_values}}. For the clusters accumulated to the left, placing one Ritz value in each cluster is not enough to significantly decrease the error, and the minimal polynomial of the matrix is in the same sense not well approximated, despite the seemingly analogous position of Ritz values; see the CSD plots, where for $k=10$ the solid lines and the associated dash-dotted lines graphically coincide. To achieve the desired decrease of the error for the case of clusters accumulated to the left, CG must place additional Ritz values in the rightmost clusters, which delays convergence. At iteration 15, the five rightmost clusters contain two Ritz values each, and the dashed line representing the CSD for $k=15$ departs from the CSD representing the matrix spectrum. If the computation proceeds, then this departure would become more and more significant because more and more Ritz values will be placed in the rightmost clusters.

This mechanism has been demonstrated in \cite{GreStr92} and it was further thoroughly explained in \cite[Section~5.6]{LieStr13}. A very detailed account of the relationship between preconditioning and the clustering argument can be found in \cite[Section~3(c)]{CarStr20}. In order to avoid misunderstandings, we again emphasize that this subsection has dealt with the mathematical behavior of CG. The effects of rounding errors on the convergence are examined next. As we will see, clusters of eigenvalues will then play a specific fundamental role.

\subsection{Sensitivity of CG to rounding errors}\label{sec:fpcg}
\phantom{}
\smallskip
\begin{tcolorbox}
\emph{Main point: Particular eigenvalue distributions, specifically in cases of large outlying eigenvalues, cause CG convergence to be more susceptible to delay caused by finite precision errors. Convergence behavior for finite precision CG can be equated (up to an unimportant difference) with exact CG on a larger problem, whose eigenvalues are replaced by tight clusters.}
\end{tcolorbox}

\begin{figure}
\centering
\includegraphics[width=0.475\textwidth]{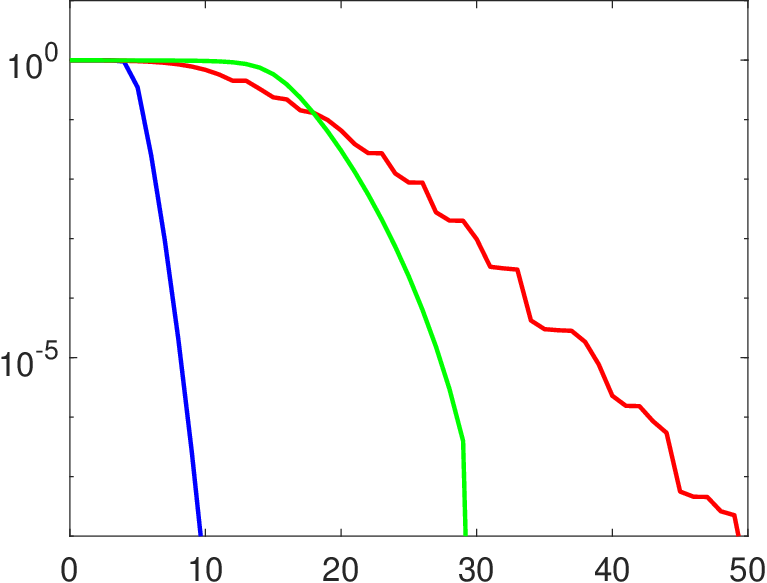}\hfill 
\includegraphics[width=0.475\textwidth]{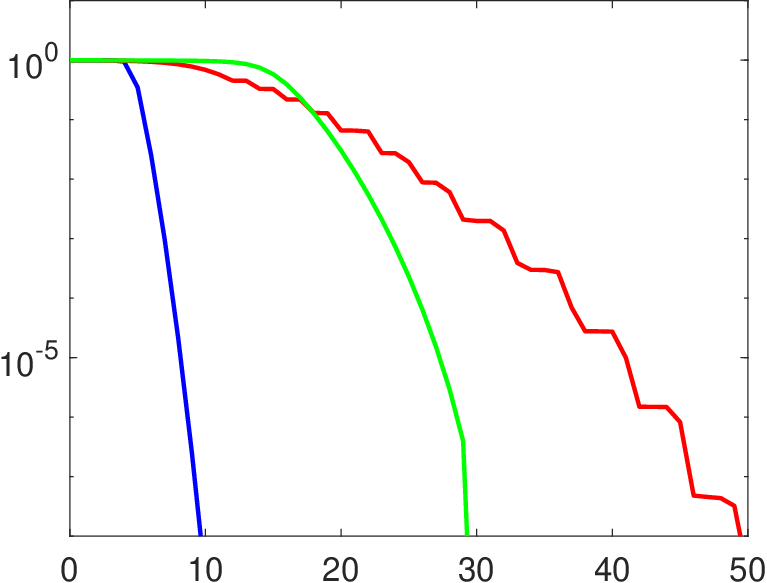}
\caption{Left: The relative error in the $A$-norm for CG in finite precision run on three problems with matrices having different eigenvalues distributions, corresponding to those in Figure~\ref{fig:cgeigdist}. Right: The relative error in the $A$-norm for exact CG run on three problems with matrices with the same eigenvalue distributions as in the left plot, but with each eigenvalue replaced by a tight cluster of 4 eigenvalues.}\label{fig:fpeigdist}
\end{figure}

\begin{figure}
\centering
\includegraphics[trim={4cm 3cm 4cm 3cm},clip,width=0.575\textwidth]{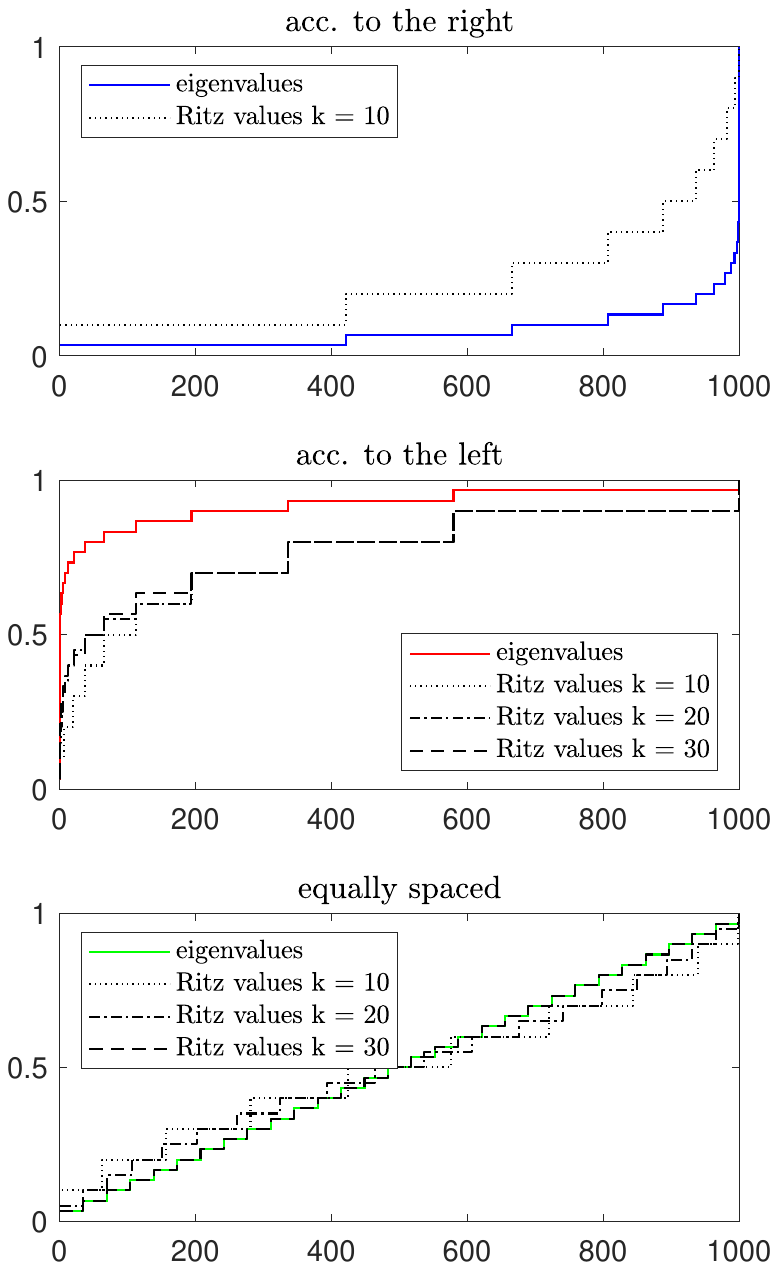}
\caption{Cumulative spectral density plots for the left part of Figure \ref{fig:fpeigdist}, i.e., the Ritz values are determined from the finite precision CG run. Compare with Figure \ref{fig:specdensexact}.}\label{fig:specdens}
\end{figure}

\smallskip
\emph{Setup:} We use the same three diagonal matrices and the same right-hand sides as in Section \ref{sec:cgeigdist}, but now we run CG in finite (double) precision. 
We plot the resulting convergence curves on the left in Figure \ref{fig:fpeigdist}, and the CSDs at certain iterations for each problem in Figure \ref{fig:specdens}.

Then, for each diagonal matrix, we create a larger matrix by replacing each eigenvalue with a tight cluster of 4 eigenvalues. The spacing between eigenvalues in a cluster is $10^{-13}$. We run exact CG for these problems and plot the resulting convergence curves on the right in Figure \ref{fig:fpeigdist}.

\smallskip
\emph{Observations:} A comparison of the green and blue curves in Figures~\ref{fig:cgeigdist} and~\ref{fig:fpeigdist} (left) shows that the convergence of CG for the matrices with eigenvalues accumulated to the right of the spectrum (blue) and for equally distributed eigenvalues (green) are essentially not affected by finite precision errors. On the other hand, for the matrix with eigenvalues accumulated to the left and only a few large outlying eigenvalues, the finite precision CG suffers from a significant delay of convergence (red curves). Moreover, a comparison of the two plots in Figure~\ref{fig:fpeigdist} shows that the behavior of finite precision CG (left) is remarkably similar to the behavior of exact CG where the eigenvalues are replaced by tight clusters (right).

\smallskip
\emph{Explanation:} There are two phenomena working against each other here. Whereas large outlying eigenvalues are desirable in exact arithmetic, they cause the problem to be more sensitive to rounding errors, which can result in convergence delay in finite precision computations. This phenomenon was investigated in~\cite{Str91}, which was inspired by the earlier discussion by Jennings~\cite{Jen77}, who related the convergence of CG to a polynomial curve fitting problem. It can be nicely viewed via the CSD plots in Figure~\ref{fig:specdens}. While for the eigenvalues accumulated to the right and for equally distributed eigenvalues there is no observable difference between the exact and finite precision CG computations (compare the top and the bottom plots in Figures~\ref{fig:fpeigdist} and~\ref{fig:specdensexact}), the CSDs associated with the eigenvalues accumulated to the left are remarkably different. There is no chance that for this problem the CSD determined by the eigenvalues can be closely approximated by the CSD generated by the Ritz values resulting from the finite precision CG computation. The plot in Figure~\ref{fig:fpeigdist} shows that with increasing iteration number more and more Ritz values have to be placed close to the rightmost outlying eigenvalues.

A theoretical explanation is provided by the work of Greenbaum~\cite{Gre89} that was further thoroughly illustrated in~\cite{GreStr92} and extended by Notay~\cite{Not93}. Additional Ritz values close to the large outlying eigenvalues have to appear in order to eliminate excessively large gradients of the CG approximation polynomials, which would otherwise occur in their neighborhood; see \cite[Section~5.9]{LieStr13} and \cite[Section~5]{MeuStr06}. For the accumulation of the eigenvalues to the right and for equally distributed eigenvalues the gradient of the CG approximation polynomial near all eigenvalues is sufficiently bounded without a need for placing additional Ritz values in their neighborhoods. This explains the numerical stability of CG for these problems; see also~\cite[Sections 3(b)(i)]{CarStr20}. It is also worth recalling the arguments in Section~\ref{sec:cgclust} above that deal with the mathematical behavior for problems with tight clusters of eigenvalues.

\subsection{Preconditioned CG and the condition number}\label{sec:pcg}
\phantom{}
\smallskip
\begin{tcolorbox}
\emph{Main point: The nonlinear adaptivity of CG to the location of the individual eigenvalues indicates that a smaller condition number does not necessarily lead to faster convergence (contrary to widespread misinterpretations of the $\kappa(A)$-bound in the literature). Therefore it is not recommended to use the minimization of the condition number as the only criterion for the choice of preconditioners. Alternatives to the condition number exist, but they require deep knowledge of CG and of the problem to be solved.
}
\end{tcolorbox}

\smallskip
\emph{Setup:} We perform two experiments. First, we define the matrix $A$ to be of the form \eqref{eq:diagmatrix} with $N=40$, $\lambda_1=10^{-3}$, $\lambda_N=100$, and $\rho=0.1$. 
Thus, the eigenvalues accumulate at the lower end of the spectrum of $A$, and we have
$$\kappa(A)=10^5.$$
We consider the diagonal preconditioner $P$ so that $P^{-1}A$ is a diagonal matrix with eigenvalues equally spaced between $\lambda_1 = 10$ and $\lambda_N = 100$, and hence
$$\kappa(P^{-1}A)=10.$$
We apply exact CG with $x_0=0$ to $Ax=b$ and $P^{-1}Ax=P^{-1}b$, where $b=[1,\dots,1]^T/\sqrt{N}$.

Second, following \cite[Section~5.3]{MorNocSie02}, we consider the boundary value problem
\begin{equation}
\label{eq:motivation}
	- \nabla \cdot (k(x) \nabla u) = 0\;\; \text{in} \;\; \Omega=(-1,1) \times (-1,1), 
 \quad u = u_D \;\; \text{on} \;\; \partial \Omega,
\end{equation}
where the domain $\Omega$ is divided into four subdomains $\Omega_1,\Omega_2,\Omega_3,\Omega_4$ corresponding to the four axis quadrants numbered counterclockwise. Let $k(x)$ be piecewise constant on the individual subdomains with $k_1 = k_3 \approx 161.45$ and $k_2 = k_4 = 1$. The Dirichlet boundary conditions are described in \cite[Section~5.3]{MorNocSie02}. We use the linear finite element discretization with the standard uniform triangulation and $N = 3969$ degrees of freedom. We apply CG in finite (double) precision with $x_0=0$ to the unpreconditioned system, and to preconditioned systems with the algebraic incomplete Cholesky factorization preconditioning (ICHOL) of the matrix $A$ with no fill-in, with ICHOL and the drop-off tolerance $10^{-2}$, and with the Laplace operator preconditioning. The corresponding matrices have the following condition numbers:
\begin{align*}
\kappa(A) &\approx 6750 &&\mbox{(unpreconditioned),}\\
\kappa(P^{-1}A) &\approx 431 &&\mbox{(ICHOL with no fill-in),}\\
\kappa(P^{-1}A) &\approx 16 &&\mbox{(ICHOL with drop-off tolerance $10^{-2}$),}\\
\kappa(P^{-1}A) &\approx 160 &&\mbox{(Laplace operator preconditioning).}
\end{align*}
The same setting was used for the motivating example in~\cite{GerMarNieStr19}, where one can find a more detailed description.

\begin{figure}
\centering
\includegraphics[width=0.475\textwidth]{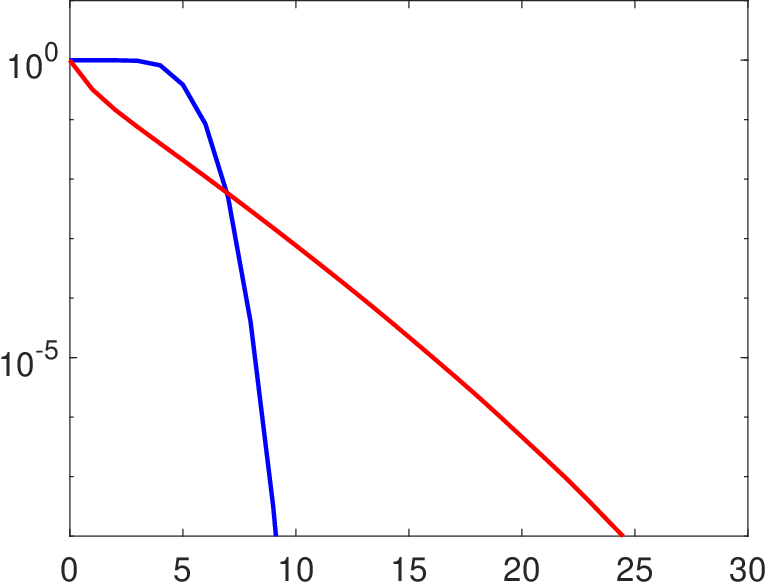}\hfill
\includegraphics[width=0.475\textwidth]{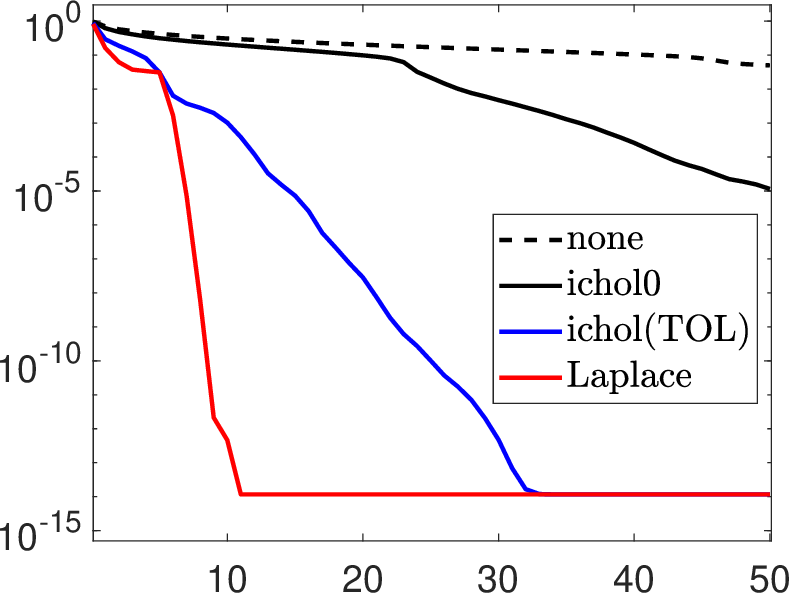}
\caption{Left: The relative error in the $A$-norm for unpreconditioned exact CG (blue) and preconditioned exact CG (red) for a linear algebraic system with $A$ of the form \eqref{eq:diagmatrix}. Right: The relative error in the $A$-norm for finite precision computations with unpreconditioned CG (dashed), ICHOL preconditioning with no fill-in (black), ICHOL preconditioning with drop-off tolerance $10^{-2}$ (blue), and Laplace preconditioning (red), applied to the discretization of \eqref{eq:motivation}.}\label{fig:pcgvscg}
\end{figure}

\smallskip
\emph{Observations:} First example: In the left part of Figure~\ref{fig:pcgvscg} we plot the relative error $A$-norms for exact CG on $Ax=b$ (blue) and $P^{-1}A=P^{-1}b$ (red). The number of iterative steps for reaching the accuracy level $10^{-8}$ is for the preconditioned system three times larger than for the unpreconditioned system, although $\kappa(P^{-1}A)$ is four orders of magnitude smaller than $\kappa(A)$. (Note that the preconditioning works well when the desired accuracy level is only on the order $10^{-2}$.) 

Second example: In the right part of Figure~\ref{fig:pcgvscg} we plot the relative error $A$-norms for finite precision CG applied to the unpreconditioned and the three preconditioned systems. A comparison of the results for the unpreconditioned system and the two ICHOL preconditioniners seems to confirm the often repeated claim that ``reducing the condition number \emph{implies} faster convergence of CG''. However, the fastest convergence occurs for the Laplace preconditioning, although in this case the condition number is an order of magnitude larger than for the ICHOL preconditioner with drop-off tolerance $10^{-2}$. It is important to note also the dramatic acceleration of CG with the Laplace preconditioning after the fifth iteration.

\smallskip
\emph{Explanation:} Both examples show that a smaller condition number does not necessarily imply faster convergence of CG. In particular, in the second example the decrease of the condition number using ICHOL with a larger drop-off tolerance is just a side effect and not the driving force of faster convergence. In this case the improvement is due to a better approximation of the true Cholesky factors of $A$. Clearly, the frequently repeated claim that the goal of preconditioning should be to reduce the condition number of the coefficient matrix (thus reducing the $\kappa(A)$-bound), and that doing so \emph{guarantees} an improvement of the CG convergence is false.\footnote{In this context we point out the monograph of Hackbusch~\cite{HacBook94} which covers, among many other topics, CG in a truly insightful way. Its Section~9.4.3 ``Convergence Analysis'' addresses also the difference between the Chebyshev method with the bound \eqref{eqn:rhoAbound1}, and CG with the application of the same bound.}
However, one should always take into account the context of the problem to be solved, in particular the desired accuracy of the computed approximation. 

A detailed explanation of the second example is out of the scope of this paper, but can be found, including the effects of rounding errors, in~\cite{GerMarNieStr19}. It shows why for the Laplace operator preconditioning the CG convergence \emph{must exactly after five initial iterations} accelerate so rapidly that after a few subsequent iterations an approximate solution is found with an accuracy close to the machine precision level. It also provides an accurate estimate for the speed of the accelerated convergence. For the ICHOL preconditioning it explains why an analogous acceleration can not take place, and also provides an accurate approximation of convergence rate.

As mentioned in the Introduction, it may seem that the design of preconditioners in practical computations using CG have to rely on the condition number because there is no viable alternative. The recent works \cite{GerMarNieStr19,GerNieStr22,Lad21,Pul21,Leu22,Lad23} show that for self-adjoint second order elliptic PDE boundary value problems with the operator of the form $-\nabla (k(x) \nabla u)$ such an alternative exists. Analysis of the motivating problem in \cite{GerMarNieStr19} suggests a possible path for further research. The convergence behavior can be anticipated \emph{a priori} based on a low-cost approximation of all eigenvalues of the preconditioned matrix before any preconditioned CG computation starts. We therefore believe that the approach in \cite{GerMarNieStr19}, together with the other given references, question the status quo. 

The latest paper in this line of development \cite{NieStr24} suggests further arguments for combining infinite dimensional operator reasoning with algebraic considerations about solving algebraic systems resulting from discretizations. It also formulates open problems that should be addressed when considering eigenvalues of preconditioned matrices in relation to the spectra of the associated infinite dimensional operators on Hilbert spaces. Section \ref{sec:intro2} below indicates that it is promising to investigate more complex approaches as an alternative to decreasing the condition number in preconditioning practical hard problems.

\subsection{Simple model problems and the practical performance of CG}
\label{sec:cgpoisson}
\phantom{}
\smallskip
\begin{tcolorbox}
\emph{Main point: Model problems that are used out of context are not indicative of the behavior of CG in solving practical problems. They can also complicate understanding of important CG features.}
\end{tcolorbox}

\smallskip
\emph{Setup:}
We consider the 2D Poisson problem $- \Delta u = f$ in $\Omega=(0,1)\times (0,1)$, 
$u = 0$ on $\partial \Omega$, and $f$ is constant. This boundary value problem is discretized using the five-point finite differences and the $50 \times 50$ grid, giving $N = 2500$ degrees of freedom. The matrix $A$ is simply generated using the MATLAB command {\tt gallery('poisson',50)}, and we use $b=[1,\dots,1]^T/\sqrt{N}$. We apply exact CG and CG in finite (double) precision, both with $x_0=0$.

\begin{figure}
\centering
\includegraphics[width=0.55\textwidth]{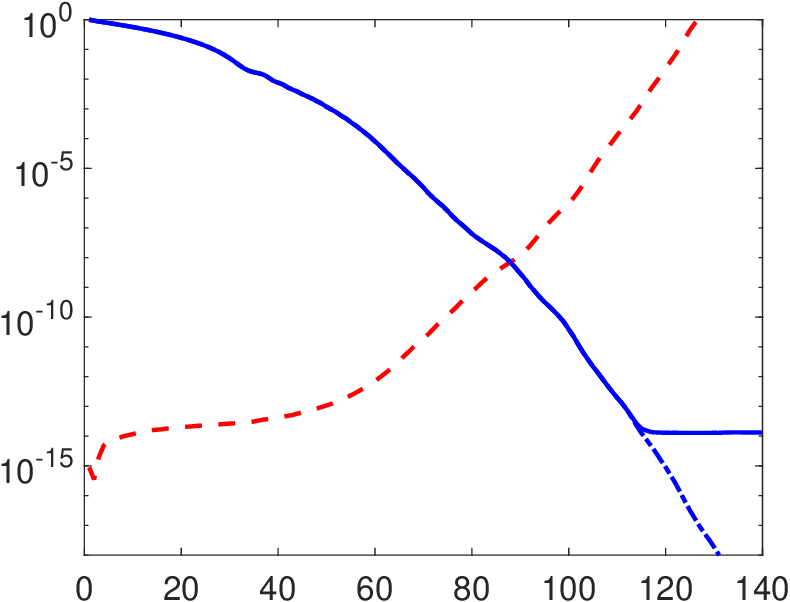}
\caption{The relative error in the $A$-norm for exact CG (dash-dotted blue) as well as for CG in finite precision (solid blue), and the loss of orthogonality among Lanczos basis vectors for CG in finite precision (dashed red), for a 2D Poisson model problem.}
\label{fig:poissonCG}
\end{figure}

\smallskip
\emph{Observations:}
In Figure~\ref{fig:poissonCG} we plot the relative $A$-norm of the error for both exact and finite precision CG. For finite precision CG we plot also the loss of orthogonality among the Lanczos basis vectors, measured by $\Vert I-V_k^T V_k\Vert_F$. Its growth seems to mirror the convergence of the relative error in the $A$-norm, which is almost the same for both exact and finite precision CG until the latter reaches its maximal attainable accuracy.

\smallskip
\emph{Explanation:}
For the given model problem, the eigenvalues are almost uniformly distributed. The loss of orthogonality among the Lanczos basis vectors is gradual. Since there is no loss of rank in the computed basis until the finite precision CG reaches its final accuracy, there is no delay of the CG convergence; see \cite[Section 5.9.4]{LieStr13} for a more detailed explanation.\footnote{Paige proved that the loss of orthogonality can occur only in the directions of the converged Ritz vectors; see \cite{Pai80}, and \cite[pp.~504-508]{MeuStr06} for an explanation of misinterpretations of this breakthrough result as ``convergence implies loss of orthogonality.'' Since there is no loss of ortogonality until CG reaches its final accuracy, no Ritz pair can approximate a matrix eigenpair with accuracy proportional to the machine precision level.} Such behavior can not be extrapolated to CG behavior in solving practical problems.

\subsection{Computational behavior of different CG algorithms}\label{sec:2vs3term}
\phantom{}
\smallskip
\begin{tcolorbox}
\emph{Main point: Rounding errors cause convergence delay and affect the attainable accuracy. The magnitude of these effects, and in the case of attainable accuracy, the mechanism, depends on the particular algorithm/implementation of CG.}
\end{tcolorbox}

\smallskip
\emph{Setup:}
We use a diagonal matrix $A$ as defined in \eqref{eq:diagmatrix} with $N=48$, $\lambda_1 = 0.1$, $\lambda_N=10^3$, and $\rho=0.25$. We also use $b=[1,\dots,1]^T/\sqrt{N}$ and $x_0=0$. We test two mathematically equivalent algorithmic variants of CG: the variant of Hestenes and Stiefel \cite{HesSti52} which uses three 2-term recurrences (see Algorithm \ref{alg:hscg}), and a different variant which uses two 3-term recurrences.

\begin{figure}
\centering
\includegraphics[width=0.55\textwidth]{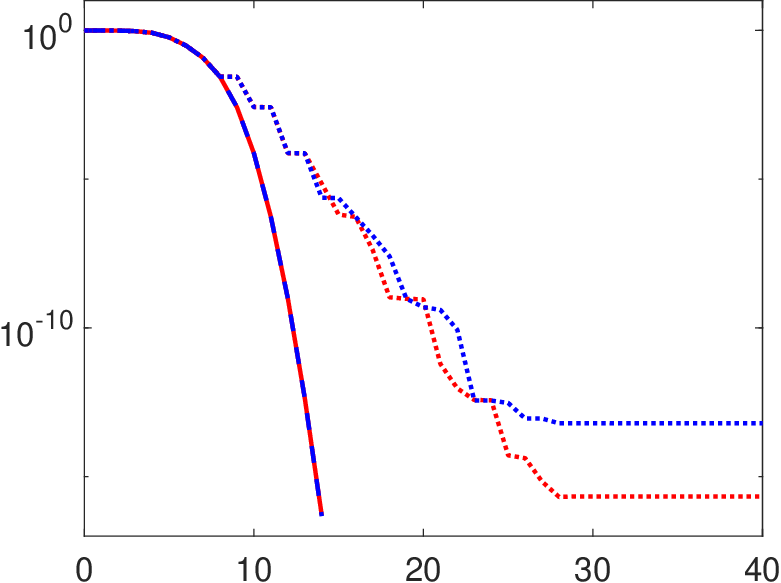}
\caption{The relative error in the $A$-norm for the two variants of CG in exact arithmetic (solid red and dashed blue) and in finite precision (2-term: dotted red; 3-term: dotted blue).}\label{fig:fp2v3}
\end{figure}

\smallskip
\emph{Observations:} The results comparing the relative $A$-norm of the error for exact and finite (double) precision computations are shown in Figure \ref{fig:fp2v3}. In exact arithmetic the $A$-norm error curves of the 2-term and 3-term variants are obviously identical. In finite precision the convergence is delayed. The delay is slightly worse in the 3-term variant, and the final accuracy level attained by this variant is over two orders of magnitude worse than the level attained by the 2-term variant.

\smallskip
\emph{Explanation:} In the 2-term recurrences the loss of accuracy is caused by a simple accumulation of local rounding errors, but in the 3-term recurrences these local rounding errors can be \emph{substantially amplified}. This behavior is analyzed together with its dependence on the initial residual in \cite{GutStr00}.

Despite maintaining mathematical equivalence to Algorithm \ref{alg:hscg}, any algorithm that reorders computations or introduces auxiliary quantities can have a different computational behavior. Examples of CG algorithms designed for high-performance parallel environments include $s$-step (communication-avoiding) CG and pipelined CG, both of which are subject to potential amplification of local rounding errors and thus more substantial delays of convergence and worse maximal attainable accuracy than Algorithm \ref{alg:hscg}; see, e.g., \cite{Car15, CarRozStrTicTum18, CooYetAguGirVan18}.

\subsection{Residual versus error, and stopping criteria for CG}\label{sec:resvserror}
\phantom{}
\smallskip
\begin{tcolorbox}
\emph{Main point: Unlike the $A$-norm (or energy norm) of the error, the residual $2$-norm does not have any physical meaning. Moreover, if the matrix condition number is large, then the residual $2$-norm is not a reliable indicator of the error. Theoretically justified estimators for the $A$-norm and the $2$-norm of the error are available at a negligible computational cost.}
\end{tcolorbox}

\smallskip
\emph{Setup:} We follow Meurant~\cite{Meu20} for constructing two linear algebraic systems of size $N=20$ such that the trajectories of the residual $2$-norm and the $A$-norm of the error in CG are prescribed. For the first, the residual norms $\Vert r_k\Vert_2$ oscillate between 1 and 2, and the errors are $\Vert e_0 \Vert_A = \Vert x-x_0\Vert_A =1$, $\Vert e_k\Vert_A = 0.4 \cdot \Vert e_{k-1}\Vert_A$  for $k=1,\ldots,N-1$. For the second, the residual norms are $\Vert r_0 \Vert_2 = 1$, $\Vert r_k\Vert_2 = 0.4\cdot \Vert r_{k-1}\Vert_2$ for $k=1,\ldots, N-1$, and the errors are $\Vert e_0 \Vert_A=1$, $\Vert e_k\Vert_A = 0.999 \cdot \Vert e_{k-1}\Vert_A$ for $k=1,\ldots,N-1$. 

To construct the systems we set $\nu_k = 1/\Vert r_k \Vert_2$, $k=1,\ldots, N-1$, and $\sigma_k = \Vert e_k \Vert_A^2 / (\Vert r_k\Vert_2 \Vert r_0\Vert_2)$,  $k=0,\ldots, N-1$, and we then set
\[
L = \begin{bmatrix}
\sigma_0 & & & \\
\sigma_1 & \sigma_1 \nu_1 & & \\
\vdots & \vdots& \ddots & \\
\sigma_{N-1} & \sigma_{N-1} \nu_1 & \cdots & \sigma_{N-1}\nu_{N-1}
\end{bmatrix}.
\]
We apply exact CG to $A= (L+\widehat{L}^T)^{-1}$ and $b=e_1$, where $\widehat{L}$ is the strictly lower triangular part of $L$ and $e_1$ is the first column of the identity matrix.

\begin{figure}
\includegraphics[width=0.475\textwidth]{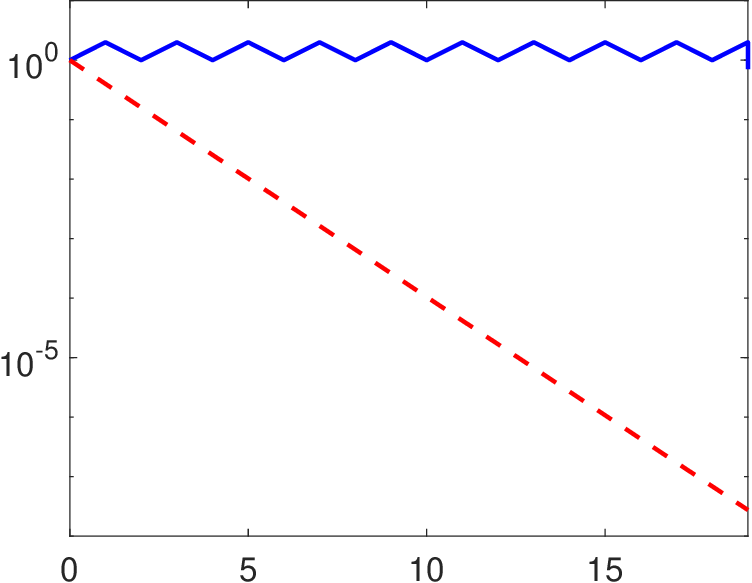}\hfill
\includegraphics[width=0.475\textwidth]{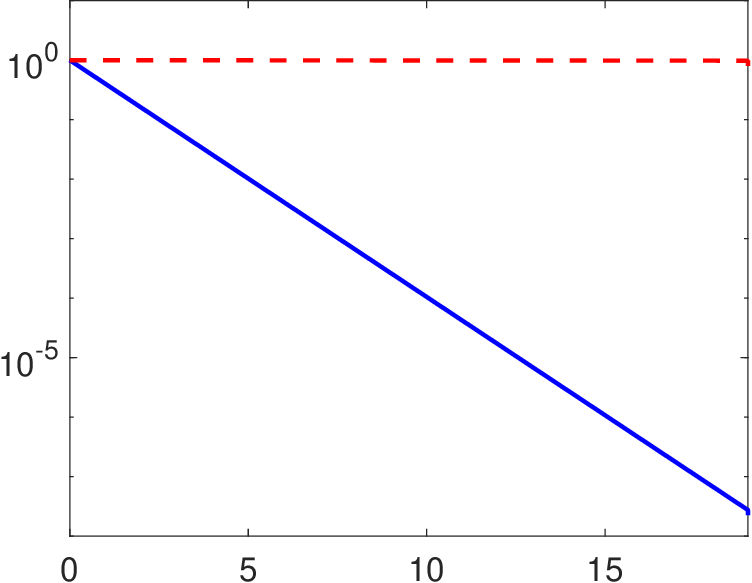}
\caption{Comparison of the relative residual $2$-norms (solid blue) and the relative $A$-norms of the error (dashed red) of exact CG for two linear algebraic systems.}\label{fig:resvserror}
\end{figure}

\smallskip
\emph{Observations:}
The relative residual 2-norm and error $A$-norms are shown in Figure \ref{fig:resvserror}. For the first system (left plot) the residual norms stagnate, apart from oscillating between 1 and 2, whereas the error norms decrease linearly. For the second system (right plot) the error norms almost stagnate (they must be strictly decreasing), whereas the residual norms are decreasing relatively quickly.

\smallskip
\emph{Explanation:}
Hestenes and Stiefel comment in their original paper \cite[Section 18]{HesSti52} that for \emph{any} prescribed sequence of residual $2$-norms, there exists a symmetric positive definite matrix $A$ and right-hand side $b$ such that CG exhibits the prescribed convergence behavior. Thus, in general circumstances, one cannot equate a small residual norm with a small error. The converse also does not hold: a large residual does not imply a large error. Note that because the convergence of CG is linked with the distribution of the eigenvalues of $A$, it is not possible to also  simultaneously prescribe the eigenvalues of $A$. This is in contrast to the GMRES method; see Example~\ref{sec:GMRESevery}. If $\kappa(A)$ is small, then the difference between the residual $2$-norm and the $A$-norm (or the $2$-norm) of the error is not substantial. This can happen in practical applications. But it can also happen that $\kappa(A)$ is large and there is no practically applicable preconditioning available that makes $\kappa(P^{-1}A)$ small. Therefore the question on what should be used for stopping criteria in CG is relevant.

The relative residual norm is inexpensive to compute, and we already have the recursively computed residual available in each iteration of CG.\footnote{When CG is used as an inner solver embedded in an outer loop as, for example, in nonlinear optimization, this can be a viable option because the inner-outer heuristic stopping criteria are consistent.} 
However, the norm of the residual has no physical meaning, while the $A$-norm of the error, which is minimized in each step by CG, represents in many applications the (discretization of) the energy norm. In practice, we of course do not know the true solution $x$, and thus we can not compute the error $e_k$. This was commented on already by Hestenes and Stiefel, who gave formulas for estimating the error norm; see~\cite[Section~4]{HesSti52}. Since then, much research has focused on developing reliable error norm estimation and associated stopping criteria for CG; see, e.g., \cite{GolMeu97, GolStr94, MeuPapTic21, StrTic02, StrTic05}. Very useful estimates of both the $A$-norm and the 2-norm of the error are available at a negligible computational cost, and they are theoretically guaranteed to hold also in cases with severe effects of rounding errors. The associated software realizations are simple and freely available, and hence there are strong arguments for considering such error estimates in practical computations.

\subsection{The trajectory of finite precision CG computations}\label{sec:cgtraj}
\phantom{}

\smallskip
\begin{tcolorbox}
\emph{Main point:
The approximate solutions produced by finite precision CG can be closely mapped to those produced by exact CG via a mapping defined by examining the rank deficiency of the Krylov subspace basis. It seems that the trajectory of the approximate solutions produced by finite precision CG remain in a narrow ``tunnel'' around those produced by exact CG.}
\end{tcolorbox}

\smallskip
\emph{Setup:} We generate a diagonal matrix $A$ as defined in \eqref{eq:diagmatrix} with $N=35$, $\lambda_1 = 0.1$, $\lambda_N=10^2$, and $\rho=0.65$. We use $b=[1,\dots,1]^T/\sqrt{N}$ and $x_0=0$. We run exact CG and CG in finite (double) double precision arithmetic. Following \cite[Section 5.9.1]{LieStr13} and \cite{Ger13}, the finite precision CG iterates are then shifted as follows: Consider the sequence
\begin{equation}
\ell(k) = \max\left\{i\hspace{2pt} \vert \hspace{2pt} \text{rank}_t(\mathcal{K}_i(A,r_0))= k \right\}, \quad k=1,2,\ldots,
\label{eq:mapping}
\end{equation}
for the (inexact) Krylov subspace $\mathcal{K}_i(A,r_0)$ computed in double precision. To compute $\text{rank}_t(\mathcal{K}_i(A,r_0))$, we use the built-in MATLAB function \texttt{rank} for determining the numerical rank with the threshold $t=10^{-1}$; compare with~\cite[Section~5.9.1]{LieStr13}, in particular, Figure 5.17.\footnote{As discussed in \cite{LieStr13}, in this way we count the number of basis vectors that are ``strongly linearly independent''. Different value of the threshold will illustrate the same phenomena with a slightly worse quantitative match.} For exact CG iterates $x_k$ and finite precision CG iterates $\bar{x}_{\ell(k)}$, we measure the ratios
\begin{equation}
    \frac{\Vert x- \bar{x}_{\ell(k)}\Vert_A}{\Vert x-x_k \Vert_A}
    \label{eq:ratio1}
\end{equation}
and
\begin{equation}
   \left\vert 1- \frac{\Vert x- \bar{x}_{\ell(k)}\Vert_A}{\Vert x-x_k \Vert_A} \right\vert.
    \label{eq:ratio2}
\end{equation}

\begin{figure}
\centering
\includegraphics[width=0.475\textwidth]{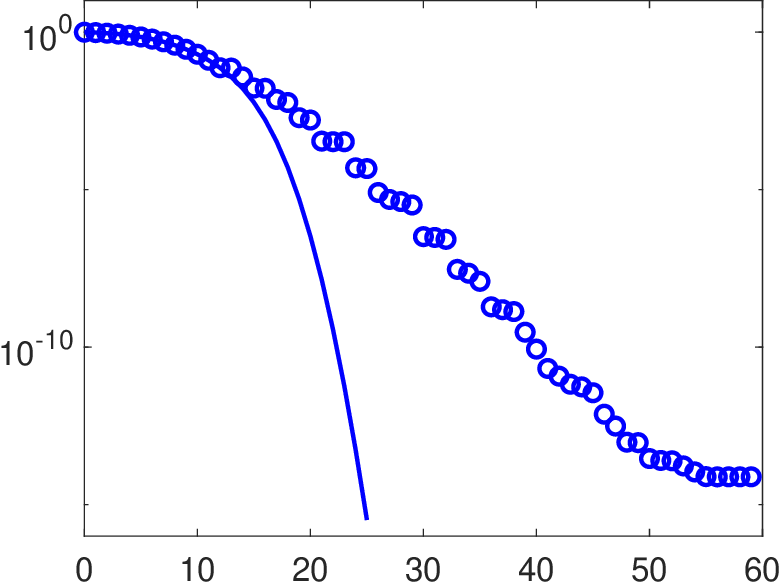}\hfill
\includegraphics[width=0.475\textwidth]{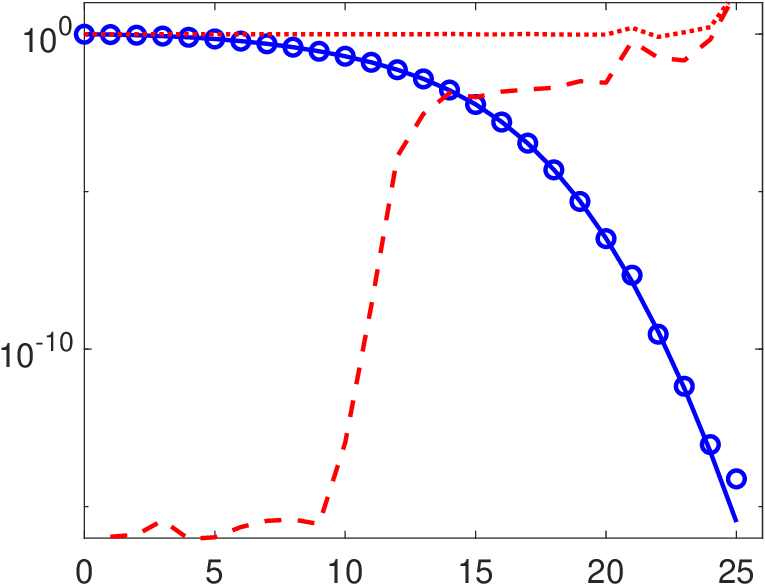}
\caption{Left: The relative error in the $A$-norm for exact CG (solid blue) and finite precision CG (blue circles). Right: The relative error in the $A$-norm for exact CG (solid blue) and \emph{shifted} finite precision CG (blue circles), the ratio \eqref{eq:ratio1} (dotted red), and the ratio \eqref{eq:ratio2} (dashed red).}\label{fig:cgtraj}
\end{figure}

\smallskip
\emph{Observations:} In the left part of Figure \ref{fig:cgtraj} we plot the relative error $A$-norms of the exact CG iterates $x_k$ (solid blue) and the finite precision CG iterates $\bar{x}_k$ (blue circles). In the right part of Figure \ref{fig:cgtraj} the solid blue curve remains the same as in the left part (apart from the change of the horizontal scale), while the blue circles now show the relative error $A$-norms for the shifted finite precision CG iterates $\bar{x}_{\ell(k)}$. We also plot the ratios \eqref{eq:ratio1} (dotted red) and \eqref{eq:ratio2} (dashed red). We see that using the mapping \eqref{eq:mapping}, the finite precision CG iterates match well with the exact CG iterates. The ratio \eqref{eq:ratio1} stays close to one throughout the computation. The ratio \eqref{eq:ratio2} starts close to machine precision and grows as the iteration proceeds, but stays below one until the finite precision CG gets close to the final accuracy level.

\smallskip
\emph{Explanation:}
The fact that the ratio \eqref{eq:ratio1} remains close to one means that the convergence trajectories for exact CG and the shifted finite precision CG are almost identical. This means that finite precision CG computations follow closely the trajectory of exact CG computations, but with progress delayed according to the approximate rank deficiency of the computed (inaccurate) Krylov subspaces. The ratio \eqref{eq:ratio2} tells us how far \eqref{eq:ratio1} is from one. Rewriting \eqref{eq:ratio2}, we can say that if
\begin{equation}
\left\vert \frac{\Vert x-x_k\Vert_A - \Vert x- \bar{x}_{\ell(k)}\Vert_A}{\Vert x-x_k\Vert_A} \right\vert \ll 1,
\label{eq:tunnel}
\end{equation}
then the trajectory of finite precision CG iterates is within a narrow ``tunnel'' around the trajectory of exact CG iterates, where the diameter of the ``tunnel'' is given by the left-hand side of \eqref{eq:tunnel}. Although \eqref{eq:tunnel} holds for this particular example, it has not been shown to hold for finite precision computations in general. We strongly suspect that such a result holds; initial investigations have been carried out in, e.g., \cite{Ger13} and \cite{GerHneKub17}. The analysis is complicated by the fact that in finite precision computations we can not easily compare the exact Krylov subspaces with the computed subspaces, since the latter ones depart from Krylov subspaces due to the effects of rounding errors.

\section{The GMRES Method}\label{sec:GMRES}

The GMRES method~\cite{SaaSch86} is well defined for any linear algebraic system $Ax=b$ with a nonsingular matrix~$A\in\R^{N\times N}$ and right-hand side $b\in\R^{N}$. If $d=d(A,r_0)$ is the grade of $r_0=b-Ax_0$ with respect to $A$, then at every step $k=1,2,\dots,d$ the GMRES method constructs a uniquely determined approximation $x_k\in x_0+{\cal K}_k(A,r_0)$ such that $r_k\perp A{\cal K}_k(A,r_0)$, or equivalently
\begin{equation}\label{eqn:GMRES}
\|r_k\|_2=\min_{z\in x_0+{\cal K}_k(A,r_0)}\|b-Az\|_2=\min_{p\in P_k(0)}\|p(A)r_0\|_2.
\end{equation}
GMRES is implemented using the Arnoldi algorithm~\cite{Arn51}, in practical applications usually in the modified Gram-Schmidt (MGS) variant, which generates an orthonormal basis $V_k$ of ${\cal K}_k(A,r_0)$ and a matrix decomposition of the form $AV_k=V_{k+1}H_{k+1,k}$, where $H_{k+1,k}\in\R^{k+1,k}$ is an unreduced upper Hessenberg matrix. Then $x_k=x_0+V_kt_k$ is determined by
\begin{equation}\label{eqn:tk}
t_k=\argmin_{t\in\R^k}\|b-A(x_0+V_kt)\|_2=\argmin_{t\in\R^k} \|\|r_0\|e_1-H_{k+1,k}t\|_2.
\end{equation}
In practical implementations of GMRES, the least squares problem on the right is solved by computing the QR decomposition of the matrix $H_{k+1,k}$. Because of the upper Hessenberg structure, this decomposition can be obtained using Givens rotations, which can be updated in every step. This process also yields an update of the value $\|r_k\|_2$ without explicitly computing $x_k$. Thus, in practical implementations of GMRES, the least squares problem is solved only when the updated value of the residual norm is below a given tolerance. While all of this is very efficient, the Arnoldi algorithm for a general nonsymmetric matrix requires full recurrences (unlike the short recurrences in CG), and hence the computational cost in terms of work and storage requirements per iteration of GMRES grows significantly. As a consequence, full recurrence GMRES is for large problems typically computationally unfeasible. A common strategy, which already appeared in~\cite[Algorithm~4]{SaaSch86}, is to restart the algorithm after a certain number of steps. We do not consider restarted GMRES here, but point out that its behavior is not fully understood and sometimes counterintuitive; see, e.g.,~\cite{Emb03}.

A pseudocode implementation of GMRES is shown in Algorithm~\ref{alg:GMRES}, and more details about the implementation of GMRES can be found, e.g., in~\cite{SaaSch86} and~\cite[Section~2.5.5]{LieStr13}. A detailed analysis why the computation of orthogonal Krylov subspace bases for general nonsymmetric matrices in general requires full instead of short recurrences is given in~\cite[Chapter~4]{LieStr13}.

\begin{algorithm}
    \caption{GMRES method (pseudocode) \label{alg:GMRES}}
    \begin{algorithmic}[1]
        \Require{Nonsingular matrix $A\in \R^{N\times N}$, right-hand side $b$, initial approximation $x_0$; convergence tolerance $\tau$; maximum number of iterations $n_{\max}$.}
        \State{$r_0 = b - A x_0$}
        \For{$k=1,2,\dots,n_{\max}$}
            \State{Compute step $k$ of the Arnoldi algorithm to obtain $AV_k=V_{k+1}{H}_{k+1,k}$.}
            \State{Update the QR factorisation of ${H}_{k+1,k}$ and compute the updated $\|r_k\|_2$.}
            \State{If $\|r_k\|_2\leq\tau$, then compute $t_k$ in~\eqref{eqn:tk}, return $x_k=x_0+V_kt_k$, and stop.}
        \EndFor
    \end{algorithmic}
\end{algorithm}

If $A$ is diagonalizable, $A=X\Lambda X^{-1}$ with $\Lambda={\rm diag}(\lambda_1,\dots,\lambda_N)$, then
\begin{align}
\|r_k\|_2 &= \min_{p\in P_k(0)}\|p(A)r_0\|_2 \,
\leq \|r_0\|_2\,\min_{p\in P_k(0)}\|p(A)\|_2
\label{eqn:GMRESnorm}\\
&\leq \kappa(X)\,\|r_0\|_2\,\min_{p\in P_k(0)}\max_{1\leq i\leq N} |p(\lambda_i)|,\label{eqn:GMRESbd}
\end{align}
where $\kappa(X)=\|X\|_2\|X^{-1}\|_2$. Greenbaum and Trefethen have called the minimization problem on the right hand side of \eqref{eqn:GMRESnorm} the \emph{ideal GMRES} approximation problem, because taking the upper bound ``disentangles the matrix essence of the process from the distracting effects of the initial vector''~\cite[p.~361]{GreTre94}. For certain matrices $A$ and iteration steps $k$, however, the value of the ideal GMRES approximation problem is much larger than even the worst-case GMRES residual norm, i.e., it can happen that
$$\max_{\substack{v\in\R^N\\ \|v\|_2=1}} \min_{p\in P_k(0)} \|p(A)v\|_2 \ll
\min_{p\in P_k(0)} \max_{\substack{v\in\R^N\\\|v\|_2=1}}\|p(A)v\|_2=
\min_{p\in P_k(0)} \|p(A)\|_2;$$
see~\cite{FabJouKniMan96,Toh97} for the first examples, and~\cite{FabLieTic13} for additional examples and a detailed study of the mathematical properties of worst-case GMRES. In such cases every further bound that is derived using ideal GMRES, in particular the one given in \eqref{eqn:GMRESbd}, will be a significant overestimate of the actual GMRES residual norm.

If $A$ is normal, then we can choose an eigenvector matrix with $\kappa(X)=1$, and it can be shown that the bound \eqref{eqn:GMRESbd} is sharp in the sense that for each step~$k$ there exists an initial residual $r_0$ (depending on $A$ and $k$) so that equality holds; see the original proofs in~\cite{GreGur94,Jou94} and~\cite{LieTic14}. Thus, for a normal matrix $A$ the location of its eigenvalues determines the worst-case behavior of GMRES and, in this worst-case sense, gives an indication of the possible actual behavior. 

If $A$ is not normal, then $\kappa(X)$ can be very large, and the value of the upper bound \eqref{eqn:GMRESbd} may be far from the actual GMRES residual norm. Moreover, by separating the initial residual from the minimization problem in the upper bound \eqref{eqn:GMRESnorm}, we have lost all information about the relation between the particular given $A$ and $r_0$. This relationship between the particular $A$ and $b$, which is often inherited from the underlying problem to be solved, can be essential for the convergence behavior of GMRES, particularly for nonnormal matrices. We point out that the dependence of the convergence of CG on the initial residual, which is expressed in the exact formula \eqref{eqn:CGmin}, is in general much less pronounced than for GMRES. It is possible to derive closed formulas analogous to \eqref{eqn:CGmin} also for the GMRES residual norms, even for nondiagonalizable matrices. But because of the generally nonorthogonal eigenvectors (or principal vectors) in the GMRES context, such formulas are more involved and harder to interpret than \eqref{eqn:CGmin}. More details can be found, e.g., in~\cite{DuiMeuSadStr14} and~\cite[Sections~5.5--5.6]{MeuDui20}.

The eigenvalues of a general nonsingular matrix~$A$ may be anywhere in the complex plane, and hence estimating the value of the polynomial min-max approximation problem in \eqref{eqn:GMRESbd} can be very challenging. A quantitative bound can be given in the simple case that the eigenvalues of $A$ are contained in a disk centered at $c\in\C$ and with radius $\rho>0$, where $\rho<|c|$ is necessary so that zero is outside the disk. Taking the polynomial $(1-z/c)^k\in P_k(0)$ then shows that
\begin{equation}\label{eqn:GMRES-disk}
\min_{p\in P_k(0)}\max_{1\leq i\leq N} |p(\lambda_i)|\leq \left(\frac{\rho}{|c|}\right)^k;
\end{equation}
see, e.g.,~\cite[Section~6.11]{Saa03}. Thus, we can expect that GMRES converges quickly when $\kappa(X)$ is small, and the eigenvalues of $A$ are contained in a small disk that is far away from the origin in the complex plane. A survey of approaches for estimating the value of the min-max approximation problem beyond this special case is given in~\cite[Sections~5.7.2--5.7.3]{LieStr13}.

It needs to be stressed that the sharpness of the bound \eqref{eqn:GMRESbd} for normal matrices does not imply that GMRES converges faster for normal matrices, or that the (departure from) normality has an easy analyzable effect on the convergence of GMRES. In fact, it can be shown that GMRES may exhibit a complete stagnation even for unitary and hence normal matrices; see~\cite{GreStr94}. On the other hand, for a nonnormal matrix the location of the eigenvalues alone, and hence the value of the min-max approximation problem in \eqref{eqn:GMRESbd}, may not give relevant information about convergence behavior of GMRES. If $A$ is not diagonalizable, then its spectral decomposition does not exist, and an analogue of~\eqref{eqn:GMRESbd} based on the Jordan canonical form is of very limited use. As shown in~\cite{GrePtaStr96,GreStr94}, any nonincreasing convergence curve is possible for GMRES for a (in general, nonnormal) matrix $A$ having any prescribed set of eigenvalues. The work in~\cite{AriPtaStr98} gives a parametrization of the set of {\em all matrices and right-hand-sides} such that GMRES provides a given convergence curve while the matrix has a prescribed spectrum; see~\cite[Section~5.7.4]{LieStr13} for a summary.

\subsection{Any nonincreasing GMRES convergence curve is possible for any eigenvalues}\label{sec:GMRESevery}
\phantom{}
\smallskip
\begin{tcolorbox}
\emph{Main point: Eigenvalues alone are in general not sufficient for describing the GMRES convergence behavior.}
\end{tcolorbox}

\smallskip
\emph{Setup:} We follow~\cite[Section 2]{GrePtaStr96} for constructing a linear algebraic system $Ax=b$, where $A\in\R^{N\times N}$ has a prescribed set of eigenvalues $\lambda_1,\dots,\lambda_N\in\C$, so that the residual norms of GMRES applied to this system with $x_0=0$ are given by a prescribed nonincreasing sequence $f_0\geq f_1\geq \cdots \geq f_{N-1} >f_N=0$.

Define $g_j=\sqrt{(f_{j-1})^2 - (f_j)^2}$ for $j=1,\ldots, N$, let $V\in\R^{N\times N}$ be any orthogonal matrix, and let $b=V[g_1,\ldots, g_N]^T$. We then construct the polynomial
\[
(z-\lambda_1)(z-\lambda_2)\cdots(z-\lambda_N)=z^n - \sum_{j=0}^{N-1} \alpha_j z^j,
\]
and its companion matrix
\[
A^B = \begin{bmatrix}
0 & \cdots & 0 & \alpha_0 \\
1 & & 0 & \alpha_1 \\
& \ddots & \vdots& \vdots \\
& & 1 & \alpha_{N-1}
\end{bmatrix}.
\]
With $B=[b,v_1,\ldots, v_{N-1}]$, where $v_j$ denotes the $j$th column of $V$, we set $A = BA^BB^{-1}$.

We use $N=21$ and consider two scenarios: In the first we prescribe the eigenvalues $\lambda_1=\cdots=\lambda_{N}=1$, and the convergence curve $f_1=\cdots=f_{N-1}=1>f_{N}=0$. In the second we prescribe the eigenvalues $\lambda_j=j$ for $j=1,\dots,N$, and a convergence curve that starts at $f_1=1$, and then decreases every 4 steps through $10^{-2}$, $10^{-4}$, and $10^{-6}$, to $10^{-8}$. In both cases we take $V=I$, and apply MATLAB's {\tt gmres} function to $Ax=b$ with $x_0$. Both matrices in this example are highly ill conditioned and highly nonnormal. Computations with MATLAB's {\tt cond} and {\tt eig} functions yield
\begin{align*}
\kappa(A) &\approx 5.4\times 10^{11}\;\;\mbox{and}\;\; \kappa(X)\approx 7.6\times 10^{10} 
&&\mbox{($A$ with $\lambda_1=\cdots=\lambda_N=1$)}, \\
\kappa(A)&\approx 3.5\times 10^{14} \;\;\mbox{and}\;\; \kappa(X)\approx 1.5\times 10^{21}
&&\mbox{($A$ with $\lambda_j=j$)}.
\end{align*}

\begin{figure}
\begin{center}
\includegraphics[width=0.55\textwidth]{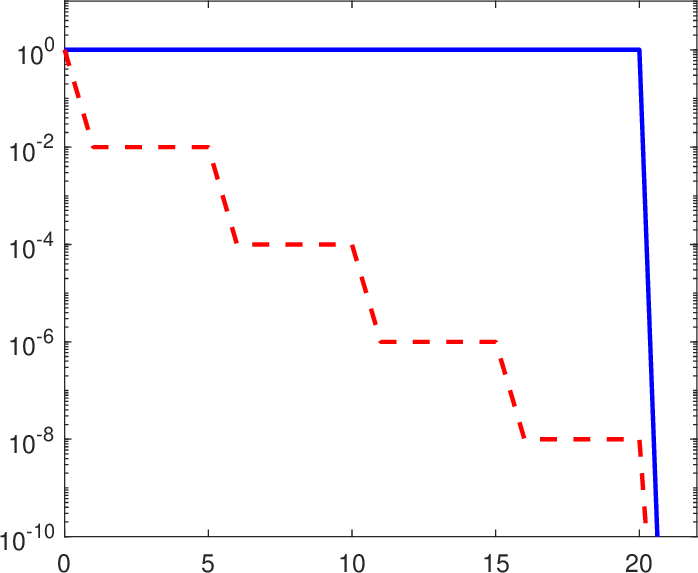}
\caption{Residual norms of GMRES applied to the two linear algebraic systems constructed with prescribed eigenvalues and convergence curves (solid blue for $\lambda_1=\cdots=\lambda_N=1$ and dashed red for $\lambda_j=j$).}\label{fig:gmresconv}
\end{center}
\end{figure}

\smallskip
\emph{Observations:} The computed GMRES residual for the two linear algebraic systems
are shown by the solid blue and dashed red curves in Figure \ref{fig:gmresconv}. As expected, they follow the prescribed convergence curves. 

\smallskip
\emph{Explanation:} 
The convergence curves in Figure \ref{fig:gmresconv} illustrate the proven theorems from~\cite{AriPtaStr98,GrePtaStr96,GreStr94}, and hence should not be surprising from a purely mathematical point of view. On the other hand, the curves illustrate the potentially intriguing behavior of GMRES for nonnormal matrices. 

The much slower convergence of GMRES (solid blue curve) occurs for the matrix having a single eigenvalue of algebraic multiplicity $N$. This clearly demonstrates that for general nonnormal matrices the goal of preconditioners for GMRES needs to be more than just ``clustering of the eigenvalues''.\footnote{Recall that this argument requires a thorough reconsideration also for CG; see Sections~\ref{sec:cgclust} and~\ref{sec:pcg}.} We also note that the slower convergence occurs for the somewhat better conditioned matrix.

We would like to mention a point that to our knowledge has not been sufficiently investigated yet: The paper~\cite{AriPtaStr98} contains a complete parameterization of the set of \emph{all pairs} $(A,b)$ for which GMRES with $x_0=0$ generates the prescribed convergence curve, while each of the matrices $A$ has the prescribed fixed eigenvalues. Thus, the parameterization contains not only artificially constructed examples, but also matrices and right hand sides that arise in practical applications.
Given the GMRES convergence curve for a particular (practical) linear algebraic system $\widehat{A}x=\widehat{b}$, we therefore know all pairs $(A,b)$ for which GMRES follows (with $x_0 = 0$) the practically observed behavior, while each of the matrices $A$ has the same eigenvalues as $\widehat{A}$. The set of the pairs $(A,b)$ is worth investigating in particular when $\widehat{A}$ is far from normal, but nevertheless its eigenvalues actually describe the GMRES convergence behavior.

\subsection{GMRES convergence for normal matrices}\label{sec:GMRESnormal}
\phantom{}
\smallskip
\begin{tcolorbox}
\emph{Main point: For a normal matrix the eigenvalues, and hence the bound \eqref{eqn:GMRESbd}, give a reasonable descriptive information of the convergence behavior of GMRES. If \eqref{eqn:GMRESbd} is estimated from above by \eqref{eqn:GMRES-disk}, then the resulting linear bound cannot capture a possible acceleration of convergence, which can occur due to adaptation of GMRES to the data.}
\end{tcolorbox}

\smallskip
\emph{Setup:} We consider normal (in fact, diagonal) matrices with eigenvalues that are rather uniformly distributed in certain disks. In order to construct such matrices we start with real nonsymmetric $(N\times N)$-matrices with normally distributed random entries that are generated with {\tt randn} in MATLAB. We then compute the eigenvalues of these matrices in MATLAB and form diagonal (and hence normal) matrices $D_N$ with these eigenvalues scaled by $1/\sqrt{N}$. We apply MATLAB's {\tt gmres} function with $x_0=0$ to linear algebraic systems with the matrices 
$$D_N,\quad 1.2 I+D_N, \quad\mbox{and}\quad 2I+0.5 D_N.$$ 
The right-hand sides $b$ are normalized random vectors, also generated with {\tt randn} in MATLAB. We use $N=100$, and we repeat the computation 100 times.

\begin{figure}
\begin{center}
\includegraphics[width=0.475\textwidth]{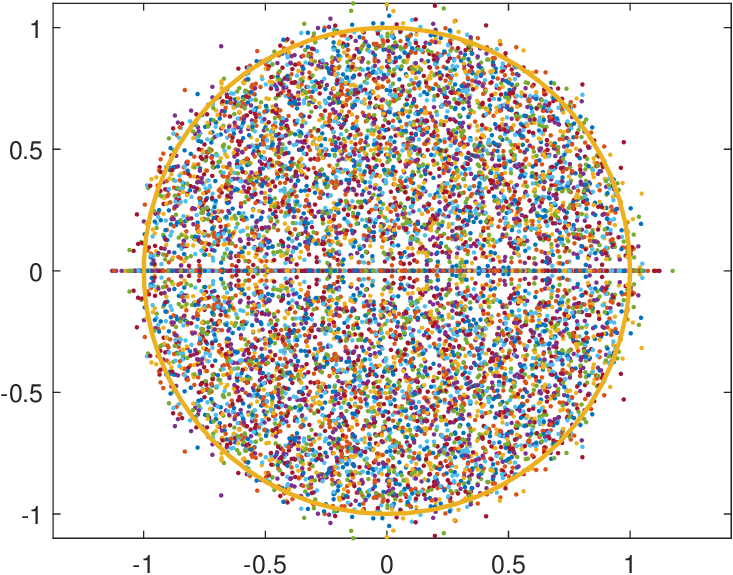}\hfill
\includegraphics[width=0.475\textwidth]{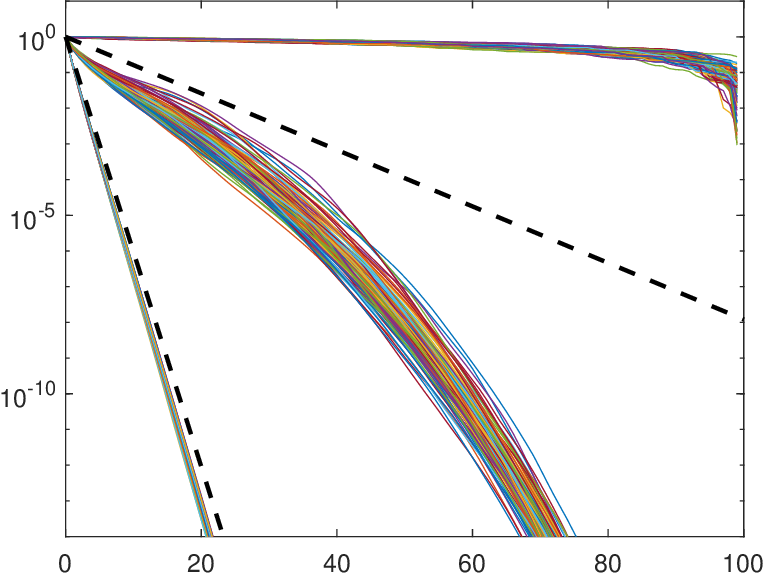}
\caption{Left: Eigenvalues of 100 matrices $D_N$ and the boundary of the unit disk. Right: Relative GMRES residual norms for linear algebraic systems with matrices $D_N$, $1.2 I+D_N$, and $2I+0.5 D_N$, and the upper bounds from~\eqref{eqn:GMRES-disk}.}\label{fig:GMRESrandn}
\end{center}
\end{figure}

\smallskip
\emph{Observations:} The left part of Figure~\ref{fig:GMRESrandn} shows the eigenvalues of the 100 matrices $D_N$ and the boundary of the unit disk. The analogous illustrations for the other matrices look similar, with the centers of the disks moved to $c=1.2$ and $c=2$, and the radius $\rho=0.5$ in the last case. We observe three distinctly different types of GMRES convergence behavior in the right part of Figure~\ref{fig:GMRESrandn}:
\begin{itemize}
\item For the matrices $D_N$ the method makes almost no progress until the very end of the iteration.
\item For the matrices $1.2 I+D_N$ the method converges faster than for the matrices $D_N$, the convergence accelerates in later iterations, and the convergence curves exhibit larger variations.
\item For the matrices $2I+0.5 D_N$ the method converges linearly and very quickly.
\end{itemize}
The two dashed lines in the right part of Figure~\ref{fig:GMRESrandn} show the values $(\rho/|c|)^k$ from the upper bound \eqref{eqn:GMRES-disk} corresponding to the disks with center $c=1.2$ and radius $\rho=1$, and with center $c=2$ and radius $\rho=0.5$. In the first case the bound captures the rate of convergence in the initial iterations, and in the second case it perfectly matches the convergence curve.

\smallskip
\emph{Explanation:} According to Girko's Circular Law, the eigenvalues of the matrices $D_N$ become uniformly distributed in the unit disk for $N\rightarrow\infty$; see the Introduction of~\cite{Ede97} for a summary of results in this context. This explains the eigenvalue distributions that are shown in the left part of Figure~\ref{fig:GMRESrandn}.

The worst-case behavior of GMRES for normal matrices is completely described by the values of the polynomial min-max approximation problem on the matrix eigenvalues, and these values also give an indication of the actual behavior; see the discussion of the bound~\eqref{eqn:GMRESbd}. Obviously, since the GMRES iteration polynomials are normalized at the origin, the worst-case GMRES (and usually also the actual GMRES) makes almost no progress when the eigenvalues are spread through the unit disk, as for the matrices $D_N$.

Most of the normal matrices $1.2I+D_N$ have their eigenvalues in the disk centered at $c=1.2$ and with radius $\rho=1$. For these matrices the worst-case bound~\eqref{eqn:GMRES-disk} gives a reasonable description of the actual GMRES convergence particularly in the initial phase of the iteration. However, since the bound~\eqref{eqn:GMRES-disk} is linear, its combination with~\eqref{eqn:GMRESbd} cannot describe the (nonlinear) GMRES acceleration observable in later iterations; cf. Section~\ref{sec:CG} for similar observations in the context of the $\kappa(A)$-bound for the CG method. In comparison with CG, understanding the convergence acceleration of GMRES is challenging (and largely open) even for normal matrices. The polynomial min-max approximation problem on a (finite) discrete set in the complex plane is much more  difficult to handle than the same problem on a (finite) discrete set on the real line.  

Finally, most of the normal matrices $2I+0.5 D_N$ have their eigenvalues located in the disk centered at $c=2$ and with radius $\rho=0.5$. For these matrices GMRES converges so fast that no convergence acceleration becomes visible.

\subsection{MGS-GMRES is normwise backward stable}\label{sec:GMRESorth}
\phantom{}
\smallskip
\begin{tcolorbox}
\emph{Main point: In MGS-GMRES, complete loss of orthogonality of the computed Krylov subspace basis means convergence of the normwise relative backward error to the maximal attainable accuracy, and vice versa. Therefore, the MGS-GMRES method is normwise backward stable.}
\end{tcolorbox}

\smallskip
\emph{Setup:} The experiment investigates the relation between the loss of orthogonality of the Krylov subspace basis vectors computed in finite precision arithmetic using the modified Gram-Schmidt (MGS) variant of the Arnoldi algorithm (measured by $\|I-V_k^TV_k\|_F$), and the convergence of the normwise relative backward error  $\|b-Ax_k\|_2/(\|b\|_2+\|A\|_2\|x_k\|_2)$ in the corresponding MGS-GMRES method.\footnote{For a discussion of the practical relevance of the normwise relative backward error when solving linear algebraic systems we refer to any good textbook on numerical linear algebra. Its use as a stopping criterion for GMRES is discussed in Section~\ref{sec:GMRESstopping}.}

We consider the example matrices {\tt fs1836} and {\tt sherman2} from Matrix Market\footnote{\url{https://math.nist.gov/MatrixMarket/}}. The matrix {\tt fs1836} is of size $N=183$ and has a condition number of approximately $1.0\times 10^7$. The matrix {\tt sherman2}is of size $N=1000$ and has a condition number of approximately $2.4\times 10^7$. Both matrices are diagonalzable, and the condition number of their eigenvector matrices computed by MATLAB are approximately $1.7\times 10^{11}$ and $9.6\times 10^{11}$, respectively. We use the right-hand side $b=Ax$, where $x=[1,\dots,1]^T/\sqrt{N}$. Since MATLAB's {\tt gmres} function is based on the Householder variant of the Arnoldi algorithm, we use in this experiment our own MGS-GMRES implementation, starting with $x_0=0$.

\begin{figure}
\begin{center}
\includegraphics[width=0.475\textwidth]{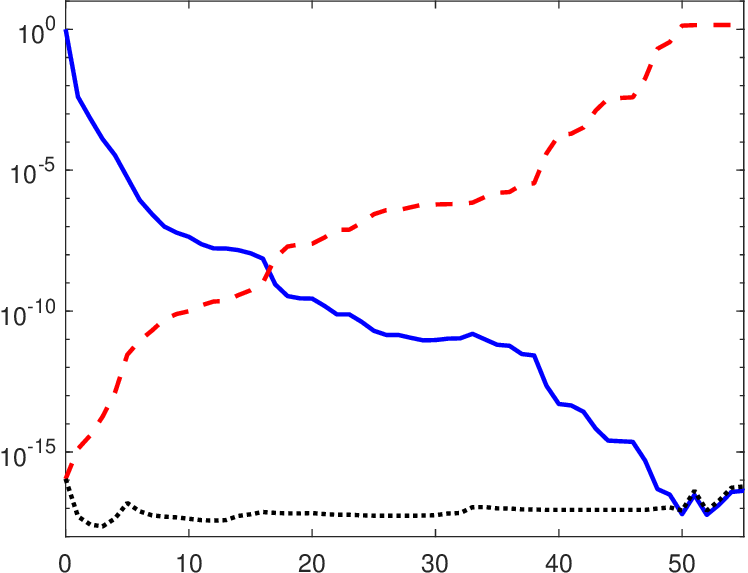} \hfill
\includegraphics[width=0.475\textwidth]{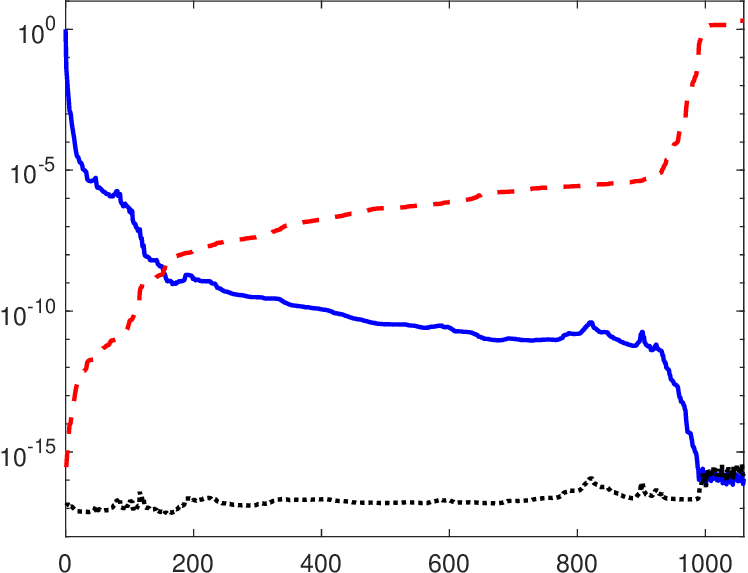}
\caption{MGS-GMRES normwise relative backward errors (solid blue) and loss of orthogonality (dashed red), and the product of the two (dotted) for linear algebraic systems with the matrices {\tt fs1836} (left) and {\tt sherman2} (right).}\label{fig:bwerr_loss}
\end{center}
\end{figure}

\smallskip
\emph{Observations:} As shown in Figure~\ref{fig:bwerr_loss}, throughout the iteration the product
$$\|I-V_k^TV_k\|_F\,\times\,\frac{\|b-Ax_k\|_2}{\|b\|_2+\|A\|_2\|x_k\|_2}$$
is almost constant, and close to the machine precision (approximately $10^{-16}$). The orthogonality of the basis vectors is completely lost only when the normwise backward error of the MGS-GMRES iteration has reached its maximal attainable accuracy level. We point out that this is a significant difference to the finite precision behavior of CG and other methods based on short recurrences. Then, apart from very special situations (see Section~\ref{sec:cgpoisson}), not only a loss of orthogonality but also a loss of rank in the computed subspace may occur, which leads to a delay of convergence; see Section~\ref{sec:fpcg} and the detailed explanations in, e.g., \cite[Section~5.9]{LieStr13} and \cite{MeuStr06}.

\smallskip
\emph{Explanation:} The full explanation of the phenomenon observed in Figure~\ref{fig:bwerr_loss} is based on a detailed analysis of rounding errors that occur in MGS-GMRES; see the papers~\cite{DrkGreRozStr95,PaiStr02,PaiRozStr06} or the summary in~\cite[Section~5.10]{LieStr13}. In this analysis it is shown (under some technical assumptions, e.g., that $A$ is not too close to being singular) that the loss of orthogonality in the Krylov subspace basis $V_k$ computed in finite precision arithmetic using the MGS variant of the Arnoldi algorithm is essentially controlled by the condition number $\kappa([\gamma v_1,AV_kD_k])$, where $\gamma\in\R$ and $D_k\in\R^{k\times k}$ are suitable chosen scalings. Since $r_k=r_0-AV_kt_k=\|r_0\|v_1-AV_kt_k$, the conditioning of the matrix $[\gamma v_1,AV_kD_k]$ can be related to the residual norm $\|b-Ax_k\|_2$, and in a second step also to the normwise relative backward error $\|b-Ax_k\|_2/(\|b\|_2+\|A\|_2\|x_k\|_2)$. This yields a rigorous proof of the numerically observed behavior of MGS-GMRES. It is worth noting that working on this challenge led to revisiting the theoretical foundations of (scaled) total least squares problems; see, e.g., \cite{PaiStr02a}.

Note that the (more costly) GMRES implementation based on the Householder variant of the Arnoldi algorithm, which is used, e.g., in MATLAB, is also normwise backward stable~\cite{DrkGreRozStr95}. In the classical Gram-Schmidt (CGS) variant, however, the orthogonality is lost too quickly to guarantee backward stability of the corresponding CGS-GMRES implementation; see~\cite{GirLanRozEsh05} for a rounding error analysis of the CGS orthogonalization algorithm.

\subsection{GMRES convergence for approximately computed preconditioning}\label{sec:pGMRES}
\phantom{}
\smallskip
\begin{tcolorbox}
\emph{Main point: 
In practical computations, preconditioning can not be performed exactly due to rounding errors. More substantially, theoretical preconditioning often has to be intentionally approximated with a rather relaxed accuracy in order to prevent prohibitive cost. Therefore theoretical results, which hold for exact preconditioners, must be used with caution when applied to practical heuristics. Moreover, when using an (approximate) preconditioner, the norm of the true residual of the final computed approximation should be checked.}
\end{tcolorbox}

\emph{Setup:} We set up a linear algebraic system ${\cal A}x=b$, where 
$${\cal A}=\begin{bmatrix}A & B^T\\B & 0\end{bmatrix}$$
comes from a discretization of a Navier-Stokes model problem in IFISS 3.6~\cite{ElmRamSil14}, and $b=[1,\dots,1]^T/\sqrt{N}$. We run the {\tt navier\_testproblem} with the (mostly default) parameters.\footnote{Cavity; regularized; 16x16 grid; uniform grid; Q1-P0; viscosity: 1/100; hybrid; Picard: 1; Newton: 1; nonlinear tolerance: 1.1*eps; uniform streamlines.} The matrix $A\in\R^{n\times n}$ is nonsymmetric, and $B\in\R^{m\times n}$ has full rank $m$. For our chosen model problem parameters we have $n=578$ and $m=256$. 

We consider the block diagonal preconditioner
$${\cal P}=\begin{bmatrix}A & 0\\0 & S\end{bmatrix},$$
where $S=BA^{-1}B^T$ is the Schur complement, and we are interested in the behavior of GMRES for the preconditioned system ${\cal P}^{-1}{\cal A}x={\cal P}^{-1}b$. The top block of the preconditioner ${\cal P}$ is given by the explicitly known matrix $A$. In order to simulate exact preconditioning, we compute the bottom block $S$ using MATLAB's backslash operator for the inversion of~$A$. We use the matrix resulting from this computation as the ``exact'' preconditioner ${\cal P}$. We then apply our own implementation of MGS-GMRES with $x_0=0$ to ${\cal A}x=b$ and to the exactly preconditioned system ${\cal P}^{-1}{\cal A}x={\cal P}^{-1}b$, where we compute ${\cal P}^{-1}{\cal A}$ and ${\cal P}^{-1}b$ again using MATLAB's backslash operator. 
 
In order to simulate the effects of intentionally inexact preconditioning, we illustrate how the GMRES convergence behavior changes when instead of explicitly computing ${\cal P}^{-1}{\cal A}$ and ${\cal P}^{-1}b$, we apply inner GMRES iterations for approximating solutions of linear algebraic systems with ${\cal P}$ in every step of the outer GMRES iteration (the so-called inner-outer iterations). 
We stop the inner iterations (starting with the initial vector $x_0=0$) when the relative residual norm reaches the respective tolerances $10^{-8}$, $10^{-4}$, and $10^{-2}$. In this way we obtain three GMRES convergence curves for approximately preconditioned systems that simulate the decreasing accuracy of performing preconditioning.

We stress that in this example we only illustrate that inexactness in performing preconditioning can cause a substantial departure of the GMRES convergence behavior from the behavior guaranteed by theoretical results which assume that the preconditioning is performed exactly. By no means we aim to suggest practical preconditioning strategies, or to study inner-outer iterations with particular preconditioners. (Note that in practice a Krylov subspace method combined with multigrid preconditioning is often a viable approach.)

\begin{figure}[h]
\begin{center}
\includegraphics[width=0.475\textwidth]{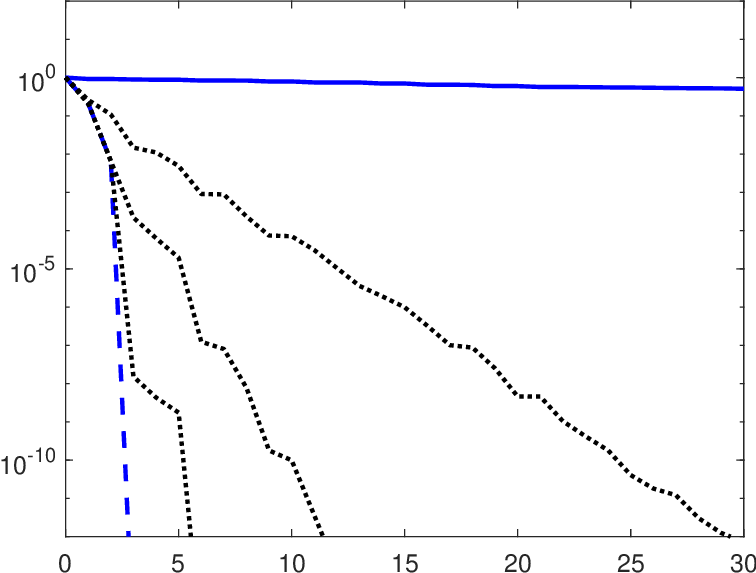}\hfill
\includegraphics[width=0.475\textwidth]{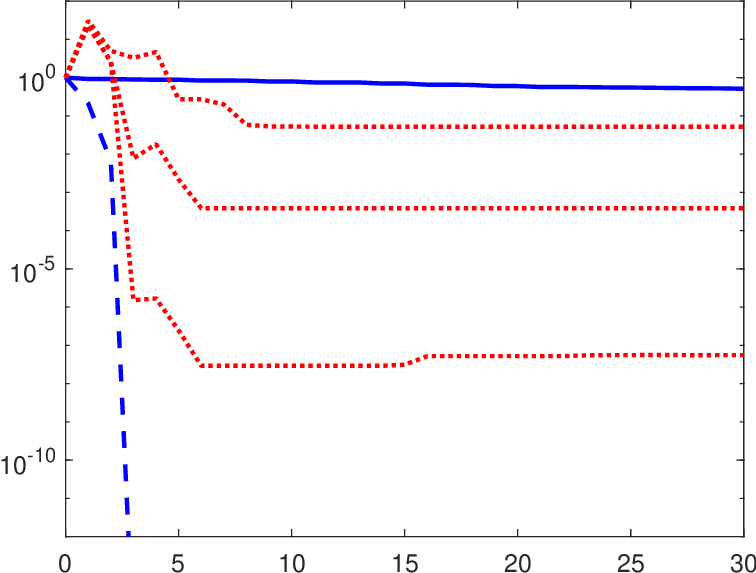}
\caption{Both parts show the relative resdiual norms $\|r_k\|_2\|/\|r_0\|_2$ of GMRES for the unpreconditioned Navier-Stokes system (solid blue) and for the exactly preconditioned system (dashed blue). The left part also shows the preconditioned relative residual norms  $\|{\cal P}^{-1}(b-{\cal A}x_k)\|_2/\|{\cal P}^{-1}r_0\|_2$ of GMRES for the approximately preconditioned systems with inner iteration tolerances $10^{-8}$, $10^{-4}$, and $10^{-2}$ (dotted black). Relaxing the tolerance slows down the GMRES convergence. The right part also shows the relative residiual norms $\|r_k\|_2\|/\|r_0\|_2$ of GMRES for the approximately preconditioned systems with inner iteration tolerances $10^{-8}$, $10^{-4}$, and $10^{-2}$ (dotted red). Relaxing the tolerance worsens the final accuracy.}\label{fig:gmres_navier}
\end{center}
\end{figure}

\smallskip
\emph{Observations:} In the left part of Figure~\ref{fig:gmres_navier} we show the relative residual norms $\|r_k\|_2/\|r_0\|_2$ of GMRES applied to ${\cal A}x=b$ (solid blue) and to the exactly pre\-cond\-i\-tioned sys\-tem ${\cal P}^{-1}{\cal A}x={\cal P}^{-1}b$ (dashed blue), as well as the \emph{preconditioned relative residual norms} $\|{\cal P}^{-1}(b-{\cal A}x_k)\|_2/\|{\cal P}^{-1}r_0\|_2$ of GMRES applied to the approximately preconditioned systems with the three tolerances for the inner solves (dotted black). We see that for the unpreconditioned system GMRES makes virtually no progress, and that for exact preconditioner it converges in three steps. The speed of convergence measured by the preconditioned relative residual norm slows down when we relax the accuracy of applying the preconditioner.

The solid blue and dashed blue curves in right part of Figure~\ref{fig:gmres_navier} are the same as in the left part. The other three curves (dotted red) show the actual relative residual norms $\|r_k\|_2/\|r_0\|_2$ when GMRES is applied to the approximately preconditioned systems with the three different tolerances for the inner iterations as above. We observe that these curves are not monotonically decreasing, and that in each case the maximal attainable accuracy is approximately on the accuracy level of the inner iteration.

\smallskip
\emph{Explanation:} As shown in the instructive and widely cited paper~\cite{MurGolWat00}, the minimal polynomial of the (nonsingular) preconditioned matrix ${\cal P}^{-1}{\cal A}$ is given by $(z-1)(z^2-z-1)$, and hence this matrix has the three distinct eigenvalues $1$ and $(1\pm \sqrt{5})/2$. Thus, in exact computations, GMRES applied to the exactly preconditioned system ${\cal P}^{-1}{\cal A}x={\cal P}^{-1}b$ converges to the exact solution in at most three steps; see the dashed blue lines in Figure~\ref{fig:gmres_navier}. As clearly pointed out in~\cite{MurGolWat00}, the degree of the minimal polynomial is essential for this property. For mathematical convergence of GMRES to the exact solution in at most~$k$ steps, it is not sufficient that~$A$ has only~$k$ distinct eigenvalues. For this the matrix additionally must be diagonalizable; cf. the example with $\lambda_1=\cdots=\lambda_N=1$ in Section~\ref{sec:GMRESevery}.

In practical computations we usually do not form the preconditioner explicitly, nor do we form the preconditioned matrix. Instead, we use inner iterations for the linear algebraic systems with ${\cal P}$. In addition, these inner iterations are usually based only on an approximation of ${\cal P}$, which is obtained, for example, by approximating a Schur complement. Clearly, the exact mathematical properties of the preconditioned matrix no longer hold for the approximate preconditioning, and therefore cannot be used for a rigorous analysis of the GMRES convergence behavior for the approximately preconditioned systems. The inexactness of the preconditioner can negatively affect the convergence behavior, both in terms of the rate of convergence and the maximal attainable accuracy, as illustrated in Figure~\ref{fig:gmres_navier}. 

Many publications, in particular in the context of linear algebraic systems in saddle point form, give bounds on the eigenvalues of inexactly preconditioned matrices, with the goal to show that under inexact preconditioning the few eigenvalues with large multiplicity are replaced by a few clusters of eigenvalues. Related to the clustering heuristics, it is worth considering the next two points. First, for highly nonnormal matrices even small perturbations can make these clusters very large. Second, in order to be meaningful for understanding the GMRES convergence behavior, such a perturbation argument for the eigenvalues needs to be complemented by additional arguments, since there is no guarantee that an approximately preconditioned matrix remains diagonalizable, and the convergence behavior for nonnormal matrices does not depend on the eigenvalues only; see Section \ref{sec:GMRESevery}. In particular cases, where the theoretical results on the distribution of eigenvalues indeed explain the results of practical computations, one should try to identify the particular properties (in the words of Lanczos ``the inner nature'') of the problem that makes this possible, starting with the mathematical model of the computationally investigated phenomenon. (See also the experiment in the next Section~\ref{sec:GMRESminpoly}.)

It is important to note that when applied to a preconditioned system, whether the preconditioner is computed exactly or not, GMRES minimizes the Euclidean norm of the preconditioned residual norm in every step. And since $\|{\cal P}^{-1}(b-{\cal A}x_k)\|_2=\|b-{\cal A}x_k\|_{({\cal P}{\cal P}^*)^{-1}}$, this can be interpreted as the minimization of the residual in a norm that depends on the preconditioner. As pointed out in~\cite[p.~193]{ElmSilWat14} in the context of the MINRES method, ``one must be careful not to select a preconditioner that wildly distorts this norm". Such a ``distortion of the norm'' can mean that the preconditioned residual norms can be significantly different from the residual norms $\|r_k\|_2 = \|b- Ax_k\|_2$ measuring convergence behavior with respect to the unpreconditioned system. There is no guarantee that the latter are monotonically decreasing, since preconditioned GMRES iterations are only optimal for the preconditioned system. 

Some implementations of preconditioned GMRES, for example MATLAB's {\tt gmres} function, have only the preconditioned residual norms as their standard output, and many publications containing preconditioned GMRES computations only report the behavior of the preconditioned residual norms. The right part of Figure~\ref{fig:gmres_navier} indicates that whenever inexact preconditioning is used, it should be accompanied by a basic analysis of the actual residual norms used in stopping the iterations, or at least by computing $\|r_k\|_2$ at the end of the iteration in order to check the attained accuracy level with respect to the given (unpreconditioned) system.

\subsection{GMRES and the minimal polynomial}\label{sec:GMRESminpoly}
\phantom{}
\smallskip
\begin{tcolorbox}
\emph{Main point:
Whenever the behavior of Krylov subspace methods is linked with approximation of the minimal polynomial of the matrix, we must rigorously specify the way in which the accuracy of such an approximation is going to be measured. In particular, the roots of the GMRES iteration polynomial do not need to match those of the minimal polynomial.}
\end{tcolorbox}

\smallskip
\emph{Setup:} 
We consider the \emph{Grcar matrix} (see, e.g.,~\cite[Example~3]{Tre92})
$$A=\begin{bmatrix}
 1 & 1 & 1 & 1 & & & \\
-1 & 1 & 1 & 1 & 1 & &\\
& \ddots & \ddots & \ddots & \ddots & \ddots &\\
& & \ddots & \ddots & \ddots & \ddots & 1\\
& & & \ddots & \ddots & \ddots & 1\\
& & & & \ddots & \ddots & 1\\
& & & & & -1 & 1
\end{bmatrix}\in\R^{500\times 500}$$
This matrix is diagonalizable and well conditioned with $\kappa(A)\approx 3.63$. However, as usual for nonsymmetric Toeplitz matrices (see, e.g.,~\cite{ReiTre92}), the eigenvectors of~$A$ are very ill conditioned, so that~$A$ is highly nonnormal. Using MATLAB's {\tt eig} and {\tt cond} functions yields an eigenvector matrix~$X$ with $\kappa(X)\approx 7.2\times 10^{38}$. We use $b=[1,\dots, 1]^T/\sqrt{N}$ and apply our own implementation of MGS-GMRES to $Ax=b$ with $x_0=0$.

\begin{figure}[t!]
\begin{center}
\includegraphics[width=0.475\textwidth]{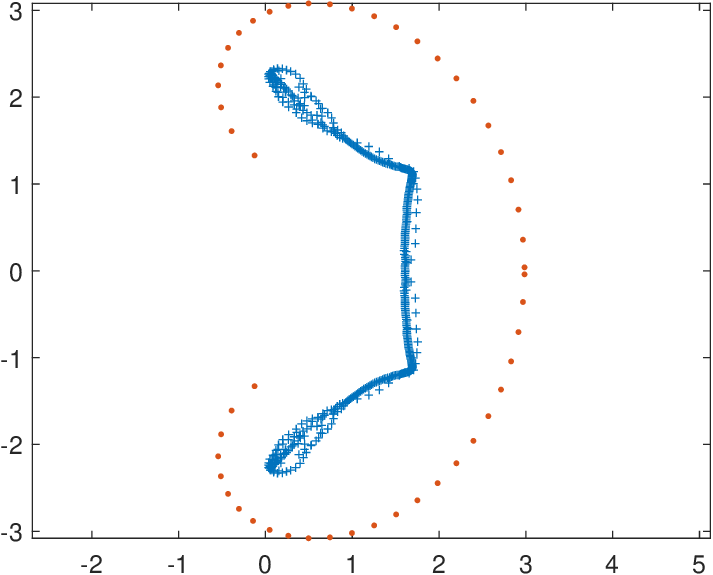}\hfill
\includegraphics[width=0.475\textwidth]{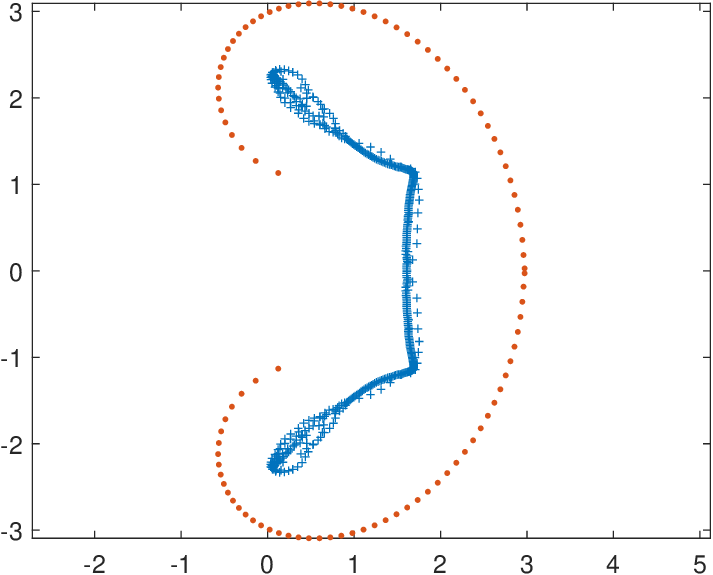}\\
\medskip
\includegraphics[width=0.475\textwidth]{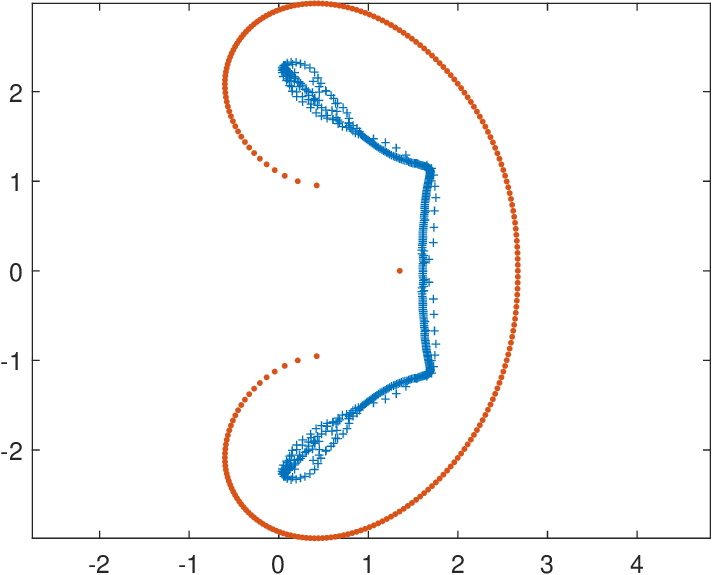}\hfill
\includegraphics[width=0.475\textwidth]{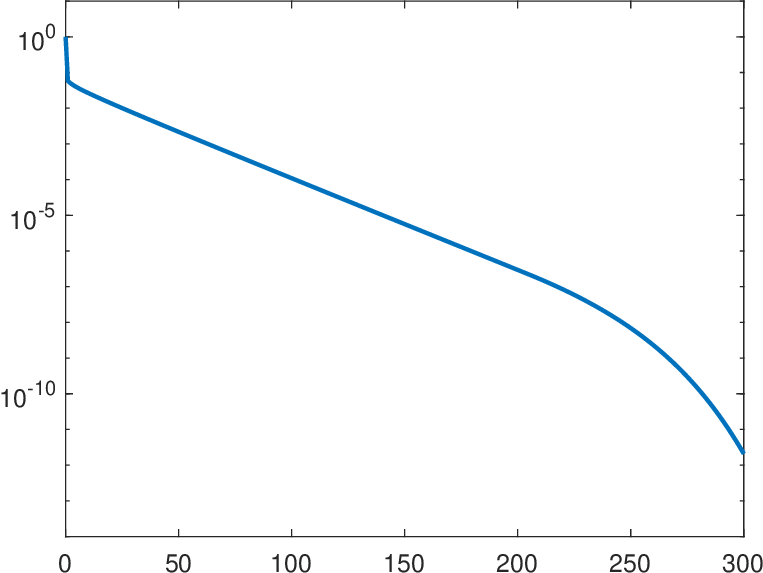}
\caption{Eigenvalues of the $500\times 500$ Grcar matrix $A$ computed by MATLAB (blue pluses) and the harmonic Ritz values (red dots) at GMRES steps 50, 100 and 200 (top left, top right, and bottom left), and relative GMRES residual norms for $Ax=b$ with $b=[1,\dots,1]^T/\sqrt{N}$ and $x_0=0$ (bottom right).}\label{fig:GMRESgrcar}
\end{center}
\end{figure}

\smallskip
\emph{Observations:}
The relative GMRES residual norms are shown in the bottom right part of Figure~\ref{fig:GMRESgrcar}. In the first iteration, the relative residual norm drops from $1.0$ to approximately $0.05$, and then the relative residual norms decrease almost linearly for the following approximately 250 steps. In the other three parts of Figure~\ref{fig:GMRESgrcar} the (blue) pluses show the (approximate) eigenvalues of~$A$ computed by MATLAB's {\tt eig} function. (Note that because of the severe ill-conditioning of the eigenvalue problem, this computation is, for the eigenvalues with large imaginary parts, affected by rounding errors. A thorough discussion of the spectrum of the Grcar matrix and of its approximation can be found in \cite{TreEmbBook05}.) The (red) dots show the roots of the GMRES polynomials, also called the harmonic Ritz values (see, e.g.,~\cite[Section~5.7.1]{LieStr13} for mathematical characterizations), at iterations $50, 100$, and $200$ (top left, top right, and bottom left, respectively). During the iteration, these roots fill up more and more of the same curve that ``surrounds'' the eigenvalues of~$A$, but overall they fail to move any closer towards the eigenvalues, and hence towards the roots of the minimal polynomial of $A$.

\emph{Explanation:}
In the CG method, the polynomial of degree $k$ providing the minimum in \eqref{eqn:CGmin} indeed approximates the minimal polynomial of the matrix $A$ (assuming that all $\eta_i$ are nonzero) {\em in the sense of solving the simplified Stietjes moment problem} (see~\cite{PozStr19}) or, equivalently, {\em in the sense of determining the nodes of the associated $k$-point Gauss quadrature}; see, e.g.,~\cite[Section~3.5]{LieStr13} and~\cite[Section~5.2]{MSB15}. As demonstrated in the examples in Section~\ref{sec:CG}, this does not mean that there exists a \emph{simple revealing relationship} between the eigenvalues of $A$ and the Ritz values $\theta_\ell^{(k)}$ in \eqref{eqn:CG_Ritz_values}, not even for tightly clustered eigenvalues. Of course, there exist many beautiful properties due to the underlying orthogonal polynomials, and many results about convergence of Ritz values to the eigenvalues. But this is not the same as defining a meaningful measure in relation to approximation of the minimal polynomial.

For GMRES the situation is even more complicated, since analogues of \eqref{eqn:CGmin} and \eqref{eqn:CG_Ritz_values} are not available for the GMRES polynomial. There exist generalizations of the Gauss quadrature through the matrix formulation of the Vorobyev moment problem (see, e.g, \cite[Section~3.7.4]{LieStr13} and \cite{Str09}) which apply to the Arnoldi algorithm for approximating eigenvalues and to the closely related FOM method for solving linear systems; see, e.g., \cite[Section 6.4]{Saa03}. In this sense, and {\em only in this sense}, the GMRES polynomial at step $k$ approximates the minimal polynomial of~$A$. However, there is no apparent way how these very involved relationships can give meaningful insights into the location of the $k$ harmonic Ritz values in relation to the roots of the minimal polynomial of~$A$. And indeed, the harmonic Ritz values may remain far from the eigenvalues of $A$ throughout the iteration, although the GMRES residual norms decrease reasonably quickly. This important point is illustrated in Figure~\ref{fig:GMRESgrcar}. 

In many publications it is claimed that the main idea behind Krylov subspace methods is to approximate the minimal polynomial of~$A$. The arguments given above explain why this interpretation is misleading.\footnote{The interpretation may be motivated by the widespread use of Krylov subspaces for solving eigenvalue problems with~$A$, and it somewhat resonates with the original paper of Krylov published in 1931~\cite{Kry31}, which deals with computing eigenvalues.} In the context of solving large linear algebraic systems, it can can lead to further misconceptions for two main reasons:

\emph{Reason 1:}  The number of iterations performed in practice is typically many orders of magnitude smaller than the degree of the minimal polynomial in question.
This remains true also in cases where an ``ideal'' preconditioner guarantees that, in theory, the degree of the minimal polynomial of the exactly preconditioned matrix is very small. As discussed in Section~\ref{sec:pGMRES}, in practice we do not precondition exactly, which results in replacing the few eigenvalues (with large multiplicities) by a few clusters of eigenvalues. The argument continues that we can utilize these clusters for approximating the minimal polynomial of the inexactly preconditioned matrix. However, even if such a matrix is diagonalizable (which is in general not obvious), its minimal polynomial has a very large degree (usually equal to the size of the problem). Approximating clusters of eigenvalues by single roots of the iteration polynomial does not work, apart from particular cases, regardless of how tight the clusters are. This is explained in Section~\ref{sec:cgclust} for a much easier case. It may work under some specific circumstances and restrictions that have to be clearly stated whenever the argument is used.

\emph{Reason 2:} The mathematical term ``approximation" should be used only with a precise description of the measure that is used for evaluating the accuracy of the approximation. In the context of Krylov subspace methods the flaw is not in using the term ``approximation" in relation to the minimal polynomial of the system matrix. The flaw is either in not specifying any measure at all, or in a vague association of such a measure with the locations of the roots of the iteration polynomials.

\subsection{Residual versus normwise backward error, and stopping criteria for GMRES}\label{sec:GMRESstopping}
\phantom{}
\smallskip
\begin{tcolorbox}
\emph{Main point: The residual $2$-norm, which is mimimized by GMRES, is commonly used as a stopping criterion for GMRES. The normwise relative backward error represents an important and practically relevant alternative.}
\end{tcolorbox}

\smallskip
\emph{Setup:} 
We consider linear algebraic systems with the matrices {\tt fs1836} and {\tt sherman2} from Matrix Market (see Section~\ref{sec:GMRESorth}), and right-hand sides $b=[1,\dots,1]^T/\sqrt{N}$. We apply our own MGS-GMRES implementation to $Ax=b$, starting with $x_0=0$, and compute the relative residual norms $\|r_k\|_2/\|b\|_2$ and the normwise relative backward errors $\|r_k\|_2/(\|b\|_2+\|A\|_2\|x_k\|_2)$. For the computation of $\|A\|_2$ we use MATLAB's {\tt norm} function.

\begin{figure}
\begin{center}
\includegraphics[width=0.475\textwidth]{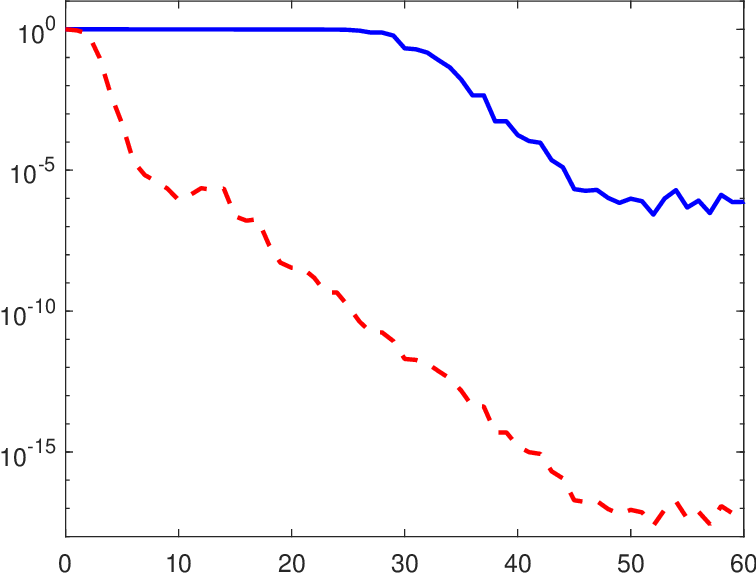} \hfill
\includegraphics[width=0.475\textwidth]{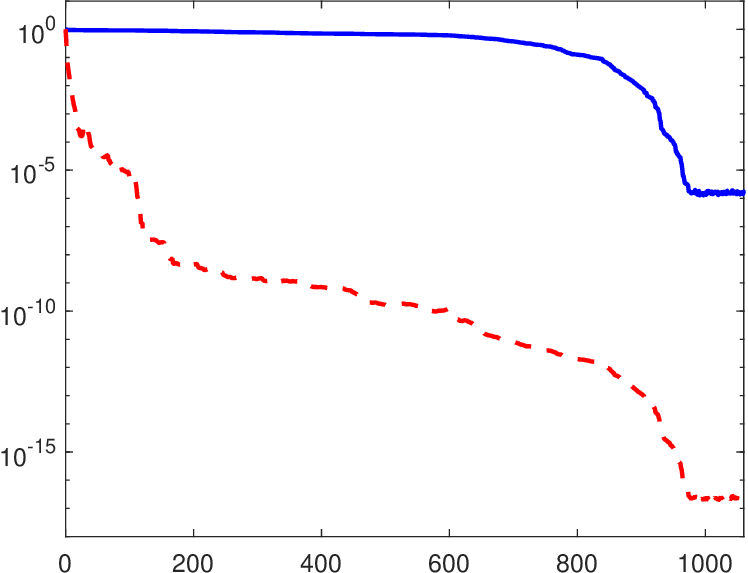}
\caption{Residual $2$-norms (solid blue) and normwise relative backward errors (dashed red) of GMRES for linear algebraic systems with {\tt fs1836} (left) and {\tt sherman2} (right).} \label{fig:gmres_stop}
\end{center}
\end{figure}

\smallskip
\emph{Observations:} We can see in Figure~\ref{fig:gmres_stop} that the behavior of the relative residual norm and the normwise relative backward error may be very different. In both cases the relative residual norm initially stagnates and eventually reaches a level of approximately $10^{-6}$. The normwise relative backward error, on the other hand, decreases quickly from the start of the iteration and reaches the machine precision level.  (Note that the backward error curves for {\tt fs1836} and {\tt sherman2} look a bit different than in Figure~\ref{fig:bwerr_loss}, since here we use different right-hand sides.) 

\smallskip
\emph{Explanation:} Given any approximation $x_k$ of the exact solution $x$, it is useful to ask which linear algebraic system is solved exactly, i.e., for which $\Delta A$ and $\Delta b$ we have 
$$(A+\Delta A)x_k=b+\Delta b;$$
see Wilkinson~\cite{WilBook65,WilBook63}. In many applications both $A$ and $b$ result from approximations, measurements or discretizations. From an application point of view it is therefore often sufficient to know that $x_k$ solves a nearby system, i.e., a system where the perturbations $\Delta A$ and $\Delta b$ are small. The perturbation restricted to the right-hand side gives
$$Ax_k=b+\Delta b,\quad\mbox{where}\quad \Delta b=-r_k=Ax_k-b.$$
Thus, if $\|r_k\|_2$ is small, then $x_k$ solves a system with the same matrix and a slightly perturbed right-hand side. The relative size of the perturbation is $\|r_k\|_2/\|b\|_2$. For $x_0=0$ this is the relative residual norm, which commonly is used as a stopping criterion for GMRES. For $x_0\neq 0$, however, the relative residual norm $\|r_k\|_2/\|r_0\|_2$ does not have a backward error interpretation, and when $\|r_0\|_2\gg \|b\|_2$ its use as a measure of convergence requires caution; see Section~\ref{sec:GMRESnonzero} below. 

As shown by Rigal and Gaches~\cite{RigGac67}, the normwise relative backward error (plotted in Figure~\ref{fig:gmres_stop}) gives the size of the smallest possible perturbations $\Delta A_{\min}$ and $\Delta b_{\min}$ for which $x_k$ solves the perturbed system exactly. More precisely, defining  
$$\beta(x_k):=\min\left\{\beta\;:\;(A+\Delta A)x_k=b+\Delta b,\;\|\Delta A\|_2\leq \beta\|A\|_2,\;
\|\Delta b\|_2\leq \beta\|b\|_2\right\}$$
we have
$$\beta(x_k)=\frac{\|r_k\|_2}{\|b\|_2+\|A\|_2\|x_k\|_2}=
\frac{\|\Delta A_{\min}\|_2}{\|A\|_2}=\frac{\|\Delta b_{\min}\|_2}{\|b\|_2}.$$
Consequently, the backward error analysis yields a valid argument for using the normwise relative backward error instead of the relative residual norm as a stopping criterion for GMRES. 

In our example, the relative residual norms in Figure~\ref{fig:gmres_stop} eventually stagnate close to $10^{-6}$. This shows that in both cases $x_k$ eventually solves a system with the matrix $A$ and a rather large perturbation of the right-hand side. But since the normwise relative backward errors decreases to the level of the machine precision, we know that $x_k$ eventually solves a system with very small perturbations of $A$ and $b$.

\subsection{GMRES convergence for different right-hand sides}\label{sec:GMRESinitres}
\phantom{}
\smallskip
\begin{tcolorbox}
\emph{Main point: The convergence behavior of GMRES can depend strongly on the right-hand side, and hence convergence analysis based only on the matrix may not be descriptive.}
\end{tcolorbox}

\smallskip
\emph{Setup:}
The first example is a variation of the \emph{Frank matrix}, which is a test matrix of upper Hessenberg form generated by {\tt gallery('frank',N,N)} in MATLAB. We ``flip'' this matrix and consider
$$F_N=\begin{bmatrix}
    N & N-1  &   N-2  &   \cdots  &   1\\
    N-1 & N-1  &   N-2  &  \cdots  &   1\\
     & N-2  & N-2    &  \cdots   &   1\\
     &   &   \ddots  &   \ddots  &   \vdots\\
     &   &     &   1  &   1
\end{bmatrix}\in \R^{N\times N}.$$
We use $N=16$, which yields a highly nonnormal matrix. Computations with MATLAB's {\tt eig} and {\tt cond} functions yield an eigenvector matrix $X$ with $\kappa(X)\approx 1.1\times 10^{13}$. We apply our implementation of MGS-GMRES with $x_0=0$ to the systems $F_{N}x=b^{(j)}$, where $b^{(1)}=[1,\dots,1]^T/\sqrt{N}$, and $b^{(2)}$ is a normalized random vector with normally distributed entries, generated using {\tt randn} in MATLAB.

The second example is a discretized convection-diffusion problem that was studied in~\cite{LieStr05}; see also~\cite[Section~5.7.5]{LieStr13}. Here the SUPG discretization with stabilization parameter $\delta$ of the problem
$$-\nu (u_{xx}+u_{yy})+u_y=0\;\;\mbox{in}\;\;\Omega=(0,1)\times (0,1),\quad
u=g\;\;\mbox{on}\;\;\partial\Omega,$$
leads to linear algebraic systems $Ax=b$ with $A=A(h,\delta,\nu)$ and $b=b(h,\delta,g)$. We use the discretization size $h=1/25$ (leading to $A\in\R^{N\times N}$ with $N=h^{-2}=625$) and fixed parameters $\nu=0.01,\,\delta=0.3$, but 25 different boundary conditions $g$. These boundary conditions set $g=0$ everywhere on $\partial\Omega$ except for a certain part of the right side of $\partial\Omega$; see~\cite[Example~2.2]{LieStr05} for details. The essential point is that we have only one matrix $A$, but 25 different right-hand sides $b^{(1)},\dots,b^{(25)}$. The matrix $A$ is highly nonnormal (here $\kappa(X)\approx 2.5\times 10^{17}$), and we again apply our implementation of MGS-GMRES with $x_0=0$.

\begin{figure}
\begin{center}
\includegraphics[width=0.475\textwidth]{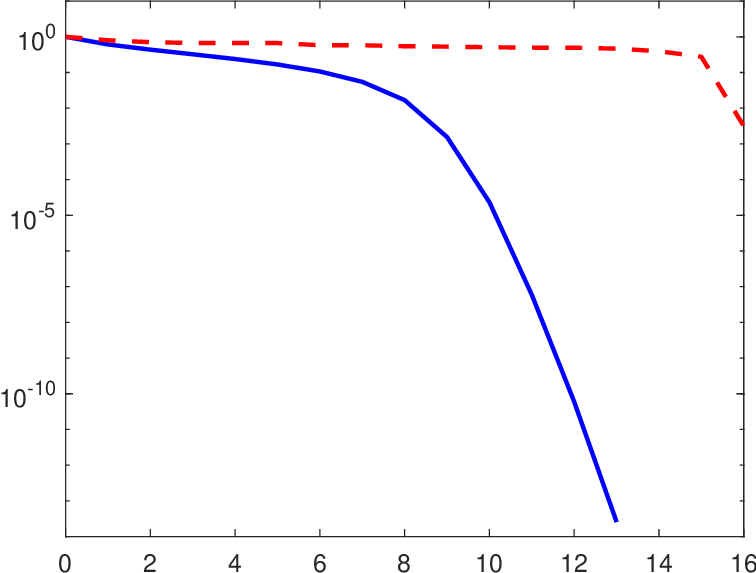}\hfill
\includegraphics[width=0.475\textwidth]{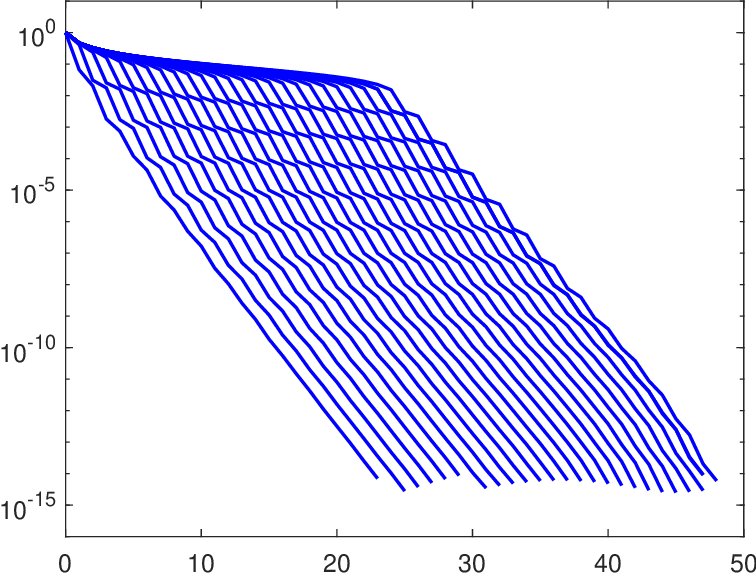}
\caption{Relative residual norms of GMRES for linear algebraic systems with the flipped Frank matrix (left) and two different right-hand sides, and a matrix from the SUPG discretization convection-diffusion model problem and 25 different right-hand sides (right).}\label{fig:frank_supg}
\end{center}
\end{figure}

\smallskip
\emph{Observations:} The relative GMRES residual norms for the two different matrices and the corresponding different right-hand sides are shown in Figure~\ref{fig:frank_supg}. For the flipped Frank matrix the solid blue curve corresponds to $b^{(1)}$, and dashed red curve to $b^{(2)}$. Apparently, GMRES converges much faster for $b^{(1)}$ than for $b^{(2)}$.

In the convection-diffusion problem we have for each $j=1,\dots,25$ a right-hand side $b^{(j)}$ for which GMRES has an initial phase of slow convergence (almost stagnation) for exactly $j-1$ steps. After the initial phase, the GMRES convergence speed is almost the same for all right-hand sides.

\smallskip
\emph{Explanation:}
The matrices used in this example are highly nonnormal, and the convergence bound \eqref{eqn:GMRESbd} is therefore of little use. Careful analyses of the relationship between $A$ and $b$ on a case-by-case basis are required to understand the GMRES convergence. For the SUPG discretized convection-diffusion problem such an analysis reveals that the length of the initial stagnation phase of GMRES for different boundary conditions depends on how many steps it takes to propagate the boundary information across the discretized domain by repeated multiplication with the matrix~$A$; see~\cite{LieStr05}. That analysis does not explain, however, how the convergence after the initial phase of stagnation depends on the parameters of the problem.

\subsection{GMRES convergence for nonzero initial approximate solutions}\label{sec:GMRESnonzero}
\phantom{}
\smallskip
\begin{tcolorbox}
\emph{Main point: Using a nonzero $x_0$ and measuring GMRES convergence with the relative residual norms can lead to an illusion of fast convergence. The relative residual norms can decrease rapidly due to dominant information created in $r_0$ by matrix-vector multiplication. Such information may not be present in the right-hand side, therefore the error can at the same time remain large, which can also lead to a premature stopping of the computation.}
\end{tcolorbox}

\smallskip
\emph{Setup:} We consider a linear algebraic system $Ax=b$, where  $A\in\R^{240\times 240}$ is the matrix {\tt steam1} from Matrix Market and $b=[1,\dots,1]^T/\sqrt{N}$. The matrix $A$ is nonsymmetric with $\kappa(A)\approx 2.8\times 10^7$. We compute $x=A^{-1}b$ using MATLAB's backslash operator. We apply our implementation of MGS-GMRES to $Ax=b$ with $x_0=0$, and with a normalized random $x_0$ generated using {\tt randn} in MATLAB. 

\begin{figure}
\begin{center}
\includegraphics[width=0.475\textwidth]{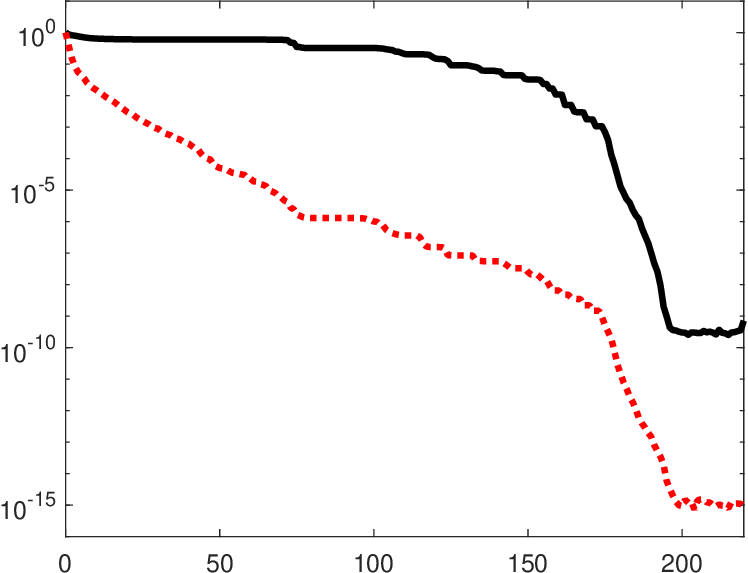}\hfill
\includegraphics[width=0.475\textwidth]{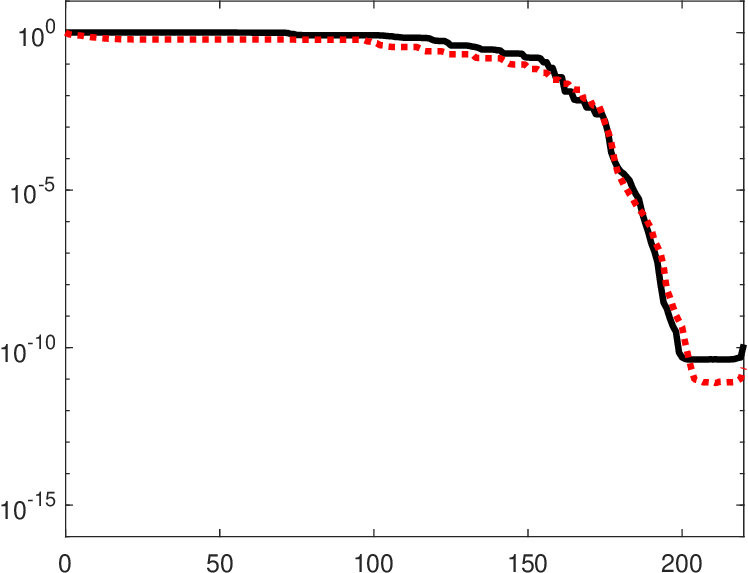}
\caption{Relative residual norms (left) and relative error norms (right) of GMRES applied to a linear algebraic system with the matrix {\tt steam1}, starting with $x_0=0$ (solid black) and a normalized random $x_0$ (dotted red).}\label{fig:GMRESnonzeroX0}
\end{center}
\end{figure}

\smallskip
\emph{Observations:} The relative residual norms $\|r_k\|_2/\|r_0\|_2$ and the relative error norms $\|x-x_k\|_2/\|x-x_0\|_2$ of GMRES are shown in the left and right plots of Figure~\ref{fig:GMRESnonzeroX0}, respectively. The solid black curves correspond to $x_0=0$, and the dashed red curves to the normalized random $x_0$.

The left plot shows that relative residual norms for the normalized random~$x_0$ decrease quickly, and reach approximately the machine precision level after about 200 iterations. For $x_0=0$ the relative residual norms decrease much slower, and eventually reach a level of approximately $10^{-10}$. On the other hand, the right plot shows that the relative error norms for both initial vectors are almost the same throughout the iteration, and approximately the same level of accuracy is reached in the end.

\smallskip
\emph{Explanation:} When the matrix~$A$ is ill conditioned, a nonzero $x_0$ may lead to $\|r_0\|_2=\|b-Ax_0\|_2\gg\|b\|_2$. (In our example we have $\|r_0\|_2\approx 4.5\times 10^6$ for the normalized random $x_0$.) The vector $r_0$ then contains an \emph{artificially created bias} towards the dominant information in the matrix~$A$ (such as large eigenvalues or singular values) that may not be related to the solution $x=A^{-1}b$. Elimination of this bias by GMRES can lead to a fast reduction of the residual norms particularly in the initial phase of the iteration. This creates an illusion of fast convergence (measured by the residual norms), and can lead to a premature stopping although $x_k$ is actually not close to $x$. Similar examples with different matrices can be seen in~\cite[Figure~5.13]{LieStr13} and~\cite{PaiStr02}. 

In order to avoid an illusion of fast convergence and a premature stopping, while still using a nonzero $x_0$, one can use the rescaling $\zeta_{\min} x_0$, where $\zeta_{\min}=(b^TAx_0)/\|Ax_0\|_2^2$ solves the approximation problem $\min_{\zeta}\|b-\zeta Ax_0\|_2$; see~\cite[p.~318]{LieStr13}. In our example $|\zeta_{\min}|\approx 2.5\times 10^{-8}$ holds for the normalized random $x_0$, and the residual and error norm curves of GMRES started with $\zeta_{\min} x_0$ are indistinguishable from the solid black curves in Figure~\ref{fig:GMRESnonzeroX0}.

\section{Concluding Remarks}\label{sec:Conclusions}

As an algorithmic idea, Krylov subspace methods are in their seventies, and their mathematical roots are closely related to much older objects like continued fractions, moments, and quadrature. It might seem that these methods are fully understood, but many questions about their behavior still remain open.
Focusing on the CG and the GMRES methods, we have presented an approach towards understanding Krylov subspace methods through a sequence of easily reproducible examples. 
In some cases this has led to clarification of common misunderstandings about the methods or to the formulation of open problems.

As argued in this paper, we believe that nonlinearity is the main mathematical asset as well as the beauty of Krylov subspace methods. The core misunderstanding in this context is the confusion of the nonlinear and data-adaptive behavior of the methods with linear convergence bounds. In general, the linearization of nonlinear phenomena is a highly useful technique that is applied throughout mathematics. However, it is important to realize that linearization can capture the behavior of strongly nonlinear phenomena only locally. Asymptotic convergence factors are based on linearization at infinity. When solving finite linear algebraic systems and methods with a mathematical finite termination property, this should be used with caution. Formally we do not have any limit at infinity, and practically we are interested in a few initial iterations (in comparison to the large number of unknowns). When solving linear operator equations in infinite dimensional Hilbert spaces, the convergence rate can become faster than any linear rate as the iterations proceed.

This paper illustrates that the behavior of Krylov subspace methods for some model problems can not be extrapolated to their behavior for practical problems. The condition number bounds, or arguments using clustered eigenvalues, can be useful under specific assumptions and with accepting the associated limitations. For a nice example, which was not discussed in this paper, we refer to the very insightful discussion of spectral equivalence of operators and operator preconditioning by Faber, Manteuffel, and Parter~\cite{FabManPar90}. The goal in this area is to construct preconditioners that guarantee condition number bounds independent of the discretization (and the PDE problem) parameters. This is certainly desirable and useful, but it might not be the final step in guaranteeing fast convergence of Krylov subspace methods. As they point out, 

\medskip
\begin{quote}
``For a fixed [discretization parameter], using a preconditioning strategy based upon an equivalent operator may not be superior to classical methods [ ... ]. Equivalence alone is not sufficient for a good preconditioning strategy. One must also choose an equivalent operator for which the bound [guaranteeing fast convergence] is small. The above observations indicate that a more precise measure of the `closeness' of two operators is required to evaluate preconditioning strategies.''
\end{quote}
\medskip

\noindent
Collaboration of researchers working on analysis of Krylov subspace methods, researchers in other fields (such as numerical PDEs), and practitioners who use the methods and are aware of the wider context, can lead to new ideas and paths of research. The approach started in the remarkable paper by Nielsen, Tveito and Hackbusch~\cite{NieTveHac09}, and continued later in \cite{GerMarNieStr19} and other works quoted in Section~\ref{sec:fpcg} can serve as an example.

The class of Krylov subspace methods is nowadays frequently seen as a computational toolbox whose mathematical investigation is more or less finished. The remaining open problems, if known to a broader audience at all, are being accepted as difficult and not worth investigating. We argue that Krylov subspace methods should still be seen as mathematical objects worth studying. Any progress in their understanding, even of their mathematical fundamentals, will bring us a step further in exploiting their full nonlinear computational potential.

\section{Appendix}

\subsection{Modern relevance of early works about Krylov subspace methods}\label{sec:intro1}

This paper is not focused on the history Krylov subspace methods. For that purpose we refer to the remarkable commented collection of works between 1948 and 1976 in the survey by Golub and O'Leary~\cite{GolOle89}, and, e.g., to the historical notes in~\cite{LieStr13}. We believe, however, that a thorough knowledge of the early works on Krylov subspace methods is very useful for understanding the strengths and weaknesses in the present use of these methods, and for obtaining perspectives for possible future developments. 

The papers of Lanczos, Hestenes and Stiefel published in 1952--53~\cite{HesSti52,Lan50,Lan52,Lan53} went far beyond the construction of a numerical algorithm for solving linear algebraic systems. They also considered the approximation of eigenvalues, and covered related mathematical background including orthogonal polynomials, continued fraction expansions, and Gauss quadrature. The paper~\cite{Lan52} considers ``purification'' of the starting vector of components along the large eigenvalues of the matrix $A$ using Chebyshev polynomials, followed by application of a version of CG. It explains in detail the sharp distinction between the action of the Chebyshev polynomials and CG polynomials, where the former addresses the continuum of an interval while the latter ``attenuate due to the nearness of their zeros'' to the eigenvalues of $A$; cf.~\cite[p.~46 and Fig.~2]{Lan52}. The paper~\cite{Lan53} suggests in this context polynomial preconditioning much earlier than other publications, but this work has been largely overlooked; see also the commentaries~\cite{Saa98,Ste91}. The papers~\cite{HesSti52,Lan50,Lan52,Lan53} also address quality measures of the computed iterative approximations. Remarkably, the early researchers also analyzed computer implementations of the methods, including the effects of rounding errors.

The papers of Karush~\cite{Kar52} and Hayes~\cite{Hay54} from 1952 and 1954, respectively, prove superlinear convergence of CG for $\mathcal{A}=\mathrm{Id}+\mathcal{C}$, where $\mathrm{Id}$ is the identity and $\mathcal{C}$ is a compact operator, concluding that the CG rate of convergence for such operators exceeds any given linear rate as the iteration proceeds. This result was, decades later, repeatedly rediscovered.  A bit later (with the original print in Russian from 1958), Vorobyev put this new development into the context of a computationally feasible formulation of the method of moments~\cite{VorBook65}; see also the insightful paper by Brezinski~\cite{Bre96} on its generalizations and the relationship with other related concepts, including the approach by Lanczos.

Unfortunately, many results of these early works remain almost unnoticed, and the common state-of-the-art literature even contains views that contradict some of them. On the other hand, a lack of theoretical understanding does not prohibit the practical application of the methods. Thus, the following quote from Forsythe's paper from 1953~\cite{For53} is to some extend still valid for the class of Krylov subspace methods today: 

\medskip
\begin{quote}
``It is remarkable how little is really understood about most of the methods for solving $Ax=b$. They are being used nevertheless, and yield
answers. This disparity between theory and practice appears to be typical of the gulf between the {\em science} and the {\em art} of numerical analysis.''
\end{quote}
\medskip

\noindent
These words, written in 1952 and published in 1953, were prophetic for the situation in the area of Krylov subspace methods for decades to come. 

Two common beliefs frequently repeated in the numerical analysis literature are that CG was originally considered as a direct method, and that it was, after unsatisfactory computational experiences, abandoned for about two decades as a numerical approach for linear systems. As mentioned above, CG was in the original papers considered as an iterative method. The story of considering it as direct method started later.  Moreover, while the interest in CG  within the numerical analysis community temporarily declined in the 1950s and 1960s (and was only revived later by Reid, Axelsson, Evans, Concus, Golub, O'Leary, Paige, Saunders, Meijerink, van der Vorst, and others in the 1970s), it actually was used by practitioners in the 1950s and 1960s; see, e.g., the collection of works given in~\cite{GolOle89}. Many important theoretical ideas appeared independently of numerical linear algebra in different fields; see, e.g., the developments in computational physics and chemistry~\cite{Gor68} (cf. also the independent developments in~\cite{GolWel69}), \cite{Rei79,SchSch66}, as well as in the works on matching moments leading to model reduction of large dynamical systems in~\cite{Gra74}.

The history of Krylov subspace methods illustrates discontinuities in the development of science and even loss of knowledge, as well as an unhealthy lack of communication between different fields.

\subsection{Preconditioning and analysis of Krylov subspace methods}\label{sec:intro2}

Practical computations with Krylov subspace methods require preconditioning. From a purely algebraical point of view, this is typically interpreted as the transformation of the original linear algebraic system into an equivalent system that is, for the given iterative method, more easily solvable. Another view is the acceleration of convergence of the computed approximations through transformation of the iterative process. For symmetric positive definite matrices and CG, the goal of preconditioning is typically stated as (see, e.g.,~\cite{Ben02}): ``Hopefully, the transformed (preconditioned) matrix will have a smaller condition number, and/or eigenvalues clustered around one.'' When this can be achieved, and a small condition number guarantees an acceptable approximate solution within a few iterations (which can happen, e.g., for preconditioners constructed via domain decomposition methods incorporating coarse components) the matter seems to be resolved. However, if in such cases tight bounds are available for the extreme eigenvalues of the preconditioned matrix, then it can be computationally competitive or even more efficient to apply the Chebyshev semiiterative method instead of CG. This alternative approach is easily parallelizable and numerically stable. Similar to the quote above, the following insight of Forsythe~\cite[p. 318]{For53} is still valid and not fully appreciated today:

\medskip
\begin{quote}
``The belief is widespread that the condition of a system $Ax = b$ has a decisive influence on the convergence of an iterative solution [ ... ]; this can not always be true. Even when it is true for an iterative process, it may be possible actually to take advantage of the poor condition of $Ax = b$ in converting the slow process into an accelerated method which converges rapidly. There is a great need for clarification of the group of ideas associated with `condition'.''
\end{quote}
\medskip

\noindent
As illustrated in Section~\ref{sec:pcg} in this paper, the nonlinear adaptivity of CG to the location of the individual eigenvalues indicates that a smaller condition number does not necessarily lead to faster convergence. One of the examples in Section~\ref{sec:pcg} is taken from the paper~\cite{GerMarNieStr19}, which contains a substantial theoretical result on approximating {\em all eigenvalues of the preconditioned matrix}, and hence shows that information about all eigenvalues for important PDE problems can be available a priori at a negligible computational cost. In addition, recent results of Colbrook, Horning, and Townsend~\cite{ColHor22,ColHorTow21} show how to compute smoothed approximations of spectral measures for infinite dimensional self-adjoint operators. 

A large family of preconditioners in numerical PDEs is motivated by the spectral or norm equivalence of operators in infinite dimensional Hilbert spaces (hence the name {\em operator preconditioning}); see, e.g., the beautiful and rarely quoted early papers by Concus and Golub published in 1973~\cite{ConGol73}, and by Faber, Manteuffel, and Parter published in 1990~\cite{FabManPar90}. Infinite dimensional operator preconditioning and algebraic preconditioning in the context of using CG are linked together in~\cite{MSB15}. In particular, it is shown that any algebraic preconditioning can be put into the operator preconditioning framework. This reveals the close mathematical connection between the choice of the inner product in the infinite dimensional Hilbert spaces in operator preconditioning and the choice of the discretization basis determining the associated linear algebraic problem. A survey of the abstract framework of operator preconditioning based on decomposition of infinite dimensional Hilbert subspaces is given in~\cite{HPS20}.

Many results in the literature are devoted to problem-specific preconditioning. The overall state-of-the-art in the area of preconditioning still remains well-characterized by the following quote from the book of Saad published in 2003~\cite[p. 297]{Saa03}:

\medskip
\begin{quote}
``Finding a good preconditioner to solve a given sparse linear system is often viewed as a combination of art and science. Theoretical results are rare and some methods work surprisingly well, often despite expectations.''
\end{quote}
\medskip

\noindent
This underlines the need for further research on preconditioning, which can not be separated from further analysis of Krylov subspace methods.

\section*{Acknowledgements} We thank Tyler Chen, Justus Ramme, Daniel Szyld, and Jan Zur for helpful suggestions, and Jan Pape\v{z} for help with Figure \ref{fig:pcgvscg}.

\bibliography{CarLieStr24.bbl}

\begin{thebibliography}{100}

\bibitem{AriPtaStr98}
{\sc M.~Arioli, V.~Pt\'{a}k, and Z.~Strako\v{s}}, {\em Krylov sequences of
  maximal length and convergence of {GMRES}}, BIT, 38 (1998), pp.~636--643.

\bibitem{Arn51}
{\sc W.~E. Arnoldi}, {\em The principle of minimized iteration in the solution
  of the matrix eigenvalue problem}, Quart. Appl. Math., 9 (1951), pp.~17--29.

\bibitem{AshManSay90}
{\sc S.~F. Ashby, T.~A. Manteuffel, and P.~E. Saylor}, {\em A taxonomy for
  conjugate gradient methods}, SIAM J. Numer. Anal., 27 (1990), pp.~1542--1568.

\bibitem{Ben02}
{\sc M.~Benzi}, {\em Preconditioning techniques for large linear systems: a
  survey}, J. Comput. Phys., 182 (2002), pp.~418--477.

\bibitem{Bre96}
{\sc C.~Brezinski}, {\em The methods of {V}orobyev and {L}anczos}, Linear
  Algebra Appl., 234 (1996), pp.~21--41.

\bibitem{CarStr20}
{\sc E.~Carson and Z.~Strako\v{s}}, {\em On the cost of iterative
  computations}, Philos. Trans. Roy. Soc. A, 378 (2020), pp.~20190050, 22.

\bibitem{Car15}
{\sc E.~C. Carson}, {\em Communication-avoiding {K}rylov subspace methods in
  theory and practice}, PhD thesis, University of California, Berkeley, 2015.

\bibitem{CarRozStrTicTum18}
{\sc E.~C. Carson, M.~Rozlo\v{z}n\'{\i}k, Z.~Strako\v{s}, P.~Tich\'{y}, and
  M.~T\r{u}ma}, {\em The numerical stability analysis of pipelined conjugate
  gradient methods: historical context and methodology}, SIAM J. Sci. Comput.,
  40 (2018), pp.~A3549--A3580.

\bibitem{ColHor22}
{\sc M.~Colbrook and A.~Horning}, {\em Specsolve: Spectral methods for spectral
  measures}, arXiv preprint 2201.01314,  (2022).

\bibitem{ColHorTow21}
{\sc M.~Colbrook, A.~Horning, and A.~Townsend}, {\em Computing spectral
  measures of self-adjoint operators}, SIAM Rev., 63 (2021), pp.~489--524.

\bibitem{ConGol73}
{\sc P.~Concus and G.~H. Golub}, {\em Use of fast direct methods for the
  efficient numerical solution of nonseparable elliptic equations}, SIAM J.
  Numer. Anal., 10 (1973), pp.~1103--1120.

\bibitem{CooYetAguGirVan18}
{\sc S.~Cools, E.~F. Yetkin, E.~Agullo, L.~Giraud, and W.~Vanroose}, {\em
  Analyzing the effect of local rounding error propagation on the maximal
  attainable accuracy of the pipelined conjugate gradient method}, SIAM J.
  Matrix Anal. Appl., 39 (2018), pp.~426--450.

\bibitem{DrkGreRozStr95}
{\sc J.~Drko\v{s}ov\'{a}, A.~Greenbaum, M.~Rozlo\v{z}n\'{\i}k, and
  Z.~Strako\v{s}}, {\em Numerical stability of {GMRES}}, BIT, 35 (1995),
  pp.~309--330.

\bibitem{DruKni95}
{\sc V.~Druskin and L.~Knizhnerman}, {\em Krylov subspace approximation of
  eigenpairs and matrix functions in exact and computer arithmetic}, Numer.
  Linear Algebra Appl., 2 (1995), pp.~205--217.

\bibitem{DuiMeuSadStr14}
{\sc J.~Duintjer~Tebbens, G.~Meurant, H.~Sadok, and Z.~Strako\v{s}}, {\em On
  investigating {GMRES} convergence using unitary matrices}, Linear Algebra
  Appl., 450 (2014), pp.~83--107.

\bibitem{Ede97}
{\sc A.~Edelman}, {\em The probability that a random real {G}aussian matrix has
  {$k$} real eigenvalues, related distributions, and the circular law}, J.
  Multivariate Anal., 60 (1997), pp.~203--232.

\bibitem{ElmRamSil14}
{\sc H.~C. Elman, A.~Ramage, and D.~J. Silvester}, {\em I{FISS}: a
  computational laboratory for investigating incompressible flow problems},
  SIAM Rev., 56 (2014), pp.~261--273.

\bibitem{ElmSilWat14}
{\sc H.~C. Elman, D.~J. Silvester, and A.~J. Wathen}, {\em Finite elements and
  fast iterative solvers: with applications in incompressible fluid dynamics},
  Numerical Mathematics and Scientific Computation, Oxford University Press,
  Oxford, second~ed., 2014.

\bibitem{Emb03}
{\sc M.~Embree}, {\em The tortoise and the hare restart {GMRES}}, SIAM Rev., 45
  (2003), pp.~259--266.

\bibitem{FabJouKniMan96}
{\sc V.~Faber, W.~Joubert, E.~Knill, and T.~Manteuffel}, {\em Minimal residual
  method stronger than polynomial preconditioning}, SIAM J. Matrix Anal. Appl.,
  17 (1996), pp.~707--729.

\bibitem{FabLieTic13}
{\sc V.~Faber, J.~Liesen, and P.~Tich\'{y}}, {\em Properties of worst-case
  {GMRES}}, SIAM J. Matrix Anal. Appl., 34 (2013), pp.~1500--1519.

\bibitem{FabManPar90}
{\sc V.~Faber, T.~A. Manteuffel, and S.~V. Parter}, {\em On the theory of
  equivalent operators and application to the numerical solution of uniformly
  elliptic partial differential equations}, Adv. in Appl. Math., 11 (1990),
  pp.~109--163.

\bibitem{Fle76}
{\sc R.~Fletcher}, {\em Conjugate gradient methods for indefinite systems}, in
  Numerical analysis (Proc 6th Biennial Dundee Conf., Univ. Dundee, Dundee,
  1975), Springer, Berlin, 1976, pp.~73--89. Lecture Notes in Math., Vol. 506.

\bibitem{For53}
{\sc G.~E. Forsythe}, {\em Solving linear algebraic equations can be
  interesting}, Bull. Amer. Math. Soc., 59 (1953), pp.~299--329.

\bibitem{FreGolNac92}
{\sc R.~W. Freund, G.~H. Golub, and N.~M. Nachtigal}, {\em Iterative solution
  of linear systems}, Acta Numer., 1 (1992), pp.~57--100.

\bibitem{FreNac91}
{\sc R.~W. Freund and N.~M. Nachtigal}, {\em Q{MR}: a quasi-minimal residual
  method for non-{H}ermitian linear systems}, Numer. Math., 60 (1991),
  pp.~315--339.

\bibitem{Gan34}
{\sc F.~R. Gantmacher}, {\em On the algebraic analysis of {K}rylov's method of
  transforming the secular equation}, Trans. Second Math. Congress, II (1934),
  pp.~45--48.
\newblock In Russian. Title translation as in~\cite{GanBook59}.

\bibitem{GanBook59}
\leavevmode\vrule height 2pt depth -1.6pt width 23pt, {\em The Theory of
  Matrices. {V}ols. 1, 2}, Chelsea Publishing Co., New York, 1959.

\bibitem{Ger13}
{\sc T.~Gergelits}, {\em Analysis of {K}rylov subspace methods}, 2013.
\newblock Master's Thesis, Charles University, Faculty of Mathematics and
  Physics.

\bibitem{GerHneKub17}
{\sc T.~Gergelits, I.~Hn{\v{e}}tynkov{\'a}, and M.~Kub{\'\i}nov{\'a}}, {\em
  Relating computed and exact entities in methods based on {Lanczos}
  tridiagonalization}, in International Conference on High Performance
  Computing in Science and Engineering, Springer, 2017, pp.~73--87.

\bibitem{GerMarNieStr19}
{\sc T.~Gergelits, K.-A. Mardal, B.~F. Nielsen, and Z.~Strako\v{s}}, {\em
  Laplacian preconditioning of elliptic {PDE}s: localization of the eigenvalues
  of the discretized operator}, SIAM J. Numer. Anal., 57 (2019),
  pp.~1369--1394.

\bibitem{GerNieStr22}
{\sc T.~Gergelits, B.~r.~F. Nielsen, and Z.~Strako\v{s}}, {\em Numerical
  approximation of the spectrum of self-adjoint operators in operator
  preconditioning}, Numer. Algorithms, 91 (2022), pp.~301--325.

\bibitem{GirLanRozEsh05}
{\sc L.~Giraud, J.~Langou, M.~Rozlo\v{z}n\'{\i}k, and J.~van~den Eshof}, {\em
  Rounding error analysis of the classical {G}ram-{S}chmidt orthogonalization
  process}, Numer. Math., 101 (2005), pp.~87--100.

\bibitem{GolMeu97}
{\sc G.~H. Golub and G.~Meurant}, {\em Matrices, moments and quadrature. {II}.
  {H}ow to compute the norm of the error in iterative methods}, BIT, 37 (1997),
  pp.~687--705.
\newblock Direct methods, linear algebra in optimization, iterative methods
  (Toulouse, 1995/1996).

\bibitem{GolOle89}
{\sc G.~H. Golub and D.~P. O'Leary}, {\em Some history of the conjugate
  gradient and {L}anczos algorithms: 1948--1976}, SIAM Rev., 31 (1989),
  pp.~50--102.

\bibitem{GolStr94}
{\sc G.~H. Golub and Z.~Strako\v{s}}, {\em Estimates in quadratic formulas},
  Numer. Algorithms, 8 (1994), pp.~241--268.

\bibitem{GolVdV97}
{\sc G.~H. Golub and H.~A. van~der Vorst}, {\em Closer to the solution:
  iterative linear solvers}, in The state of the art in numerical analysis
  (York, 1996), vol.~63 of Inst. Math. Appl. Conf. Ser. New Ser., Oxford Univ.
  Press, New York, 1997, pp.~63--92.

\bibitem{GolVar61a}
{\sc G.~H. Golub and R.~S. Varga}, {\em Chebyshev semi-iterative methods,
  successive over-relaxation iterative methods, and second order {R}ichardson
  iterative methods. {I}}, Numer. Math., 3 (1961), pp.~147--156.

\bibitem{GolVar61b}
\leavevmode\vrule height 2pt depth -1.6pt width 23pt, {\em Chebyshev
  semi-iterative methods, successive over-relaxation iterative methods, and
  second order {R}ichardson iterative methods. {II}}, Numer. Math., 3 (1961),
  pp.~157--168.

\bibitem{GolWel69}
{\sc G.~H. Golub and J.~H. Welsch}, {\em Calculation of {G}auss quadrature
  rules}, Math. Comp. 23 (1969), 221-230; addendum, ibid., 23 (1969),
  pp.~A1--A10.

\bibitem{Gor68}
{\sc R.~G. Gordon}, {\em Error bounds in equilibrium statistical mechanics}, J.
  Math. Phys., 9 (1968), pp.~655--663.

\bibitem{Gra74}
{\sc W.~B. Gragg}, {\em Matrix interpretations and applications of the
  continued fraction algorithm}, Rocky Mountain J. Math., 4 (1974),
  pp.~213--225.

\bibitem{Gre79}
{\sc A.~Greenbaum}, {\em Comparison of splittings used with the conjugate
  gradient algorithm}, Numer. Math., 33 (1979), pp.~181--193.

\bibitem{Gre89}
\leavevmode\vrule height 2pt depth -1.6pt width 23pt, {\em Behavior of slightly
  perturbed {L}anczos and conjugate-gradient recurrences}, Linear Algebra
  Appl., 113 (1989), pp.~7--63.

\bibitem{GreBook97}
\leavevmode\vrule height 2pt depth -1.6pt width 23pt, {\em Iterative Methods
  for Solving Linear Systems}, vol.~17 of Frontiers in Applied Mathematics,
  SIAM, Philadelphia, PA, 1997.

\bibitem{GreGur94}
{\sc A.~Greenbaum and L.~Gurvits}, {\em Max-min properties of matrix factor
  norms}, SIAM J. Sci. Comput., 15 (1994), pp.~348--358.

\bibitem{GrePtaStr96}
{\sc A.~Greenbaum, V.~Pt\'{a}k, and Z.~Strako\v{s}}, {\em Any nonincreasing
  convergence curve is possible for {GMRES}}, SIAM J. Matrix Anal. Appl., 17
  (1996), pp.~465--469.

\bibitem{GreStr92}
{\sc A.~Greenbaum and Z.~Strako\v{s}}, {\em Predicting the behavior of finite
  precision {L}anczos and conjugate gradient computations}, SIAM J. Matrix
  Anal. Appl., 13 (1992), pp.~121--137.

\bibitem{GreStr94}
\leavevmode\vrule height 2pt depth -1.6pt width 23pt, {\em Matrices that
  generate the same {K}rylov residual spaces}, in Recent advances in iterative
  methods, vol.~60 of IMA Vol. Math. Appl., Springer, New York, 1994,
  pp.~95--118.

\bibitem{GreTre94}
{\sc A.~Greenbaum and L.~N. Trefethen}, {\em G{MRES}/{CR} and
  {A}rnoldi/{L}anczos as matrix approximation problems}, SIAM J. Sci. Comput.,
  15 (1994), pp.~359--368.

\bibitem{GutStr00}
{\sc M.~H. Gutknecht and Z.~Strako\v{s}}, {\em Accuracy of two three-term and
  three two-term recurrences for {K}rylov space solvers}, SIAM J. Matrix Anal.
  Appl., 22 (2000), pp.~213--229.

\bibitem{HacBook94}
{\sc W.~Hackbusch}, {\em Iterative Solution of Large Sparse Systems of
  Equations}, vol.~95 of Applied Mathematical Sciences, Springer-Verlag, New
  York, 1994.
\newblock Translated and revised from the 1991 German original.

\bibitem{Hay54}
{\sc R.~M. Hayes}, {\em Iterative Methods for Solving Linear Problems in
  Hilbert Space}, PhD thesis, Univ. of California at Los Angeles, Los Angeles,
  CA, 1954.

\bibitem{HesSti52}
{\sc M.~R. Hestenes and E.~Stiefel}, {\em Methods of conjugate gradients for
  solving linear systems}, J. Research Nat. Bur. Standards, 49 (1952),
  pp.~409--436.

\bibitem{HPS20}
{\sc J.~Hrn\v{c}\'{\i}\v{r}, I.~Pultarov\'{a}, and Z.~Strako\v{s}}, {\em
  Decomposition into subspaces preconditioning: abstract framework}, Numer.
  Algorithms, 83 (2020), pp.~57--98.

\bibitem{Jen77}
{\sc A.~Jennings}, {\em Influence of the eigenvalue spectrum on the convergence
  rate of the conjugate gradient method}, J. Inst. Math. Appl., 20 (1977),
  pp.~61--72.

\bibitem{Jou94}
{\sc W.~Joubert}, {\em A robust {GMRES}-based adaptive polynomial
  preconditioning algorithm for nonsymmetric linear systems}, SIAM J. Sci.
  Comput., 15 (1994), pp.~427--439.

\bibitem{Kar52}
{\sc W.~Karush}, {\em Convergence of a method of solving linear problems},
  Proc. Amer. Math. Soc., 3 (1952), pp.~839--851.

\bibitem{Kry31}
{\sc A.~N. Krylov}, {\em On the numerical solution of the equation by which the
  frequency of small oscillations is determined in technical problems}, Izv.
  Akad. Nauk SSSR, Ser. Fiz.-Mat., 4 (1931), pp.~491--539.
\newblock In Russian. Title translation as in~\cite{GanBook59}.

\bibitem{Kui00}
{\sc A.~B.~J. Kuijlaars}, {\em Which eigenvalues are found by the lanczos
  method?}, SIAM J. Matrix Anal. Appl., 22 (2000), pp.~306--321.

\bibitem{Lad23}
{\sc M.~Ladeck\'{y}, R.~J. Leute, A.~Falsafi, I.~Pultarov\'{a}, L.~Pastewka,
  T.~Junge, and J.~Zeman}, {\em An optimal preconditioned {FFT}-accelerated
  finite element solver for homogenization}, Appl. Math. Comput., 446 (2023),
  pp.~Paper No. 127835, 19.

\bibitem{Lad21}
{\sc M.~Ladeck\'{y}, I.~Pultarov\'{a}, and J.~Zeman}, {\em Guaranteed two-sided
  bounds on all eigenvalues of preconditioned diffusion and elasticity problems
  solved by the finite element method}, Appl. Math., 66 (2021), pp.~21--42.

\bibitem{Lan50}
{\sc C.~Lanczos}, {\em An iteration method for the solution of the eigenvalue
  problem of linear differential and integral operators}, J. Research Nat. Bur.
  Standards, 45 (1950), pp.~255--282.

\bibitem{Lan52}
\leavevmode\vrule height 2pt depth -1.6pt width 23pt, {\em Solution of systems
  of linear equations by minimized iterations}, J. Research Nat. Bur.
  Standards, 49 (1952), pp.~33--53.

\bibitem{Lan53}
\leavevmode\vrule height 2pt depth -1.6pt width 23pt, {\em Chebyshev
  polynomials in the solution of large-scale linear systems}, in Proceedings of
  the {A}ssociation for {C}omputing {M}achinery, {T}oronto, 1952, Sauls
  Lithograph Co. (for the Association for Computing Machinery), Washington, D.
  C., 1953, pp.~124--133.

\bibitem{Leu22}
{\sc R.~J. Leute, M.~Ladeck\'{y}, A.~Falsafi, I.~J\"{o}dicke, I.~Pultarov\'{a},
  J.~Zeman, T.~Junge, and L.~Pastewka}, {\em Elimination of ringing artifacts
  by finite-element projection in {FFT}-based homogenization}, J. Comput.
  Phys., 453 (2022), pp.~Paper No. 110931, 20.

\bibitem{LieStr05}
{\sc J.~Liesen and Z.~Strako\v{s}}, {\em G{MRES} convergence analysis for a
  convection-diffusion model problem}, SIAM J. Sci. Comput., 26 (2005),
  pp.~1989--2009.

\bibitem{LieStr13}
\leavevmode\vrule height 2pt depth -1.6pt width 23pt, {\em Krylov subspace
  methods}, Numerical Mathematics and Scientific Computation, Oxford University
  Press, Oxford, 2013.
\newblock Principles and analysis.

\bibitem{LieTic04}
{\sc J.~Liesen and P.~Tich{\'y}}, {\em Convergence analysis of {K}rylov
  subspace methods}, GAMM Mitt. Ges. Angew. Math. Mech., 27 (2004),
  pp.~153--173.

\bibitem{LieTic14}
\leavevmode\vrule height 2pt depth -1.6pt width 23pt, {\em Max-min and min-max
  approximation problems for normal matrices revisited}, Electron. Trans.
  Numer. Anal., 41 (2014), pp.~159--166.

\bibitem{LinSaaYan16}
{\sc L.~Lin, Y.~Saad, and C.~Yang}, {\em Approximating spectral densities of
  large matrices}, SIAM Rev., 58 (2016), pp.~34--65.

\bibitem{Luz31}
{\sc N.~N. Luzin}, {\em On {K}rylov's method for transforming the secular
  equation}, Izv. Akad. Nauk SSSR, Ser. Fiz.-Mat., 7 (1931), pp.~903--958.
\newblock Title translation as in~\cite{GanBook59}.

\bibitem{MSB15}
{\sc J.~M\'{a}lek and Z.~Strako\v{s}}, {\em Preconditioning and the conjugate
  gradient method in the context of solving {PDE}s}, vol.~1 of SIAM Spotlights,
  Society for Industrial and Applied Mathematics (SIAM), Philadelphia, PA,
  2015.

\bibitem{MeuBook06}
{\sc G.~Meurant}, {\em The {L}anczos and conjugate gradient algorithms. {F}rom
  theory to finite precision computations}, vol.~19 of Software, Environments,
  and Tools, Society for Industrial and Applied Mathematics (SIAM),
  Philadelphia, PA, 2006.

\bibitem{Meu20}
\leavevmode\vrule height 2pt depth -1.6pt width 23pt, {\em On prescribing the
  convergence behavior of the conjugate gradient algorithm}, Numer. Algorithms,
  84 (2020), pp.~1353--1380.

\bibitem{MeuDui20}
{\sc G.~Meurant and J.~Duintjer~Tebbens}, {\em Krylov methods for nonsymmetric
  linear systems---from theory to computations}, vol.~57 of Springer Series in
  Computational Mathematics, Springer, Cham, [2020] \copyright 2020.

\bibitem{MeuPapTic21}
{\sc G.~Meurant, J.~Pape\v{z}, and P.~Tich\'{y}}, {\em Accurate error
  estimation in {CG}}, Numer. Algorithms, 88 (2021), pp.~1337--1359.

\bibitem{MeuStr06}
{\sc G.~Meurant and Z.~Strako{\v{s}}}, {\em The {L}anczos and conjugate
  gradient algorithms in finite precision arithmetic}, Acta Numer., 15 (2006),
  pp.~471--542.

\bibitem{MorNocSie02}
{\sc P.~Morin, R.~H. Nochetto, and K.~G. Siebert}, {\em Convergence of adaptive
  finite element methods}, SIAM Rev., 44 (2002), pp.~631--658.
\newblock Revised reprint of ``Data oscillation and convergence of adaptive
  FEM'' [SIAM J. Numer. Anal. {\bf 38} (2000), no. 2, 466--488 (electronic);
  MR1770058 (2001g:65157)].

\bibitem{MurGolWat00}
{\sc M.~F. Murphy, G.~H. Golub, and A.~J. Wathen}, {\em A note on
  preconditioning for indefinite linear systems}, SIAM J. Sci. Comput., 21
  (2000), pp.~1969--1972.

\bibitem{NieStr24}
{\sc B.~F. Nielsen and Z.~Strako\v{s}}, {\em A simple formula for the
  generalized spectrum of second order self-adjoint differential operators},
  SIAM Rev.,  (2024).

\bibitem{NieTveHac09}
{\sc B.~F. Nielsen, A.~Tveito, and W.~Hackbusch}, {\em Preconditioning by
  inverting the {L}aplacian: an analysis of the eigenvalues}, IMA J. Numer.
  Anal., 29 (2009), pp.~24--42.

\bibitem{Not93}
{\sc Y.~Notay}, {\em On the convergence rate of the conjugate gradients in
  presence of rounding errors}, Numer. Math., 65 (1993), pp.~301--317.

\bibitem{Pai80}
{\sc C.~C. Paige}, {\em Accuracy and effectiveness of the {L}anczos algorithm
  for the symmetric eigenproblem}, Linear Algebra Appl., 34 (1980),
  pp.~235--258.

\bibitem{PaiRozStr06}
{\sc C.~C. Paige, M.~Rozlo\v{z}n\'{\i}k, and Z.~Strako\v{s}}, {\em Modified
  {G}ram-{S}chmidt ({MGS}), least squares, and backward stability of
  {MGS}-{GMRES}}, SIAM J. Matrix Anal. Appl., 28 (2006), pp.~264--284.

\bibitem{PaiSau75}
{\sc C.~C. Paige and M.~A. Saunders}, {\em Solution of sparse indefinite
  systems of linear equations}, SIAM J. Numer. Anal., 12 (1975), pp.~617--629.

\bibitem{PaiStr02a}
{\sc C.~C. Paige and Z.~Strako{\v{s}}}, {\em Bounds for the least squares
  distance using scaled total least squares}, Numer. Math., 91 (2002),
  pp.~93--115.

\bibitem{PaiStr02}
\leavevmode\vrule height 2pt depth -1.6pt width 23pt, {\em Residual and
  backward error bounds in minimum residual {K}rylov subspace methods}, SIAM J.
  Sci. Comput., 23 (2002), pp.~1898--1923.

\bibitem{PozStr19}
{\sc J.~W. Pearson and J.~Pestana}, {\em Preconditioners for {K}rylov subspace
  methods: an overview}, GAMM-Mitt., 43 (2020), pp.~e202000015, 35.

\bibitem{Pul21}
{\sc I.~Pultarov\'{a} and M.~Ladeck\'{y}}, {\em Two-sided guaranteed bounds to
  individual eigenvalues of preconditioned finite element and finite difference
  problems}, Numer. Linear Algebra Appl., 28 (2021), pp.~Paper No. e2382, 15.

\bibitem{ReiTre92}
{\sc L.~Reichel and L.~N. Trefethen}, {\em Eigenvalues and pseudo-eigenvalues
  of {T}oeplitz matrices}, Linear Algebra Appl., 162/164 (1992), pp.~153--185.

\bibitem{Rei71}
{\sc J.~K. Reid}, {\em On the method of conjugate gradients for the solution of
  large sparse systems of linear equations}, in Large sparse sets of linear
  equations ({P}roc. {C}onf., {S}t. {C}atherine's {C}oll., {O}xford, 1970),
  Academic Press, London-New York, 1971, pp.~231--254.

\bibitem{Rei79}
{\sc W.~P. Reinhardt}, {\em $l^2$ discretization of atomic and molecular
  electronic continua: {M}oment, quadrature and $j$-matrix techniques}, Comp.
  Phys. Comm., 17 (1979), pp.~1--21.

\bibitem{RigGac67}
{\sc J.-L. Rigal and J.~Gaches}, {\em On the compatibility of a given solution
  with the data of a linear system}, J. Assoc. Comput. Mach., 14 (1967),
  pp.~543--548.

\bibitem{Saa81}
{\sc Y.~Saad}, {\em Krylov subspace methods for solving large unsymmetric
  linear systems}, Math. Comp., 37 (1981), pp.~105--126.

\bibitem{Saa98}
\leavevmode\vrule height 2pt depth -1.6pt width 23pt, {\em Commentary on
  {L}anczos's ``{C}hebyshev polynomials in the solution of large-scale linear
  systems''}, in {C}ornelius {L}anczos {C}ollected {P}ublished {P}apers with
  {C}ommentaries, {V}ol. VI, North Carolina State University, Raleigh, North
  Carolina, 1998, pp.~3--527--3--528.

\bibitem{Saa03}
\leavevmode\vrule height 2pt depth -1.6pt width 23pt, {\em Iterative methods
  for sparse linear systems}, Society for Industrial and Applied Mathematics,
  Philadelphia, PA, second~ed., 2003.

\bibitem{SaaSch86}
{\sc Y.~Saad and M.~H. Schultz}, {\em G{MRES}: a generalized minimal residual
  algorithm for solving nonsymmetric linear systems}, SIAM J. Sci. Statist.
  Comput., 7 (1986), pp.~856--869.

\bibitem{SchSch66}
{\sc L.~Schlessinger and C.~Schwartz}, {\em Analyticity as a useful
  computational tool}, Phys. Rev. Lett., 16 (1966), pp.~1173--1174.

\bibitem{Ste91}
{\sc G.~W. Stewart}, {\em Lanczos and linear systems}, Preprint UMIACS-TR-91,
  University of Maryland, 1991.

\bibitem{Str91}
{\sc Z.~Strako\v{s}}, {\em On the real convergence rate of the conjugate
  gradient method}, Linear Algebra Appl., 154/156 (1991), pp.~535--549.

\bibitem{Str09}
{\sc Z.~Strako\v{s}}, {\em Model reduction using the {V}orobyev moment
  problem}, Numer. Algorithms, 51 (2009), pp.~363--379.

\bibitem{StrTic02}
{\sc Z.~Strako\v{s} and P.~Tich\'{y}}, {\em On error estimation in the
  conjugate gradient method and why it works in finite precision computations},
  Electron. Trans. Numer. Anal., 13 (2002), pp.~56--80.

\bibitem{StrTic05}
\leavevmode\vrule height 2pt depth -1.6pt width 23pt, {\em Error estimation in
  preconditioned conjugate gradients}, BIT, 45 (2005), pp.~789--817.

\bibitem{Toh97}
{\sc K.-C. Toh}, {\em G{MRES} vs.\ ideal {GMRES}}, SIAM J. Matrix Anal. Appl.,
  18 (1997), pp.~30--36.

\bibitem{Tre92}
{\sc L.~N. Trefethen}, {\em Pseudospectra of matrices}, in Numerical analysis
  1991 ({D}undee, 1991), vol.~260 of Pitman Res. Notes Math. Ser., Longman Sci.
  Tech., Harlow, 1992, pp.~234--266.

\bibitem{TreBau97}
{\sc L.~N. Trefethen and D.~Bau, III}, {\em Numerical linear algebra}, Society
  for Industrial and Applied Mathematics (SIAM), Philadelphia, PA, 1997.

\bibitem{TreEmbBook05}
{\sc L.~N. Trefethen and M.~Embree}, {\em Spectra and Pseudospectra. {T}he
  Behavior of Nonnormal Matrices and Operators}, Princeton University Press,
  Princeton, NJ, 2005.

\bibitem{SluVor86}
{\sc A.~van~der Sluis and H.~A. van~der Vorst}, {\em The rate of convergence of
  conjugate gradients}, Numer. Math., 48 (1986), pp.~543--560.

\bibitem{Van92}
{\sc H.~A. van~der Vorst}, {\em Bi-{CGSTAB}: a fast and smoothly converging
  variant of {B}i-{CG} for the solution of nonsymmetric linear systems}, SIAM
  J. Sci. Statist. Comput., 13 (1992), pp.~631--644.

\bibitem{VanBook03}
\leavevmode\vrule height 2pt depth -1.6pt width 23pt, {\em Iterative {K}rylov
  methods for large linear systems}, vol.~13 of Cambridge Monographs on Applied
  and Computational Mathematics, Cambridge University Press, Cambridge, 2003.

\bibitem{VorBook65}
{\sc Y.~V. Vorobyev}, {\em Methods of Moments in Applied Mathematics},
  Translated from the Russian by Bernard Seckler, Gordon and Breach Science
  Publishers, New York, 1965.

\bibitem{WesSon80}
{\sc P.~Wesseling and P.~Sonneveld}, {\em Numerical experiments with a multiple
  grid and a preconditioned {L}anczos type method}, in Approximation methods
  for {N}avier-{S}tokes problems ({P}roc. {S}ympos., {U}niv. {P}aderborn,
  {P}aderborn, 1979), vol.~771 of Lecture Notes in Math., Springer, Berlin,
  1980, pp.~543--562.

\bibitem{WilBook65}
{\sc J.~H. Wilkinson}, {\em The Algebraic Eigenvalue Problem}, Monographs on
  Numerical Analysis, The Clarendon Press Oxford University Press, New York,
  1988.
\newblock 1st ed. published 1965.

\bibitem{WilBook63}
\leavevmode\vrule height 2pt depth -1.6pt width 23pt, {\em Rounding Errors in
  Algebraic Processes}, Dover Publications Inc., New York, 1994.
\newblock Reprint of the 1963 original published by Prentice-Hall, Englewood
  Cliffs, NJ.

\end{thebibliography}

\end{document}